\documentclass[10pt,a4,twoside,rm]{article}
\usepackage[vcentermath]{youngtab}
\usepackage{young}
\usepackage{ytableau}
\usepackage{theorem}
\usepackage{amssymb}
\usepackage[widespace]{fourier}
\usepackage{fancyhdr}
\usepackage{float}
\usepackage[bb=esstix,bbscaled=.5,cal=cm,calscaled=.92]{mathalfa}
\usepackage{bbm}
\usepackage{graphicx}
\usepackage[ansinew,latin1]{inputenc}
\usepackage[active]{srcltx}
\usepackage{amsfonts}
\usepackage{nccmath}
\usepackage{titlesec}
\usepackage{pb-diagram}
\usepackage{color}
\usepackage{epsfig}
\usepackage[sans]{dsfont}
\usepackage{multicol}
\usepackage[left=4.5cm,top=4cm,right=4.5cm,bottom=2cm]{geometry}
\usepackage{amsmath}
\usepackage{yfonts}
\usepackage{pb-diagram}
\usepackage[refpage]{nomencl}
\usepackage{ifthen}
\usepackage{lipsum}
\usepackage{commath}

\makenomenclature

\renewcommand\thesection{\arabic{section}} 
\renewcommand\thesubsection{\arabic{section}.\arabic{subsection}} 
\titleformat{\section}[block]{\scshape\centering}{\thesection.}{1em}{} 
\titleformat{\subsection}{\bfseries}{\thesubsection.}{1em}{} 

\newcommand{\seq}{{\rm seq}_n}

\newcommand{\kk}{\mathcal{K}}
\newcommand{\A}{\mathcal{A}}

\newcommand{\B}{\mathcal{B}_{r,n}}
\newcommand{\BB}{\mathbb{B}}
\newcommand{\Comp}{{\mathcal Comp}_n}
\newcommand{\Par}{{\mathcal Par}_n}
\newcommand{\CC}{ \mathbb C }
\newcommand{\Y}{\mathcal{Y}_{r,n}}
\newcommand{\YY}{\mathcal{Y}_{r,n}(q)}

\newcommand{\End}{{\rm End}}
\newcommand{\spa}{{\rm Span}}

\newcommand\es{\mathbbm{s}}
\newcommand\et{\mathbbm{t}}
\newcommand\eu{\mathbbm{u}}
\newcommand\ev{\mathbbm{v}}

\newcommand\ea{\mathbbm{a}}
\newcommand\eb{\mathbbm{b}}

\newcommand\bs{\mathbf{s}}

\newcommand\bu{\mathbf{u}}
\newcommand\bv{\mathbf{v}}

\newcommand{\s}{\mathfrak{s}}
\newcommand{\U}{\mathfrak{u}}
\newcommand{\V}{\mathfrak{v}}
\newcommand{\T}{\mathfrak{t}}
\newcommand{\R}{\mathfrak{r}}

\newcommand{\bT}{\pmb{\mathfrak{t}}}
\newcommand{\Bs}{\pmb{\mathfrak{s}}}

\newcommand{\Bu}{\pmb{\mathfrak{u}}}
\newcommand{\Bv}{\pmb{\mathfrak{v}}}

\newcommand{\aaa}{{\mathfrak{a}}}
\newcommand{\bbb}{{\mathfrak{b}}}
\newcommand{\ccc}{{\mathfrak{c}}}

\newcommand{\MP}{{ {Par}}_{r,n}}

\newcommand{\MC}{{ {Comp}}_{r,n}}

\newcommand{\Si}{\mathfrak{S}}
\newcommand{\std}{{\rm Std}}
\newcommand{\Tab}{{\rm Tab}}

\newcommand{\h}{{h}}
\newcommand{\HH}{ \mathcal{H}}
\newcommand{\TT}{{\mathfrak T}}
\newcommand{\E}{ {\mathcal E}_n(q)}
\newcommand{\Ea}{ {\mathcal E}_n^{\alpha}(q)}

\newcommand{\EE}{ {\mathcal E}_n}
\newcommand{\Elambda}{E_{\blambda}}

\newcommand\bS{\Sigma}
\newcommand\blambda{{\boldsymbol\lambda}}

\newcommand\bnu{{\boldsymbol\nu}}
\newcommand\be{\mathbb{E}}
\newcommand\bmu{{\boldsymbol\mu}}
\newcommand\YYlambda{\mathcal{Y}_{{ \alpha }}(q)}

\newcommand\beLambda{\be_{\Lambda}}
\theorembodyfont{ } \theoremstyle{marginbreak}

\theoremstyle{plain}

\newtheorem{teo}{Theorem}
\newtheorem{coro}[teo]{Corollary}

\newtheorem{defi}[teo]{Definition}
\newtheorem{exa}[teo]{Example}
\newtheorem{lem}[teo]{Lemma}
\newtheorem{propos}[teo]{Proposition}

\newtheorem{obs}[teo]{Remark}

\newenvironment{demo}
{\textsc{Proof.}} {\quad \hfill $\Box$}
\newenvironment{demoemb}
{\textsc{Proof of Theorem \ref{embeddding}.}} {\quad \hfill $\Box$}

\begin{document}

\title{\bf \normalsize Cell structures for the Yokonuma-Hecke algebra and the algebra of braids and ties }
\author{\sc  jorge espinoza{\thanks{Supported in part by Beca Doctorado Nacional 2013-CONICYT 21130109}} \, and steen ryom-hansen\thanks{Supported in part by FONDECYT grants 1121129 and 1171379  } }
\date{}   \maketitle
\begin{abstract}
\noindent \textsc{Abstract. }
We construct a faithful tensor representation for the Yokonuma-Hecke algebra $ \Y$, and use it to give
a concrete isomorphism between $ \Y $ and Shoji's modified Ariki-Koike algebra.
We give a cellular basis for $ \Y $ and show that the Jucys-Murphy elements for $ \Y $ are JM-elements in
the abstract sense. Finally, we construct a cellular basis for the Aicardi-Juyumaya algebra of braids
and ties.
\end{abstract}


\medskip
\noindent
MCS2010: 33D80.

\section{Introduction}
In the present paper, we study the representation theory of
the Yokonuma-Hecke algebra $ \Y $ in type $A$ and of the related Aicardi-Juyumaya algebra $ \EE $ of braids and ties.
In the past few years, quite a few papers have been dedicated to the study of both algebras.

\medskip
The Yokonuma-Hecke algebra $ \Y $
was first introduced in the sixties by Yokonuma \cite{Y} for general types 
but the recent activity on $ \Y $ was initiated by Juyumaya who in \cite{JJ3} gave a new presentation of $\Y$.
It is a deformation of the group algebra of the wreath product $ C_r \wr \Si_n $ of the cyclic group $ C_r $ and
the symmetric group $ \Si_n $. On the other hand, it is quite different from the more familiar deformation of $ C_r \wr \Si_n $, the Ariki-Koike algebra
$ \widetilde{\mathcal{H}}_{r,n} $. For example, the usual Iwahori-Hecke algebra $ {\mathcal H}_n$ of type $A$
appears canonically as a quotient of $ \Y$, whereas it appears canonically as subalgebra of $ \widetilde{\mathcal{H}}_{r,n} $.

\medskip
Much of the impetus to the recent development on $ \Y $ comes from knot theory. In
the papers \cite{CJKL}, \cite{CL}, \cite{JJ2} and \cite{JL}
a Markov trace on $ \Y $ and its associated knot invariant $ \Theta $ is studied.

\medskip
The Aicardi-Juyumaya algebra $ \EE$
of braids and ties, along with its diagram calculus, was introduced in \cite{AJ1} and \cite{JJ1} via a presentation derived from the
presentation of $ \Y$. The algebra $ \EE$ is also related to knot theory. Indeed,
Aicardi and Juyumaya constructed in \cite{AJ2} a Markov trace on $ \EE$,
which gave rise to a
three parameter knot invariant $ \Delta$. There seems to be no simple relation between $ \Theta $ and $ \Delta$.

\medskip
A main aim of our paper is to show that $ \Y$ and $ \EE$ are cellular algebras in the sense of Graham and Lehrer, \cite{GL}.
On the way we give a concrete isomorphism between $ \Y $ and Shoji's modified Ariki-Koike algebra $ {\mathcal{H}}_{r,n} $.
This gives a new proof of a result of Lusztig \cite{L2} and Jacon-Poulain d'Andecy \cite{JdA}, showing that
$ \Y $ is in fact a sum of matrix algebra over Iwahori-Hecke algebras of type $ A$.

\medskip
For the
parameter $ q = 1$, it was shown in Banjo's work  \cite{Ba} that the algebra $ \EE  $ is a special case of P. Martin's ramified partition algebras. Moreover,
Marin showed in \cite{M} that $ \EE $ for $ q=1 $ is isomorphic to a sum of matrix algebras over a certain wreath product algebra,
in the spirit of Lusztig's and Jacon-Poulain d'Andecy's Theorem.
He raised the question whether this result could be proved for general parameters.
As an application of our cellular basis for $ \EE$ we do obtain such a structure Theorem for $ \EE$, thus answering in the positive Marin's question.

\medskip
Recently it was shown in \cite{CJKL} and \cite{PDW} that the Yokonuma-Hecke algebra invariant $ \Theta$
can be described via a formula involving the HOMFLYPT-polynomial and the linking number.
In particular, when applied to
classical knots, $ \Theta$ and the HOMFLYPT-polynomial coincide (this was already known for some time).
Given our results on $ \EE $ it would be interesting to investigate whether a similar result would hold for $ \Delta$.

\medskip
Roughly our paper can be divided into three parts. The first part, sections 2 and 3,
contains the construction of a faithful tensor space module $ V^{\otimes n} $ for $ \Y$. The construction of $ V^{ \otimes n}$
is a generalization of the $ \EE $-module structure on $ V^{ \otimes n}$ that was defined in \cite{Rh} and it allows us
to conclude that $ \EE $ is a subalgebra of $ \Y $ for $ r \ge n $, and for {\it any} specialization of the ground ring.
The tensor space module $ V^{ \otimes n} $ is also related to the strange Ariki-Terasoma-Yamada action, \cite{ATY} and \cite{SaS},
of the Ariki-Koike algebra on $  V^{ \otimes n}$,
and thereby to the action
of Shoji's modified Ariki-Koike algebra $ \mathcal{H}_{r,n} $ on $ V^{ \otimes n}$, \cite{Sho}.
A speculating remark concerning this last point was made in \cite{Rh}, but the appearance of Vandermonde determinants in the
proof of the faithfulness of the action of $ \Y $ in $ V^{ \otimes n}$
makes the remark much more precise. The defining relations of the modified Ariki-Koike algebra also involve Vandermonde
determinants and from this we obtain the proof of the isomorphism $ \Y \cong \mathcal{H}_{r,n} $ by viewing both algebras as subalgebras of
${ \rm End} (V^{ \otimes n}) $. Via this, we get a new proof of Lusztig's and Jacon-Poulain d'Andecy's isomorphism Theorem for $ \Y$,
since it is in fact equivalent to a similar isomorphism Theorem for $ \mathcal{H}_{r,n} $, obtained independently by Sawada-Shoji and
Hu-Stoll.

\medskip

\medskip
The second part of our paper, section\textcolor{black}{s} 4 and 5,
contains the proof that $ \Y $ is a cellular algebra in the sense of
Graham-Lehrer, via a concrete combinatorial construction of a
cellular basis for it, generalizing Murphy's standard basis for the
Iwahori-Hecke algebra of type $A$. The fact that  $ \Y $ is cellular
could also have been deduced from the isomorphism $ \Y \cong
\mathcal{H}_{r,n} $ and from the fact that $ \mathcal{H}_{r,n} $ is
cellular, as was shown by Sawada and Shoji in \cite{SS}. Still, the
usefulness of cellularity depends to a high degree on having a
concrete cellular basis in which to perform calculations, rather
than knowing the mere existence of such a basis, and our
construction should be seen in this light.

\medskip
Cellularity is a particularly strong language for the study of modular, that is non-semisimple representation theory, which occurs in our situation when
the parameter $ q $ is specialized to a root of unity.
But here our applications go in a different direction and depend on a nice compatibility property of our cellular basis with respect to
a natural subalgebra of $ \Y$.
We get from this that
the elements $ m_{\s \s } $ of the cellular basis for $ \Y $, given by one-column standard multitableaux $ \s$, correspond to certain
idempotents that appear in Lusztig's presentation of $ \Y$ in \cite{L1} and \cite{L2}.
Using the faithfulness of the tensor space module $ V^{\otimes n} $ for $ \Y $ we get via this Lusztig's idempotent presentation of $ \Y$.
Thus the second part of the paper depends logically on the first part.

\medskip
In section 5 we treat the Jucys-Murphy's elements for $ \Y$. They were already introduced and studied by
Chlouveraki and Poulain d'Andecy in \cite{MCH}, but here we show that they are JM-elements in the abstract sense defined by Mathas, with
respect to the cell structure that we found.

\medskip
The third part of our paper, section 6,  contains the construction of a cellular basis for $ \EE $. This construction does not depend logically on
the results of parts 1 and 2, but is still strongly motivated by them. The
generic representation theory of $ \EE $ was already studied in \cite{Rh} and was shown to be
a blend of the symmetric group and the Hecke algebra representation theories and this is reflected in the cellular basis.
The cellular basis is also here a variation of Murphy's standard basis but the details of the construction are
substantially more involved than in the $ \Y$-case.

\medskip
As an application of our cellular basis we show that $ \EE$ is isomorphic to a direct sum of matrix algebras over certain wreath product algebras
$ \mathcal{H}^{wr}_{\alpha} \! \! $, depending on a partition $ \alpha$.
An essential ingredient in the proof of this result is a compatibility property of our cellular basis for $ \EE$ with respect to these subalgebras.
It appears to be a key feature of Murphy's standard basis and its generalizations that they carry compatibility properties of this kind,
see for example \cite{MHar}, \cite{EG} and \cite{BdVE}, and thus our work can be viewed as a manifestation of this phenomenon.

\medskip
We thank the organizers of the Representation Theory Programme at the Institut Mittag-Leffler for
the possibility to present this work in April 2015.  We thank G. Lusztig for pointing out on that occasion that
the isomorphism Theorem for $ \Y $ is proved already in \cite{L2}.
We thank J. Juyumaya for pointing out a missing condition in a first version of Theorem \ref{embeddding}.
{\color{black}We also} thank M. Chlouveraki, L. Poulain d'Andecy {\color{black}and S. Rostam} for their comments and D. Plaza for many useful discussions on the topic of the paper.

\medskip
\textcolor{black}{Finally, but not least, it is a great pleasure to thank the anonymous referee for a very detailed list of criticisms and
suggestions that helped us greatly improve the text.}
\section{Notation and Basic Concepts}

In this section we set up the fundamental notation and introduce the objects we wish to investigate.

\medskip
Throughout the paper we fix the rings $R:=\mathbb{Z}[q,q^{-1},\xi,r^{-1}, \Delta^{-1}]$
and $S:=\mathbb{Z}[q,q^{-1}]$,
where
$q$ is an indeterminate, $r$ is a positive integer, $\xi:=e^{ 2 \pi i/r} \in \mathbb C $
and $ \Delta $ is the Vandermonde determinant $ \Delta:=
 \prod_{0 \leq i < j \leq
r-1} (\xi^i -\xi^j)$.

We shall need the quantum integers $ [m]_q$ defined for $m \in  \mathbb Z$ by
$ [m]_q :=\dfrac{q^{2m}-1}{q^2-1}$
if $q \neq 1 $ and $[m]_q := m $ if $q = 1$.


\medskip
Let $\Si_n$ \nomenclature[01]{$\Si_n$}{The symmetric group on $n$
letters} be the symmetric group on $n$ letters.  We choose the
convention that it acts on $\mathbf{n}:= \{ 1,2\ldots, n \} $ on the
right. Let $\bS_n:=\{s_1,\ldots, s_{n-1}\} $ be the set of simple
transpositions in $ \Si_n $, that is $ s_i = (i,i+1) $. Thus,
$\Si_n$ is the Coxeter group on $ \bS_n$ subject to the relations
\begin{alignat}{3}
s_is_j&=s_js_i &&\quad\mbox{ for }\; |i-j|>1\label{s1}\\
s_is_{i+1}s_i&=s_{i+1}s_is_{i+1}&&\quad\mbox{ for } i=1,2,\ldots,n-2 \label{ss11}\\
s_i^2&=1&&\quad\mbox{ for } i=1,2,\ldots,n-1.\label{sss111}
\end{alignat}
We let $ \ell(\cdot) $ denote the usual length function on $\Si_n$.
\begin{defi}
Let $n$ be a positive integer. The Yokonuma-Hecke algebra, denoted
$\Y=\Y(q)$\nomenclature[02]{$\Y$}{The Yokonuma-Hecke algebra}, is
the associative $R$-algebra generated by the elements
$g_1,\ldots,g_{n-1},t_1,\ldots,t_n$, subject to the following
relations:
\begin{alignat}{3}
t_i^r&=1 &&\quad\mbox{ for all }\; i\label{r1}\\
t_it_j&=t_jt_i&&\quad\mbox{ for all }\; i,j\label{r2}\\
t_j g_i &=g_i t_{j s_i }  &&\quad\mbox{ for all } i,j \label{r3}\\
g_ig_j&=g_jg_i&&\quad\mbox{ for }\; |i-j|>1\label{r4}\\
g_ig_{i+1}g_i&=g_{i+1}g_ig_{i+1}&&\quad\mbox{ for all }\;
i=1,\ldots,n-2\label{r5}
\end{alignat}
together
with the quadratic relation
\begin{align}
g_i^2&=1+(q-q^{-1})e_ig_i\qquad \mbox{ for all }\; i\label{r6}
\end{align}
where
\begin{align}e_i:=\dfrac{1}{r}\sum_{s=0}^{r-1}t_i^st_{i+1}^{-s}.
\end{align}
\end{defi}
Note that since $ r $ is invertible in $ R $, the element $ e_i \in  \Y(q) $ makes sense.

\medskip
One checks that $ e_i $ is an idempotent and that
$ g_i $ is invertible in $ \Y(q) $ with inverse
\begin{equation}{\label{inverse}}
g^{-1}_i=g_i+(q^{-1}-q)e_i.
\end{equation}

The study of the representation theory $\Y(q) $ is one of the main
themes of the present paper.  $ \Y(q) $ can be considered as a
generalization of the usual Iwahori-Hecke algebra ${\mathcal H}_n
=\HH_n(q)$ \nomenclature[03]{$\HH_n$}{The Iwahori-Hecke algebra} of
type $A_{n-1} $ since $ {\mathcal Y}_{1,n}(q)  = \HH_n(q)$. In
general $ \HH_n(q)$ is a canonical quotient of $\Y(q) $ via the
ideal generated by all the $ t_i -1$'s. On the other hand, as a
consequence of the results of the present paper, $ \HH_n(q)$ also
appears as a subalgebra of $\Y(q) $ although not canonically.

\medskip
$\Y(q) $ was introduced by Yokonuma in the sixties as the endomorphism algebra of a
module for the Chevalley group of type $ A_{n-1} $, generalizing
the usual Iwahori-Hecke algebra construction, see \cite{Y}.
This also gave rise to a presentation for $ \Y(q) $. A different presentation for $ \Y(q)$, widely used in the literature, was found
by Juyumaya. The presentation given above
appeared
first in
\cite{MCH} and differs slightly from Juyumaya's presentation. In Juyumaya's presentation another variable $u $ is used
and the quadratic relation (\ref{r6}) takes the form $ \tilde{g}_i^2 =1+(u-1)e_i(\tilde{g}_i  +1)$. The relationship
between the two presentations is given by $ u = q^2 $ and
\begin{equation} \tilde{g}_i =g_i + (q-1)e_i g_i,
\end{equation}
or
equivalently $g_i =\tilde{g}_i + (q^{-1}-1)e_i \tilde{g}_i$,
see eg. \cite{CJKL}.

\medskip
In this paper we shall be interested in the general, not necessarily semisimple, representation theory of $ \Y(q) $ and shall therefore need base change of
the ground ring.
Let $ \kk $ be a commutative ring, with elements $ q, \xi \in \kk^{\times}$. Suppose moreover that $ \xi $ is an $r$'th root of unity
and that $ r $ and  $\prod_{0 \leq i < j \leq
r-1} (\xi^i -\xi^j)$ are
invertible in $ \kk $ (for example $ \kk $ a field with $ r, \xi \in \kk^{\times} $ and $ \xi $ a primitive $r$'th root of unity).
Then
we can make $ \kk $ into an $ R $-algebra by mapping $ q \in R $ to $ q \in \kk $, and $ \xi \in R $ to $ \xi \in \kk $.
This gives rise to the specialized Yokonuma-Hecke algebra
$$ \Y^{\kk}(q) =\Y(q) \otimes_R \kk. $$
\nomenclature[04]{$\Y^{\kk}$}{The specialized Yokonuma-Hecke
algebra}

\medskip
Let $w \in \Si_n$ and
suppose that $w=s_{i_1}s_{i_2}\cdots s_{i_m}$ is a reduced
expression for $w$. Then by the relations the element
$g_w:=g_{i_1}g_{i_2}\cdots g_{i_m}$ does not depend on the choice of
the reduced expression for $w$. We use the convention that $ g_1 := 1 $.
In \cite{JJ2} Juyumaya proved that the following set
is an $R$-basis for $\Y(q)$
\begin{align}{\label{juyu}}
\mathcal{B}_{r,n}=\{t_1^{k_1}t_2^{k_2}\cdots t_n^{k_n}g_w\mid w\in
\Si_n,\;k_1,\ldots,k_n\in \mathbb{Z}/r\mathbb{Z}\}.
\end{align}
In particular, $\Y(q)$ is a free $R$-module of
rank $r^nn!$. Similarly, $\Y^{\kk}(q)$ is a free over $\kk$ of rank
$r^nn!$.

\medskip
Let us introduce some useful elements of $ \Y(q) $ (or $\Y^{\kk}(q)$).
For $1\leq i,j\leq n$ we define $ e_{ij} $ by
\begin{align}
e_{ij}:=\dfrac{1}{r}\sum_{s=0}^{r-1}t_i^st_j^{-s}.
\end{align}
These $e_{ij}$'s are idempotents and $e_{ii}=1$ and $e_{i,i+1}=e_i$.
Moreover $e_{ij}=e_{ji}$ and it is easy to verify from (\ref{r3})
that
\begin{equation}{\label{defeij}}
e_{ij}=g_ig_{i+1}\cdots g_{j-2} e_{j-1} g_{j-2}^{-1}\cdots
g_{i+1}^{-1}g_{i}^{-1}\qquad \mbox{ for } i<j.
\end{equation}
From
(\ref{r1})-(\ref{r3}) one obtains that
\begin{alignat}{3}
t_i e_{ij}&=t_je_{ij}&&\mbox{ for all } i,j \label{eg2} \\
e_{ij}e_{kl}&=e_{kl}e_{ij}&&\mbox{ for all } i,j,k,l\\
 e_{ij} g_k &=g_k e_{i s_k,  j s_k } &&\mbox{ for all }  i,j \mbox{ and } k =1,\ldots,n-1 \label{eg1}.
\end{alignat}
For any nonempty subset $I\subset \mathbf{n}$ we extend the
definition of $e_{ij}$ to \textcolor{black}{$E_I$} by setting
\begin{align}
E_I:=\prod_{i,j\in I,\,i<j}e_{ij}
\end{align}
where we use the convention that $E_I:=1$ if $|I|=1$.

\medskip
We need a further generalization of this. Recall that a set of
subsets $A=\{I_1,I_2,\ldots I_k\}$ of $\mathbf{n}$ is called a
\textit{set partition} of $\mathbf{n}$ if the $I_j$'s are nonempty,
disjoint and have union $\mathbf{n}$. We refer to the $ I_i$'s as
the \textit{blocks} of $A$. The set of all set partitions of
$\mathbf{n}$ is denoted
$\mathcal{SP}_n$\nomenclature[05]{$\mathcal{SP}_n$}{The set of set
partition of $\mathbf{n}$}. There is a natural poset structure on $
{\cal SP}_n$ defined as follows. Suppose that $A=\{I_1,I_2,\ldots,
I_k\} \in {\cal SP}_n $ and $ B=\{J_1,I_2,\ldots, J_l \} \in {\cal
SP}_n$. Then we say that $ A \subseteq B $ if each $ J_j $ is a
union of some of the $ I_i$'s.

For any set partition $A=\{I_1,I_2,\ldots ,I_k\} \in
\mathcal{SP}_n$ we define
\begin{align}{\label{YokonumaIdem}}
E_A:=\prod_j E_{I_j}.
\end{align}
Extending the right action of $ \Si_n $ on $\mathbf{n}$ to
a right action on $\mathcal{SP}_n$ via
$Aw:=\{ I_1 w ,\ldots,   I_k w\}\in
\mathcal{SP}_n$ for $w\in\Si_n$, we have the following Lemma.
\begin{lem}\label{pro1}
For $A\in\mathcal{SP}_n$ and $w\in\Si_n$ as above, we have that
$$E_A g_w=g_w E_{Aw}.$$
In particular, if $w$ leaves invariant every block of $A$, or more generally
permutes certain of the blocks of $ A $ (of the same size), then $E_A$
and $g_w$ commute.
\end{lem}

\begin{demo}
This is immediate from (\ref{eg1}) and the definitions.
\end{demo}

\medskip
As mentioned above, the specialized Yokonuma-Hecke algebra $\Y^{\kk}(q) $ only exists
if $ r $ is a unit in $ {\kk}$.
The algebra of braids and ties $ \E$, introduced by  Aicardi and Juyumaya, is an algebra related to $ \Y(q) $ that exists for any ground ring. It has
a diagram calculus consisting of braids that may be decorated with socalled ties, which explains its name, see \cite{AJ1}.
Here we only give its definition in terms of generators and relations.

\begin{defi}\label{braidsties}
Let $n$ be a positive integer.  The algebra of braids and ties, $
{\cal E}_n = \E $\nomenclature[06]{$\mathcal{E}_n$}{The braids and
ties algebra} , is the associative $S$-algebra generated by the
elements $g_1,\ldots,g_{n-1},e_1,\ldots,e_{n-1}$, subject to the
following relations:
\begin{alignat}{3}
g_ig_j&=g_jg_i&&\quad\mbox{ for } |i-j|>1\label{E1}\\
g_ie_i&=e_ig_i&& \quad \mbox{ for all } i \label{E2}\\
g_ig_{j}g_i&=g_{j}g_ig_{j}&&\quad\mbox{ for }  |i-j|=1 \label{E3}\\
e_ig_{j}g_i&=g_{j}g_ie_{j}&&\quad\mbox{ for }  |i-j|=1 \label{E4}\\
e_ie_{j}g_j&=e_{i}g_je_{i}=g_{j}e_ie_{j}   &&\quad\mbox{ for }  |i-j|=1 \label{E5}\\
e_ie_j&=e_je_i&&\quad\mbox{ for all } i,j\label{E6}\\
g_ie_j&=e_jg_i&&\quad\mbox{ for } |i-j|>1\label{E7}\\
e_i^2&=e_i&&\quad\mbox{ for all } i\label{E8} \\
g_i^2&=1+(q-q^{-1})e_ig_i  && \, \, \, \, \, \,  \mbox{ for all }\; i\label{E9}.
\end{alignat}
\end{defi}

Once again,
this differs slightly from the presentation normally used for $ \E $, for example in \cite{Rh},
where the variable $u $ is used and the quadratic relation takes the form $ \tilde{g}_i^2 =1+(u-1)e_i(\tilde{g}_i  +1)$.
And once again, to change between the two presentations one uses $ u = q^2 $ and
\begin{equation}{\label{changeofpresentation}}g_i =\tilde{g}_i + (q^{-1}-1)e_i \tilde{g}_i
\end{equation}

\medskip
For any commutative ring $ \kk $ containing the invertible element $
q  $, we define the specialized algebra ${ \mathcal E}^{\kk}_n(q) $
via ${ \mathcal E}^{\kk}_n(q)  := \E \otimes_S \kk
$\nomenclature[07]{$\mathcal{E}^{\kk}_n$}{The specialized braids and
ties algebra} where $ \kk $ is made into an $S$-algebra by mapping $
q \in S $ to $ q \in \kk $.

\begin{lem}\label{pro11}
Let $ \kk $ be a commutative ring containing invertible elements $ r, \xi, \Delta   $ as above. Then there is a homomorphim
$ \varphi = \textcolor{black}{\varphi_{\kk}} : { \mathcal E}^{\kk}_n(q)  \longrightarrow \Y^{\kk}(q) $
of $ \kk$-algebras induced by
$ \varphi(g_i )  := g_i $ and $ \varphi(e_i )  := e_i $.
\end{lem}

\begin{demo}
This is immediate from the relations. We shall later on show that $\varphi $ is an embedding if $ r \ge n $.
\end{demo}

\medskip
Let $\mathbb{N}^0 $ denote the nonnegative integers. We next recall
the combinatorics of Young diagrams and tableaux. A
\textit{composition} $\mu=(\mu_1,\mu_2,\ldots,\mu_l)$ of $n \in
\mathbb{N}^0$ is a finite sequence in $\mathbb{N}^0 $ with sum $
n$\textcolor{black}{.} The $ \mu_i$'s are called the parts of $ \mu$.
A \textit{partition} of $n$ is a composition whose parts are
non-increasing. We write $\mu\models n$ and $\lambda\vdash n$ if
$\mu$ is a composition of $n$ and $\lambda$ is a partition of $n$.
In these cases we set $|\mu|:=n$ and $|\lambda|:=n$
\nomenclature[08]{$\lvert\mu\rvert$}{Size of the composition $\mu$} and define
the length of $ \mu $ or $ \lambda $ as the number of parts of $ \mu
$ or $ \lambda$. \textcolor{black}{If
$\mu=(\mu_1,\mu_2,\ldots,\mu_l)$ is a composition of length $l$ we
define the opposite composition $ \mu^{op} $ as $ \mu^{op} :=
(\mu_l,\ldots,\mu_2, \mu_1)$.} We denote by $ \Comp $
\nomenclature[09]{$\Comp$}{The set of composition of $n$} the set of
compositions of $ n $ and by $ \Par $ \nomenclature[10]{$\Par$}{The
set of partitions of $n$} the set of partitions of $n$. The
\textit{Young diagram} of a composition $\mu$ is the subset
$$[\mu]=\{(i,j)\mid 1\leq j\leq \mu_i \mbox{ and } i\geq 1 \}$$
of $\mathbb{N}^0\times\mathbb{N}^0$. The elements of $[\mu]$ are called
the \textit{nodes} of $\mu$. We represent $ [\mu] $
as an array of boxes in the plane, identifying each
node with a box. For example, if $\mu=(3,2,4)$ then
$$[\mu]=\yng(3,2,4) \, \, \, .$$
For $\mu\models n$ we define a $\mu$-\textit{tableau} as
a bijection $\mathfrak{t}:[\mu]\to\mathbf{n}$.  We identify
$\mu$-\textit{tableaux} with labellings of the nodes of $ [\mu] $: for example, if $\mu=(1,3)$ then
${\footnotesize\young(1,234)}$
is a $\mu$-tableau. If $\mathfrak{t} $ is a $ \mu $-tableau we write $ Shape(\mathfrak{t}) :=   \mu $.

We say that a $\mu$-tableau $\mathfrak{t}$ is \textit{row standard}
if the entries in $\mathfrak{t}$ increase from left to right in each
row and we say that $\mathfrak{t}$ is \textit{standard} if
$\mathfrak{t}$ is row standard and the entries also increase from
top to bottom \textcolor{black}{in each column}. The set of standard
$\lambda$-tableau\textcolor{black}{x} is denoted $\std(\lambda)$ and we write
$d_{\lambda}:=|\std(\lambda)|$ for its cardinality. For example,
${\footnotesize\young(235,14)}$ is row standard and
${\footnotesize\young(134,25)}$ is standard. For a composition of $
\mu $ of $ n $ we denote by $\mathfrak{t}^{\mu}$ the standard
tableau in which the integers $1,2,\ldots,n$ are entered in
increasing order from left to right along the rows of $[\mu]$. For
example, if $\mu=(2,4)$
then $\mathfrak{t}^{\mu}={\footnotesize\young(12,3456)}$.

The symmetric group $\Si_n$ acts on the right on the set of
$\mu$-tableaux by permuting the entries inside a given tableau. \textcolor{black}{Let $\s$ be a row standard $\lambda$-tableau. We denote by
$d(\s)$ the unique element of $\Si_n$ such that
$\s=\T^{\lambda}d(\s)$.}
The \textit{Young subgroup} \textcolor{black}{$\Si_{\mu}$} associated with \textcolor{black}{a composition $\mu$}
is the row stabilizer of $\mathfrak{t}^{\mu}$.
Let $\mu=(\mu_1,\ldots,\mu_k)$ and $\nu =(\nu_1,\ldots,\nu_l)$ be compositions. We write $\mu\unrhd \nu$ if
for all $  i\geq 1 $ we have
$$\sum_{j=1}^i\mu_j\geq \sum_{j=1}^i\nu_j$$
where we add zero parts $ \mu_i :=0 $ and $ \nu_i := 0$ at the end of $ \mu $ and $ \nu $ so that
the sums are always defined.
This is the dominance order on compositions.
We extend it to row standard tableaux as
follows.
Given a row standard tableau $\T$ of some shape and an integer $m\leq n$, we let
$\T\downarrow m$ be the tableau obtained from $\T$ by deleting all
nodes with entries greater than $m$. Then, for a pair
of $\mu$-tableaux $\s$ and $\T$ we write $\s \unrhd \T$
if $ {Shape}(\s\downarrow m)\unrhd {Shape}(\T\downarrow m)$
for all $m=1,\ldots,n$. We write $\s\rhd \T$ if
$\s\unrhd \T$ and $\s\neq \T$. This defines the dominance order on \textcolor{black}{row standard} tableaux. It is only a partial order, for example
$$\young(13,25,4)\;\rhd\; \young(24,35,1)\quad \mbox{ and }\quad
\young(13,25,4)\;\rhd\; \young(45,13,2)$$ whereas
${\footnotesize\young(24,35,1)}$ and
${\footnotesize\young(45,13,2)}$ are incomparable.

\medskip
We have that
$\T^{\lambda}\unrhd \T$ for all row standard $\lambda$-tableau
$\T$.\\

An \textit{$r$-multicomposition}, or simply a
\textit{multicomposition}, of $n$ is an ordered $r$-tuple
$\blambda=(\lambda^{(1)},\lambda^{(2)},\ldots,\lambda^{(r)})$ of
compositions $\lambda^{(k)}$ such that
$\sum_{i=1}^{r}|\lambda^{(i)}|=n$. We call $\lambda^{(k)}$ the
$k$'th component of $\blambda$, note that it may be empty. An \textit{$r$-multipartition}, or
simply a \textit{multipartition}, is a multicomposition whose
components are partitions.
The nodes of a multicomposition are labelled by tuples
$(x,y,k)$ with $k$ giving the number of the component and $ (x,y) $
the node of that component. For the multicomposition $ \blambda $ the set of nodes is denoted $ [\blambda] $.
This is the
Young diagram for $ \blambda $ and is represented graphically as the $ r$-tuple of Young diagrams of the components. For example, the
Young diagram of $\blambda=((2,3),(3,1),(1,1,1))$ is
$$\left( \yng(2,3) \;,\; \yng(3,1) \;,\; \yng(1,1,1)\; \right).$$
We denote by $\MC $ \nomenclature[11]{$\MC$}{The set of
multicomposition of $n$} the set of $r$-multicompositions of $n$ and
by $ \MP $ \nomenclature[12]{$\MP$}{The set of multipartitions of $n$}
the set of $r$-multipartitions of $n$. Let $\blambda $ be a
multicomposition of $n$.  A $\blambda$-\textit{multitableau} is a
bijection $\bT:[\blambda]\to\mathbf{n}$ which may once again be
identified with a filling of $ [\blambda] $ using the numbers from $
\mathbf{n}$. The restriction of $\bT$ to $ \lambda^{(i)} $ is called
the $i$'th component of $\bT$ \textcolor{black}{and we write $
\bT=(\T^{(1)}, \T^{(2)}, \ldots, \T^{(r)}) $ where $ \T^{(i)} $ is
the $i$'th component of $ \bT$.} We say that $\bT$ is \textit{row
standard} if all its components are row standard, and
\textit{standard} if all its components are standard tableaux. If $
\bT $ is a $ \blambda$-multitableau we write $ Shape(\bT) =
\blambda$. The set of all standard $\blambda$-multitableaux is
denoted by $\std(\blambda)$. In the examples
\begin{equation}
\bT=\left(\,\young(123,45)\;,\;\young(6,7,89)\;\right)\qquad
\Bs=\left(\,\young(278,14)\;,\;\young(56)\;,\;\young(3,9)
\;\right)\label{t}
\end{equation}
$\bT$ is a standard multitableau whereas $\Bs$ is only
a row standard tableau. We denote by $\bT^{\blambda}$ the
$\blambda$-multitableau in which $1,2,\ldots, n$ appear in order along the rows of
the first component, then along the rows of the second component,
and so on. For example, in (\ref{t}) we have that $ \bT= \bT^{\blambda}$ for
$\blambda=((3,2),(1,1,2))$.
For each multicomposition $\blambda$ we define the Young subgroup
$\Si_{\blambda}$ as the row stabilizer of
$\bT^{\blambda}$.

Let $\Bs$ be a row standard $\blambda$-multitableau. We denote by
$d(\Bs)$ the unique element of $\Si_n$ such that
$\Bs=\bT^{\blambda}d(\Bs)$. The set formed
by these elements is a complete set of right coset representatives
of $\Si_{\blambda}$ in $\Si_n$. Moreover
$$\{d(\Bs)\mid \s \mbox{ is a row standard } \blambda \mbox{-multitableau}\}$$
is a \textit{distinguished set of right coset representatives}, that is
$\ell(wd(\Bs))=\ell(w)+\ell(d(\Bs))$
for $w\in\Si_{\blambda}$.

Let $\blambda$ be a multicomposition of $n$ and let $\bT$ be a
$\blambda$-multitableau. For $j=1,\ldots,n$ we write
$p_{\bT}(j):=k$ if $j$ appears in the $k$'th component
$\T^{(k)}$ of $\bT$. We call $p_{\bT}(j) $ the
\textit{position} of $j$ in $\bT$.
When $\bT=\bT^{\blambda}$, we write $p_{\blambda}(j)$ for
$p_{\bT^{\blambda}}(j)$ and
say that a $\blambda$-multitableau $\bT$ is of the
\textit{initial kind} if
$p_{\bT}(j)=p_{\blambda}(j)$ for all $j=1,\ldots,n$.

Let $\blambda=(\lambda^{(1)},\lambda^{(2)},\ldots,\lambda^{(r)})$
and $\bmu=(\mu^{(1)},\mu^{(2)},\ldots,\mu^{(r)})$ be
multicompositions of $n$. We write $\blambda \unrhd \bmu$ if
$\lambda^{(i)}\unrhd \mu^{(i)}$ for all $i=1,\ldots,n$, this is our
dominance order on $ \MC$. If $\Bs $ and $\bT$ are row standard
\textcolor{black}{multitableaux} and $m=1,\ldots,n$ we define $\Bs
\downarrow m $ and $ \bT \downarrow m$ as for usual tableaux and
write $\Bs \unrhd \bT$ if $Shape(\Bs\downarrow m) \unrhd
Shape(\bT\downarrow m) $ for all $ m$.

It should be noted that our dominance order $ \unrhd $ is different
from the dominance order on multicompositions and multitableaux that is
used in some parts of the literature, for example in \cite{DJM}. Let
us denote by $ \succeq $ the order used in \cite{DJM}. Then we have
that
$$\left(\;\young(12,34)\;,\young(5)\;\right) \succeq
\left(\;\young(13,2)\;,\young(45)\;\right)$$ whereas these
multitableaux are incomparable with respect to $ \unrhd $. On the
other hand, if $\Bs$ and $\bT$ are multitableaux of the same shape and
$p_{\Bs}(j)=p_{\bT}(j)$ for all $j$, then we have that
$\Bs\unrhd\bT$ if and only if $\Bs\succeq\bT$.

To each $r$-multicomposition
$\blambda=(\lambda^{(1)},\ldots,\lambda^{(r)})$ we associate a
composition $ \|\blambda\| $ of length \textcolor{black}{$r$} as
follows
\begin{equation}\label{thecompositionas}
 \norm{\blambda}  := ( | \lambda^{(1)}| ,\ldots,| \lambda^{(r)}| ).
\end{equation}
\color{black}\nomenclature[13]{$\norm{\blambda}$}{The composition associated with the multicomposition $\blambda$}Let $ \Si_{  \norm{\blambda} } $ be the associated Young subgroup. Then $ w \in \Si_{  \norm{\blambda} } $ iff
  $ \bT^{\blambda} w $ is of the initial kind. For any $\blambda$-multitableau $ \Bs $ there is a
decomposition of $ d(\Bs ) $ with respect to $ \Si_{  \norm{\blambda} } $, that is
\textcolor{black}{
\begin{equation}\label{decompositionintialkind} d(\Bs ) = d(\Bs_0) w_{\Bs},\; \mbox{ where }\; d(\Bs_0) \in \Si_{  \norm{\blambda} }\;\mbox{ and }\; l( d(\Bs ) ) = l(d(\Bs_0) ) + l( w_{\Bs}).\end{equation}}
We define in this situation $ \Bs_0 = \bT^{\blambda} d( \Bs_0) $; \textcolor{black}{it} is  of the initial kind.
Let $ \bT$ be another multitableau of shape $ \bmu$ and let $ d(\bT) = d(\bT_0) w_{\bT}$
be its decomposition.
Suppose that $ \norm{\blambda} = \norm{\bmu} $ and that
$ w_{\Bs} = w_{\bT} $. Then we have the following compatibility property with respect to the dominance order
\begin{equation}
\Bs \unrhd \bT \mbox{ if and only if } \Bs_0 \unrhd \bT_0.
\end{equation}
Let $ \alpha := \norm{\blambda} $.
Let $ y \in \Si_n$ be as short as possible such that $ \Bs y $ is of the initial kind
and set
$ \T := \T^{\alpha} d(\Bs) y   $. Then
$ d(\Bs_0) = d(\T )$
and $ w_{\Bs} = y^{-1} $. If $ y = s_{i_1} s_{i_2} \ldots s_{i_k} $ is reduced expression for $ y $
then for all $ j $ we have that $ i_j $ and $ i_j +1 $ occur in distinct components of $ \Bs s_{i_1} s_{i_2}\ldots s_{i_{j-1}} $
(with $ i_j +1 $ to the left of $ i_j $) as can be seen using the inversion description of the length function on
$ \Si_n$,
and a similar property holds for $ w_{\Bs}$.
\color{black}

\section{Tensorial representation of $\Y(q)$}

In this section we obtain our first results by constructing a tensor space module for the
Yokonuma-Hecke algebra which we show is faithful. From this we deduce that the Yokonuma-Hecke algebra is in fact
isomorphic to
a specialization of the `modified Ariki-Koike' algebra,
that was introduced by Shoji in \cite{Sho} and studied for example in \cite{SS}.
\begin{defi}{\label{tensoraction}}
Let $V$ be the free $R$-module with basis $\{v_{i}^{t}\mid 1\leq
i\leq n,\; 0\leq t\leq r-1\}$.
Then we define operators $\mathbf{T}\in \End_R(V)$ and
$\mathbf{G}\in \End_R(V^{\otimes 2})$ as follows:
\begin{equation}{\label{operatorT}}
(v_i^t)\mathbf{T}:=\xi^{t}v_i^t
\end{equation}
and
\begin{equation}{\label{operatorG}}
(v_i^t\otimes v_j^s)\mathbf{G}:=\left\{\begin{array}{ll}v_j^s\otimes v_i^t&\mbox{ if }\;t\neq s\\
qv_i^t\otimes v_j^s&\mbox{ if }\;t=s,\;i=j\\
v_j^s\otimes v_i^t&\mbox{ if }\;t=s,\;i>j\\
(q-q^{-1})v_i^t\otimes v_j^s+v_j^s\otimes v_i^t&\mbox{ if
}\;t=s,\;i<j.\end{array}\right.
\end{equation}

We extend them to operators $\mathbf{T}_i$ and $\mathbf{G}_i$ acting
in the tensor space $V^{\otimes n}$ by letting $\mathbf{T}$ act in
the $i$'th factor and $\mathbf{G}$ in the $i$'th and $i+1$'st
factors, respectively.
\end{defi}

Our goal is to prove that these operators define a faithful
representation of the Yokonuma-Hecke algebra. We first prove an
auxiliary Lemma.

\begin{lem}{\label{o1}}
Let $\mathbf{E}_i$ be defined by
$\mathbf{E}_i:= \dfrac{1}{r}\sum_{m=0}^{r-1} \mathbf{T}_i^m
\mathbf{T}_{i+1}^{-m}$. Consider the map
\[(v_i^t\otimes v_j^s)\mathbf{E}:=\left\{\begin{array}{ll}0&\mbox{ if }\;t\neq s\\
v_i^t\otimes v_j^s&\mbox{ if }\;t=s. \end{array}\right.\] Then
$\mathbf{E}_i$ acts in $V^{\otimes n}$ as $\mathbf{E}$ in the
factors $(i, i+1)$ and as the identity in the rest.
\end{lem}

\begin{demo}
We have that
$$(v_j^t\otimes v_k^t)\mathbf{T}_i\mathbf{T}_{i+1}^{-1}=\xi^{t}\xi^{-t}v_j^t\otimes v_k^t=v_j^t\otimes v_k^t.$$
Thus we get immediately that $(v_i^t\otimes
v_j^s)\mathbf{E}_i=v_i^t\otimes v_j^s$ if $s=t$. Now, if $s\neq t$
we have that
$$(v_j^t\otimes v_k^s)\mathbf{T}_i\mathbf{T}_{i+1}^{-1}=\xi^{t}\xi^{-s}v_j^t\otimes v_k^t=\xi^{t-s}v_j^t\otimes v_k^t.$$
Since $0\leq t,s\leq r-1$, we have that $\xi^{t-s}\neq 1$ which
implies that $$ \sum_{m=0}^{r-1} \xi^{m(t-s)} =
(\xi^{r(t-s)}-1)/(\xi^{(t-s)}-1)=0$$ and so it follows that
$(v_i^t\otimes v_j^s)\mathbf{E}=0$ if $s\neq t$.
\end{demo}

\begin{obs}\label{equivrep}
The operators $ \mathbf{G}_i $ and $ \mathbf{E}_i$ should be
compared with the operators introduced in \cite{Rh} in order to
obtain a representation of $ \E $ in $ V^{\otimes n} $. Let us
denote by $ \widetilde{\mathbf{G}}_i $ and $
\widetilde{\mathbf{E}}_i $ the operators defined in \cite{Rh}. Then
we have that $ {\mathbf{E}}_i =  \widetilde{\mathbf{E}}_i  $ and
$$\textcolor{black}{\mathbf{G}_i}= \widetilde{\mathbf{G}}_i +(q^{-1}-1) \widetilde{\mathbf{E}}_i  \widetilde{\mathbf{G}}_i $$
corresponding to the change of presentation given in
(\ref{changeofpresentation}).
\end{obs}

\begin{teo}\label{t1}
There is a representation $\rho$ of $\Y(q)$ in $V^{\otimes n}$
given by $t_i\to \mathbf{T}_i$ and $g_i\to \mathbf{G}_i$.
\end{teo}

\begin{demo} We must check that the operators $\mathbf{T}_i$
and $\mathbf{G}_i$ satisfy the relations
$(\ref{r1}),\ldots,(\ref{r6}) $ of the Yokonuma-Hecke algebra. Here
the relations (\ref{r1}) and (\ref{r2}) are trivially satisfied
since the $\mathbf{T}_i$'s commute. The relation (\ref{r4}) is also easy to verify since the operators $\mathbf{G}_i$ and $\mathbf{G}_j$ act as  $\mathbf{G}$ in two different consecutive factors if $|i-j|>1$.

\medskip
In order to prove the braid relations (\ref{r5}) we rely on the
fact, obtained in \cite{Rh} Theo\-rem 1, that the operators
$\widetilde{\mathbf{G}}_i$'s and $\widetilde{\mathbf{E}}_i$'s
satisfy the relations for ${\cal E}_n(q) $ \textcolor{black}{(with modified quadratic relation as indicated just below Definition \ref{braidsties}).
In particular, the braid relations $ \widetilde{\mathbf{G}}_i \widetilde{\mathbf{G}}_{i+1}\widetilde{\mathbf{G}}_i =
\widetilde{\mathbf{G}}_{i+1} \widetilde{\mathbf{G}}_{i}\widetilde{\mathbf{G}}_{i+1} $ hold and also
$ \widetilde{\mathbf{E}}_i \widetilde{\mathbf{G}}_{i+1} \widetilde{\mathbf{G}}_{i}\widetilde{\mathbf{G}}_{i+1}
\widetilde{\mathbf{E}}_i  =
\widetilde{\mathbf{E}}_{i+1} \widetilde{\mathbf{G}}_{i+1} \widetilde{\mathbf{G}}_{i}\widetilde{\mathbf{G}}_{i+1}
\widetilde{\mathbf{E}}_{i+1}$ holds, as one sees from Definition \ref{braidsties}.
Via Remark \ref{equivrep} we now get that
$$\begin{array}{l}\mathbf{G}_i\mathbf{G}_{i+1}\mathbf{G}_i=
(1+(q^{-1}-1)\widetilde{\mathbf{E}}_i)(\widetilde{\mathbf{G}}_i \widetilde{\mathbf{G}}_{i+1}\widetilde{\mathbf{G}}_i+(q^{-1}-1)\widetilde{\mathbf{G}}_i
 \widetilde{\mathbf{G}}_{i+1 }\widetilde{\mathbf{E}}_{i+1}\widetilde{\mathbf{G}}_i      )(1+(q^{-1}-1)\widetilde{\mathbf{E}}_{i})
\\ \hspace{1.44cm} =(1+(q^{-1}-1)\widetilde{\mathbf{E}}_i)(\widetilde{\mathbf{G}}_i
\widetilde{\mathbf{G}}_{i+1}\widetilde{\mathbf{G}}_i+(q^{-1}-1)\widetilde{\mathbf{G}}_i
 \widetilde{\mathbf{G}}_{i+1 }\widetilde{\mathbf{E}}_{i+1} \widetilde{\mathbf{E}}_{i} \widetilde{\mathbf{G}}_i      )(1+(q^{-1}-1) \widetilde{\mathbf{E}}_{i})
\\ \hspace{1.44cm} =
(1+(q^{-1}-1)\widetilde{\mathbf{E}}_{i+1})(\widetilde{\mathbf{G}}_{i+1} \widetilde{\mathbf{G}}_{i}\widetilde{\mathbf{G}}_{i+1}+(q^{-1}-1)\widetilde{\mathbf{G}}_{i+1}
 \widetilde{\mathbf{G}}_{i}\widetilde{\mathbf{G}}_{i+1} \widetilde{\mathbf{E}}_{i} \widetilde{\mathbf{E}}_{i+1}    )(1+(q^{-1}-1)\widetilde{\mathbf{E}}_{i+1}) \\
\hspace{1.44cm}=\mathbf{G}_{i+1}\mathbf{G}_{i}\mathbf{G}_{i+1}
\end{array}$$
}and (\ref{r5}) follows as claimed. In a similar way we get that the
$\mathbf{G}_i$'s satisfy the quadratic relation (\ref{r6}).

We are then only left with the relation (\ref{r3}).
We have here three cases to consider:
\begin{align}\mathbf{T}_i\mathbf{G}_j&=\mathbf{G}_j\mathbf{T}_i\qquad |i-j|>1\label{c1}\\
\mathbf{T}_i \mathbf{G}_i&=\mathbf{G}_i \mathbf{T}_{i+1}\label{c2}\\
\mathbf{T}_{i+1} \mathbf{G}_i&=\mathbf{G}_i \mathbf{T}_i\label{c3}.
\end{align}
The case (\ref{c1}) clearly holds since the operators $\mathbf{T}_i$
and $\mathbf{G}_j$ act in different factors of the tensor product
$v_{i_1}^{j_1}\otimes v_{i_2}^{j_2}\otimes \ldots \otimes
v_{i_n}^{j_n}$. In order to verify the other two cases
we may assume that $ i= 1 $ and $n=2$. It is enough to evaluate on vectors of the form
$v_{i_1}^{j_1}\otimes v_{i_2}^{j_2}\in V^{\otimes 2}$. For $j_1=j_2$
the actions of $\mathbf{T}_1$ and $\mathbf{T}_2$ are given as the
multiplication with the same scalar and so the relations (\ref{c2})
and (\ref{c3}) also hold.

Suppose then finally that $j_1\neq j_2$. We then have that
\[(v_{i_1}^{j_1}\otimes v_{i_2}^{j_2})\mathbf{T}_1\mathbf{G}_1=\xi^{j_1}v_{i_2}^{j_2}\otimes v_{i_1}^{j_1}
=(v_{i_1}^{j_1}\otimes v_{i_2}^{j_2})\mathbf{G}_1\mathbf{T}_2
\]
and
\[(v_{i_1}^{j_1}\otimes v_{i_2}^{j_2})\mathbf{T}_2\mathbf{G}_1=\xi^{j_2}v_{i_2}^{j_2}\otimes v_{i_1}^{j_1}
=(v_{i_1}^{j_1}\otimes v_{i_2}^{j_2})\mathbf{G}_1\mathbf{T}_1
\]
and the proof of the Theorem is finished.
\end{demo}

\begin{obs}{\label{specializedaction}}
Let $ \kk $ be an $ R $-algebra as in the previous section with
corresponding specialized Yokonuma-Hecke algebra
$ \Y^{\kk}(q) $. Then we obtain a specialized tensor product representation $ \rho^{\kk}:\Y^{\kk}(q)  \rightarrow
\End_{\kk}(V^{\otimes n}) $. Indeed, the above proof amounts only to checking relations, and so carries over to $\Y^{\kk}(q) $.
\end{obs}

\begin{teo}{\label{faithful}}
$\rho$ and $  \rho^{\kk} $ are faithful representations.
\end{teo}

\begin{demo}
We first consider the faithfulness of $\rho$.
Recall Juyumaya's $R$-basis for $\Y(q)$
\[\mathcal{B}_{r,n}=\{g_{\sigma}t_1^{j_1}\cdots t_n^{j_n}\mid \sigma\in \Si_n,\;j_i\in \mathbb{Z}/r\mathbb{Z}\}.\]
For $ \sigma = s_{i_1 } \ldots s_{i_m } \in \Si_n $ written in reduced form we define  $ \mathbf{G}_{\sigma} :=
\mathbf{G}_{i_1}  \ldots \mathbf{G}_{i_m} $. To prove that $\rho$ is faithful it is enough to show that
$$\rho(\mathcal{B}_{r,n})=\{\mathbf{G}_{\sigma}\mathbf{T}_1^{j_1}\cdots
\mathbf{T}_n^{j_n}\mid \sigma\in \Si_n,\;j_k\in
\mathbb{Z}/r\mathbb{Z}\}$$ is an $R$-linearly independent subset of
${\rm End}(V^{\otimes n})$. Suppose therefore that there exists a
nontrivial linear dependence
\begin{align}\mathop{\sum_{\sigma\in \Si_n}}_{j_i\in \mathbb{Z}/r\mathbb{Z}}
\lambda_{j_1,\ldots,j_n,\sigma}\mathbf{G}_{\sigma}\mathbf{T}_1^{j_1}\cdots
\mathbf{T}_n^{j_n}&=0\qquad \label{d1}\end{align} where not every
$\lambda_{j_1,\ldots,j_n,\sigma}\in R $ is zero.

We first observe that for arbitrary $ a_i $'s and $\sigma\in \Si_n$
the action of $\mathbf{G}_{\sigma}$ on the special tensor
$v_n^{a_n}\otimes \cdots \otimes v_1^{a_1}$, having the lower indices strictly decreasing,
is particularly simple. Indeed, since $ \sigma= s_{i_1} \ldots s_{i_m} $ is
a reduced expression for $ \sigma$
we have that the action of $\mathbf{G}_{\sigma} =
\mathbf{G}_{i_1}  \ldots \mathbf{G}_{i_m} $
in that case
always
involves the third case of ({\ref{operatorG}}) and thus is given by place permutation, in other words
\begin{equation}{\label{placepermutation}} (v_n^{a_n} \otimes \cdots \otimes v_1^{a_1}) \mathbf{G}_{\sigma}
= ( v_n^{a_n}\otimes \cdots \otimes v_1^{a_1}) \sigma =
 v_{i_n}^{a_{i_n}}\otimes \cdots \otimes v_{i_1}^{a_{i_1}}
\end{equation}
for some permutation $ i_n, \dots, i_1 $ of $ n,\ldots,1 $ uniquely given by $\sigma$.
Let $\TT_n$ be the $R$-subalgebra of $\End(V^{\otimes n})$ generated
by the $  \mathbf{T}_i $'s.
For fixed $ k_1,\ldots, k_n$ we now
define $$ V_{k_1,\ldots, k_n}:= \spa_R \{ v_{k_1}^{j_1}
\otimes \cdots \otimes v_{k_n}^{j_n}\,  |\,  j_k \in
{\mathbb Z}/r {\mathbb Z} \}. $$ Then $ V_{k_1, \ldots, k_n} $
is a $ \TT_n $-submodule of $ V^{\otimes n}$.
Given (\ref{placepermutation}), to prove that the linear dependence (\ref{d1}) does not exist, it
is now enough to show that $ V_{k_1, \ldots, k_n} $ is a faithful $ \TT_n $-module.

For $ j=0,1, \ldots, r-1$ we define $ w^j_k \in V$ via
$$ w^j_k := \sum_{i=0}^{r-1} \xi^{ij} v^i_{k}. $$
Then $ \{w^i_k\,  |\ i=0,1, \ldots, r-1, k=1,\ldots,n\} $ is also an $R$-basis for $
V$, since for fixed $ k$ the base change matrix between $ \{v^i_k\,  |\ i=0,1,
\ldots, r-1\} $ and $ \{w^j_k\,  |\ j=0,1, \ldots, r-1\} $ is given by
a Vandermonde matrix with determinant $ \prod_{0 \leq i < j \leq
r-1} (\xi^i -\xi^j)$ which is a unit in $ R$. But then also $ \{ w^{j_1 }_{k_1}
\otimes \ldots \otimes w^{j_n }_{k_n} |  j_i  \in
\mathbb{Z}/r\mathbb{Z} \} $ is an $ R $-basis for $ V_{k_1,
\ldots, k_n} $.
On the other hand, for all $ j $ we have that $ \textcolor{black}{w^j_k  {\mathbf{T}}=
w^{j+1}_k }$ where the indices are understood modulo $r$. Hence, given
the nontrivial linear combination in $ \TT_n$
$$ \sum_{j_i  \in \mathbb{Z}/r\mathbb{Z}}
\lambda_{j_1,\ldots,j_n } \mathbf{T}_1^{j_1}\cdots
\mathbf{T}_n^{j_n}  $$
we get by acting with it on $ w^0_{k_1} \otimes
\ldots \otimes w^0_{k_n} $ the following nonzero element
$$  \sum_{j_i  \in \mathbb{Z}/r\mathbb{Z}} \lambda_{j_1,\ldots,j_n }
w^{j_1  }_{k_1} \otimes  \ldots \otimes w^{j_n }_{k_n}. $$
This proves the Theorem in the case of $ \rho$. The case $ \rho^{\kk} $ is proved similarly, using
that $ \prod_{0 \leq i < j \leq
r-1} (\xi^i -\xi^j)$ is a unit in $ \kk $ as well.


\end{demo}



\subsection{The modified Ariki-Koike algebra.}
In this subsection we obtain our first main result, showing that the Yokonuma-Hecke algebra is
isomorphic to a variation of the Ariki-Koike algebra, called
the modified Ariki-Koike algebra $\mathcal{H}_{r,n}$.
It was introduced by Shoji.
Given the faithful tensor representation of the previous subsection, the proof of this isomorphism Theorem is actually almost trivial, but still we think that it
is a surprising result. Indeed, the quadratic relations involving the braid group generators look quite different
in the two algebras and as a matter of fact
the usual Hecke algebra of type $ A_{n-1} $ appears
naturally as a subalgebra of the (modified) Ariki-Koike algebra, but only as quotient
of the Yokonuma-Hecke algebra.

\medskip
Let us recall Shoji's definition of the modified Ariki-Koike algebra.
He defined it over the
ring
$R_1 :={\mathbb Z}[q,q^{-1}, u_1, \ldots, u_r, \Delta^{-1}]$,
where $ q, u_1, \ldots, u_r $ are indeterminates and $ \Delta:= \prod_{i>j}( u_i- u_j)$ is the Vandermonde determinant.
We here consider the modified Ariki-Koike algebra over the ring $ R $, corresponding to
a specialization of Shoji's algebra via the homomorphism $ \varphi: R_1 \rightarrow R $ given by
$ u_i \mapsto \xi^i $ and $ q \mapsto q $.

\medskip
Let $\mathbf{A}$ be the square matrix of degree $r$ whose $ij$-entry
is given by $\mathbf{A}_{ij} = \xi^{j(i-1)}$ for $1\leq i,j\leq r$, i.e. $\mathbf{A}$
is the usual Vandermonde matrix.
Then we can write the inverse of $\mathbf{A}$ as
$\mathbf{A}^{-1}=\Delta^{-1}\mathbf{B}$, where $ \Delta =  \prod_{i>j}( \xi^i- \xi^j) $ and $ \mathbf{B}=(h_{ij})$
is the adjoint matrix of $\mathbf{A}$, and for $1\leq i\leq r$ define a
polynomial $F_i(X)\in \mathbb{Z}[\xi][X] \subseteq R[X]$ by
$$F_i(X):=\sum_{1\leq i\leq r}h_{ij}X^{j-1}.$$
\begin{defi}
The modified Ariki-Koike algebra, denoted $\mathcal{H}_{r,n}=
\mathcal{H}_{r,n}(q)$\nomenclature[14]{$\mathcal{H}_{r,n}$}{The
modified Ariki-Koike algebra}, is the associative $R$-algebra
generated by the elements $\h_2,\ldots, \h_n$ and
$\omega_1,\ldots,\omega_n$ subject to the following relations:
\begin{alignat}{3}
(\h_i-q)(\h_i+q^{-1})&=0 &&\quad\mbox{ for all }\; i\label{z1}\\
\h_i\h_j&=\h_j\h_i&&\quad\mbox{ for }\; |i-j|>1\label{z2}\\
\h_i\h_{i+1}\h_i&=\h_{i+1}\h_i\h_{i+1}&&\quad\mbox{ for all }\;
i=1,\ldots,n-2\label{z3}\\
(\omega_i-\xi^1)\cdots(\omega_i-\xi^r)&=0&&\quad\mbox{ for all }\; i\label{z4}\\
\omega_i\omega_j&=\omega_j\omega_i&&\quad\mbox{ for all }\;
i,j\label{z5}
\end{alignat}
\begin{alignat}{3}
\h_j \omega_j &=\omega_{j-1}\h_j+\Delta^{-2}\sum_{c_1<c_2}(\xi^{c_2}-\xi^{c_1})(q-q^{-1})F_{c_1}(\omega_{j-1})F_{c_2}(\omega_j)&&\label{z6}\\
\h_j \omega_{j-1} &=\omega_{j}\h_j-\Delta^{-2}\sum_{c_1<c_2}(\xi^{c_2}-\xi^{c_1})(q-q^{-1})F_{c_1}(\omega_{j-1})F_{c_2}(\omega_j)&&\label{z7}\\
\h_j\omega_l&=\omega_l \h_j\qquad l\neq j,j-1&& \label{z8}
\end{alignat}
\end{defi}

$ \mathcal{H}_{r,n}(q) $ was introduced as a way of approximating the usual Ariki-Koike algebra and
is isomorphic to it if a certain separation condition holds.
In general the two algebras are not isomorphic, but related via a, somewhat mysterious, homomorphism from the Ariki-Koike algebra
to $\mathcal{H}_{r,n}(q) $, see \cite{Sho}.

\medskip

Sakamoto and Shoji, \cite{Sho} and \cite{SS}, gave a $ \mathcal{H}_{r,n}(q) $-module structure on $ V^{ \otimes n }  $ that we now explain.
We first introduce a total order on the $v_{i}^j$'s via
\begin{align}v_{1}^1,v_{2}^1,\ldots,v_{n}^1,v_{1}^2,\ldots,v_{n}^2,\ldots,v_{1}^r,\ldots,v_{n}^r\end{align}
and denote by $v_1,\ldots,v_{rn}$ these vectors in this order. We then define the linear operator $\mathbf{H}\in \End(V^{\otimes
2})$ as follows:
\[(v_i\otimes v_j)\mathbf{H}:=\left\{\begin{array}{ll}
qv_i\otimes v_j&\mbox{ if } i=j\\
v_j\otimes v_i&\mbox{ if } i>j\\
(q-q^{-1})v_i\otimes v_j+v_j\otimes v_i&\mbox{ if } i<j. \end{array}\right.\]
We then extend this to an operator $\mathbf{H}_i$ of
$V^{\otimes n}$ by letting $\mathbf{H}$ act in the $i$'th and
$i+1$'st factors.
This is essentially Jimbo's original operator for constructing tensor representations for the usual Iwahori-Hecke algebra $\mathcal{H}_n$
of type $ A$.
The following result is shown in \cite{Sho}.

\begin{teo}
The map $\tilde{\rho}:\mathcal{H}_{r,n}(q)\to \End(V^{\otimes n})$ given
by $\h_j\to \mathbf{H}_i$, $\omega_j\to \mathbf{T}_j$ defines a
faithful representation of $\mathcal{H}_{r,n}(q) $.
\end{teo}

We are now in position to prove the following main Theorem.

\begin{teo}{\label{mainthisection}}
The Yokonuma-Hecke algebra $\YY$ is isomorphic to the modified
Ariki-Koike algebra $\mathcal{H}_{r,n}(q)$.
\end{teo}

\begin{demo}
By the previous Theorem and Theorem \ref{faithful} we can identify
$\YY$ and $ \mathcal{H}_{r,n}(q) $ with the subalgebras $\rho(\YY)$ and
$\tilde{\rho}( \mathcal{H}_{r,n}(q) ) $ of $\End(V^{\otimes n})$, respectively.
Hence, in order to prove the Theorem we must show that $\rho(\YY) = \tilde{\rho}( \mathcal{H}_{r,n}(q)  )  $.
But by definition,  we surely have that the $\mathbf{T}_i$'s belong to
both subalgebras, since $ \mathbf{T}_i  = \rho(t_i) $ and $ \mathbf{T}_i  = \tilde{\rho}( \omega_i) $.

It is therefore enough to show that the
$\mathbf{G}_i$'s from $\rho(\Y(q))$ belong to
$\tilde{\rho}(\HH_{r,n})$, and that the
$\mathbf{H}_i$'s from $\tilde{\rho}(\HH_{r,n})$ belong to $\rho(\Y(q))$.

On the other hand, the operator $\mathbf{G}$ coincides with the
operator denoted by $S$ in Shoji's paper \cite{Sho}.  But then Lemma 3.5 of that paper is the equality
$$\mathbf{G}_{i-1}=\mathbf{H}_i-\Delta^{-2}(q-q^{-1})\sum_{c_1<c_2}F_{c_1}(\mathbf{T}_{i-1})F_{c_2}(\mathbf{T}_i).$$
Thus, since
$\Delta^{-2}(q-q^{-1})\sum_{c_1<c_2}F_{c_1}(\mathbf{T}_{i-1})F_{c_2}(\mathbf{T}_i)$
belongs to both algebras $\tilde{\rho}( \mathcal{H}_{r,n}(q)     )$ and $\rho( \YY)$, the Theorem follows.
\end{demo}

\medskip
Lusztig gave in \cite{L1} a structure Theorem for $ \YY $, showing that it is a direct sum of matrix algebras over Iwahori-Hecke algebras of type $ A $. This
result was recently recovered by Jacon and Poulain d'Andecy in \cite{JdA}.
We now briefly explain how this result, via our isomorphism Theorem, is equivalent to a similar result for
$\mathcal{H}_{r,n}(q)$, obtained in \cite{HS} and \cite{Sho}

\medskip
For a composition $ \mu=(\mu_1,\mu_2,\ldots,\mu_r)$ of $ n $ of
length $ r $, we let $ \mathcal{H}_{\mu}(q) $
\nomenclature[15]{$\mathcal{H}_{\mu}$}{The Young-Hecke algebra} be
the corresponding {\it Young-Hecke algebra}, by which we mean that $
\mathcal{H}_{\mu}(q) $ is the $R$-subalgebra of $\mathcal{H}_n(q)$
generated by the $ g_i $'s for $ i \in \bS_n \cap \Si_{\mu}$. Thus
$\mathcal{H}_{\mu}(q)  \cong \mathcal{H}_{\mu_1}(q) \otimes \ldots
\otimes \mathcal{H}_{\mu_r}(q) $ where each factor $ {\mathcal
H}_{\mu_i}(q) $ is a Iwahori-Hecke algebra corresponding to the
indices given by the part $ \mu_i $. Let $ p_{\mu} $ denote the
multinomial coefficient $$ p_{ \mu}:= \binom{n}{\mu_1  \cdots  \mu_r
}. $$

\medskip
With this notation, the structure Theorem due to Lusztig and Jacon-Poulain d'Andecy is as follows
\begin{equation}{\label{Jacond'Andecy}}
 \YY \cong \bigoplus_{\mu = (\mu_1,\mu_2,\ldots,\mu_r)  \models n} \mbox{Mat}_{p_{\mu}} ({\mathcal H}_{\mu}(q))
\end{equation}
where for any $ R$-algebra $ \mathcal A $, we denote by $
\mbox{Mat}_{m} ({\mathcal A}) $ the $ m \times m $ matrix algebra
with entries in $ \mathcal A $.

\medskip
On the other hand, a similar structure Theorem was established for the modified Ariki-Koike algebra $\mathcal{H}_{r,n}(q)$, independently by Sawada and Shoji in
\cite{SS} and by Hu and Stoll in \cite{HS}:
\begin{equation}{\label{HuStoll}}
{ \mathcal{H}_{r,n}(q)} \cong \bigoplus_{\mu = (\mu_1,\mu_2,\ldots,\mu_r)  \models n} \mbox{Mat}_{p_{\mu}} ({\mathcal H}_{\mu}(q)).
\end{equation}
Thus, our isomorphism Theorem {\ref{mainthisection}} shows that above two structure Theorems are equivalent.

\medskip
We finish this section by showing the following embedding Theorem, already announced above.
It is also a consequence of our tensor space module for $ \YY$.

\begin{teo}\label{embeddding}
Suppose that $ r \geq n$. Then the homomorphism
$ \varphi: { \mathcal E}^{\mathcal{K}}_n(q)  \rightarrow \Y^{\mathcal{K}}(q) $ introduced in Lemma \ref{pro11} is an embedding.
\end{teo}

In order to prove Theorem \ref{embeddding} we need to modify the proof of Corollary 4 of \cite{Rh} to make it valid for
general $ \kk$. For this we first prove the following Lemma.

\begin{lem}\label{previolinealdep}
Let $ \kk $ be an $R$-algebra as above and
let $A=(I_1,\ldots,I_d)\in\mathcal{SP}_n$ be a set partition. Denote by $V_{A} $ the $\kk$-submodule of $ V^{\otimes n} $ spanned by the vectors
$$v_{n}^{j_n}\otimes\cdots \otimes v_{k}^{j_k}\otimes\cdots  \otimes v_{l}^{j_l} \otimes \cdots \otimes v_{1}^{j_1}\qquad 0\leq j_i \leq r-1$$
with decreasing lower indices and satisfying that $j_k=j_l$ exactly
\textcolor{black}{if} $k$ and $l$ belong to the same block $I_i$
of $ A $. Let $ E_A \in \EE^{\kk}(q)$ be the element defined the
same way as $ E_A \in \YY $, that is via formula
({\ref{YokonumaIdem}}). Then for all $v\in V_A$ we have that
$vE_{A}= v $ whereas $vE_B=0$ for $ B \in \mathcal{SP}_n $
satisfying $B\not\subseteq A$ with respect to the order $ \subseteq
$ introduced above.
\end{lem}

\begin{demo}
In order to prove the first statement it is enough to show that $e_{kl}$ acts as the identity on the basis vectors of $ V_{A} $
whenever $ k $ and $l $ belong to the same block of $ A $.
But this follows from the expression for $e_{kl}$ given in (\ref{defeij}) together with
the definition (\ref{operatorG}) of the action of $\mathbf{G}_{i}$
on $ V^{\otimes n}$ and Lemma \ref{o1}.
Just as in the proof of Theorem \ref{faithful} we use that the action of $\mathbf{G}_{i}$ on $ v \in V_A$
is just permutation of the $i$'th and $i+1 $'st factors of $ v$ since the lower indices are decreasing.

\medskip
In order to show the second statement, we first remark that the
condition $B\not\subseteq A$ means that there exist
\textcolor{black}{$k$ and $l$} belonging to the same block of $B$, but
to different blocks of $A$. In other words \textcolor{black}{$e_{kl}$}
appears as a factor of the product defining $E_B$ whereas for all
basis vectors of $ V_A $
$$v_{n}^{j_n}\otimes\cdots \otimes v_{k}^{j_k}\otimes\cdots  \otimes v_{l}^{j_l} \otimes \cdots \otimes v_{1}^{j_1} $$
we have that $ j_k \neq j_l $. Just as above, using that the action
of $\mathbf{G}_{i}$ is given by place permutation when the lower
indices are decreasing, we deduce from this that $ V_A
\textcolor{black}{e_{kl}} = 0 $ and so finally that $ V_A E_B=0$, as
claimed.
\end{demo}

\medskip
\begin{demoemb}
\textcolor{black}{It is enough to show
that the composition $ \rho^{\kk}  \circ \varphi $ is injective since we know from Theorem {\ref{faithful}} that $ \rho^{\kk} $ is faithful.}
Now recall from Theorem 2 of \cite{Rh} that the set $\{{E}_{A} {g}_w | A \in \mathcal{SP}_n, w \in \Si_n \}$
generates $ \EE(q) $ over $ \CC[q,q^{-1}] $ (it is even a basis). The proof of this does not involve any special properties of $ \CC$ and hence
$\{{E}_{A} {g}_w | A \in \mathcal{SP}_n, w \in \Si_n \}$ also generates $ \EE^{\kk}(q) $ over $ \kk$.
\medskip

Let us now consider a
nonzero element
$  \sum_{w,A}r_{w,A}{E}_{A} {G}_w  $
in $ \EE^{\kk}(q) $. Under $ \rho^{\kk}  \circ \varphi $ it is
mapped to
$  \sum_{w,A}r_{w,A}\mathbf{E}_{A} \mathbf{G}_w $ which we must show to be nonzero.

\medskip
For this we choose
$A_0 \in\mathcal{SP}_n$ satisfying $r_{w,A_0}\neq 0$ for some $w\in\Si_n$ and
minimal with respect to this under our order $\subseteq$ on $\mathcal{SP}_n$.
Let $ v \in V_{A_0} \setminus \{ 0 \}$ where $ V_{A_0} $ is defined as in the previous Lemma \ref{previolinealdep}.
\textcolor{black}{Note that the condition $ r \ge n $ ensures that $ V_{A_0} \neq 0$, so such a $ v$ does exist.}
Then the Lemma gives us that
\begin{equation}\label{lowerindices}
 v \big(\sum_{w,A}r_{w,A}\mathbf{E}_{A} \mathbf{G}_w\big) = v \big(\sum_{w}r_{w,A_0} \mathbf{G}_w \big).
\end{equation}
The lower indices of $ v $ are strictly decreasing and so each $ \mathbf{G}_w  $ acts on it by place permutation. It follows from
this that (\ref{lowerindices}) is nonzero, and the Theorem is proved.
\end{demoemb}

\begin{obs}
The above proof did not use the linear independence of
$\{{E}_{A} {g}_w | A \in \mathcal{SP}_n ,  \newline w \in \Si_n \}$ over $ \kk$. In fact, it gives a
new proof of Corollary 4 of \cite{Rh}.

\textcolor{black}{
In the special case $ \kk = {\mathbb C}[q,q^{-1}] $ and $ r=n$ the Theorem
is an immediate consequence of the faithfulness of the tensor product
$ V^{\otimes n} $ as an $ \E$-module, as proved in Corollary 4 of \cite{Rh}.
Indeed, let $ \rho_{\EE}^{{\mathbb C}[q,q^{-1}]}:
\EE(q) \rightarrow {\rm End}(V^{\otimes n} ) $
be the homomorphism associated with the $ \E $-module structure on
$ V^{\otimes n} $, introduced in \cite{Rh}.
Then the injectivity of $ \rho_{\EE}^{{\mathbb C}[q,q^{-1}]}$ together with the
factorization $ \rho_{\EE}^{{\mathbb C}[q,q^{-1}]} = \rho^{{\mathbb C}[q,q^{-1}]}\circ \varphi_{{\mathbb C}[q,q^{-1}]}$
shows directly
that $ \varphi_{{\mathbb C}[q,q^{-1}]}$ is injective. One actually checks that the
proof of Corollary 4 of \cite{Rh} remains valid for $ \kk = R $ and
$ r \ge n $, but still this is not enough to prove injectivity of $ \varphi=
\varphi_{\kk} $ for a general $ \kk$ since extension of scalars from  $ R $
to $ \kk $ is not left exact. Note that the specialization argument
of \cite{Rh} fails for general $ \kk$.}
\end{obs}


\section{Cellular basis for the Yokonuma-Hecke algebra}
The goal of this section is to construct a cellular basis for the
Yokonuma-Hecke algebra. The cellularity of the Yokonuma-Hecke algebra could also have been
obtained from the cellularity of the modified Ariki-Koike algebra, see \cite{SS}, via our isomorphism Theorem from the previous section.
We have several reasons for still giving a direct construction of a cellular basis for the Yokonuma-Hecke algebra.
Firstly, we believe that our construction is simpler and more natural than the one in \cite{SS}.
Secondly, our basis turns out to have a nice compatiblity property with the subalgebra $\TT_n $ of $ \Y(q) $ studied above, a compatibility that
we would like to emphasize.
This compatibiliy is essential for our proof of Lusztig's presentation for $ \Y(q)$, given at the end of this section.
We also need the cellular basis in order to show, in the following section, that the Jucys-Murphy operators introduced by
Chlouveraki and Poulain d'Andecy are JM-elements in the abstract sense introduced by Mathas.
Finally, several of the methods for the construction of the basis are needed in the last section
where the algebra of braids and ties is treated.

\medskip
Let us start out by recalling the definition from \cite{GL} of a cellular basis.

\begin{defi}{\label{cellular}} Let $\mathcal{R}$ be an integral domain. Suppose that $A$ is an
$\mathcal{R}$-algebra which is free as an $\mathcal{R}$-module. Suppose that $(\Lambda,
\geq)$ is a poset and that for each $\lambda\in \Lambda$ there is a
finite indexing set $T(\lambda)$ (the '$\lambda$-tableaux') and elements $c_{\s\T}^{\lambda}\in
A$ such that
$$\mathcal{C}=\{c_{\s\T}^{\lambda}\mid \lambda\in \Lambda \mbox{ and } \s,\T\in T(\lambda)\}$$
is an $\mathcal{R}$-basis of $A$. The pair $(\mathcal{C},\Lambda)$ is a
\textit{cellular basis} of $A$ if
\begin{enumerate}\renewcommand{\labelenumi}{\textbf{(\roman{enumi})}}
\item The $\mathcal{R}$-linear map $*:A\to A$ determined by
$(\mathop{c_{\s\T}^{\lambda}})^*=c_{\T\s}^{\lambda}$ for all
$\lambda\in\textcolor{black}{\Lambda}$ and all $\s,\T\in T(\lambda)$
is an algebra anti-automorphism of $A$.

\item For any $\lambda\in \Lambda,\; \T\in T(\lambda)$ and $a\in A$
there exist $r_{\V}\in R$ such that for all $\s\in T(\lambda)$
$$c_{\s\T}^{\lambda} a \equiv \sum_{\V\in T(\lambda)} r_{\V}c_{\s\V}^{\lambda} \mod{A^{\lambda}}$$
where $A^\lambda$ is the $\mathcal{R}$-submodule of $A$ with basis
$\{c_{\U\V}^{\mu}\mid \mu\in \Lambda,\mu>\lambda \mbox{ and }
\U,\V\in T(\mu)\}$.
\end{enumerate}
If $A$ has a cellular basis we say that $A$ is a \textit{cellular
algebra}.
\end{defi}

For our cellular basis for $ \YY$ we use
for $ \Lambda $ the set $\MP $ of
$r$-multipartitions of $n$, endowed with the dominance order as explained in
section 2, and for $T(\lambda) $ we use the set of standard $r$-multitableaux
$\std(\blambda)$, introduced in the same section.
For $\ast:\Y(q)\to \Y(q)$ we use the $R$-linear antiautomorphism
of $\Y(q)$ determined by $g_i^{\ast}=g_i$ and $t_k^{\ast}=t_k$ for
$1\leq i <n$ and $1\leq k \leq n$.
Note that $ \ast $ does exist as can easily be checked from the relations
defining $ \Y(q) $.

We then only have to explain the construction of the
basis element itself, for pairs of standard tableaux.
Our guideline for this is Murphy's construction of {\it the standard basis} of
the Iwahori-Hecke algebra $ {\mathcal H}_n(q)$.

\medskip
For $
\blambda
\in \MC $ we first define
\begin{equation}\label{Murphy-element}
 x_{\blambda}:=\sum_{w\in\Si_{\blambda}}q^{\ell(w)}g_{w}\, \, \,  \textcolor{black}{ \in \YY.}
\end{equation}
In the case of the Iwahori-Hecke algebra $ {\mathcal H}_n(q) $, and
$\blambda$ a usual composition, the element $ x_{\blambda}$ is the
starting point of Murphy's standard basis, corresponding to the most
dominant tableau $\T^{\blambda}$.
In our $ \YY $ case, the element $ x_{\blambda} $
will only be the first ingredient of the \textcolor{black}{cellular basis} element corresponding
to the tableau $\bT^{\blambda}  $. Let us now explain the other two
ingredients.

\medskip
For a composition $ \mu =( \mu_1, \ldots , \mu_k ) $ we define the {\it reduced composition}
$ {\rm red \, } \mu $ as the composition obtained from $ \mu $ by deleting all zero parts $ \mu_i  =0 $
from $ \mu$. We say that a composition $ \mu $ is reduced if $ \mu = {\rm red \, } \mu$.

\medskip


For any reduced composition $ \mu =( \mu_1, \mu_2, \ldots , \mu_k ) $ we introduce the set partition $ A_{\mu} :=  (I_1, I_2, \ldots, I_k) $ by
filling in the numbers consecutively, that is
\begin{equation}\label{fillinginthenumbers}
 I_1 := \{ 1, 2, \ldots, \mu_1  \}, \,  I_2 := \{  \mu_1 +1,   \mu_1 +2, \ldots,  \mu_1 + \mu_2  \}, \, etc.
\end{equation}
and for a multicomposition $ \blambda \in \MC$ we define
$A_{\blambda} := A_{ {\rm red  }  \norm{\blambda} } \in
\mathcal{SP}_n$ \nomenclature[16]{$A_{\blambda}$}{The set partition
associated with $\blambda$}. Thus we get for any $ \blambda \in \MC$
an idempotent $E_{A_{\blambda}} \in \YY $ which will be the second
ingredient of our $ \YY $-element for $\bT^{\blambda}  $. We shall
from now on use the notation
\begin{equation}
\textcolor{black}{\Elambda} := E_{A_{\blambda}}.
\end{equation}
\nomenclature[17]{$\Elambda$}{The idempotent $E_{A_{\blambda}}$}Clearly $ t_i  \textcolor{black}{\Elambda} = \textcolor{black}{\Elambda} t_i $ for all $ i$. Moreover $\textcolor{black}{\Elambda}$ satisfies
the following key property.
\begin{lem}\label{o5}\mbox{}
Let $ \blambda \in \MC $ and let $A_{\blambda}$ be the associated set partition.
Suppose that $ k $ and $ l$ belong to the same block of $A_{\blambda}$. Then
$ t_k \textcolor{black}{\Elambda} = t_l \textcolor{black}{\Elambda} $.
\end{lem}

\begin{demo}
This follows from the definitions.
\end{demo}

\medskip
From Juyumaya's basis (\ref{juyu}) it follows that $t_i$ is
a diagonalizable element on $\Y(q)$. The eigenspace projector
for the action $ t_i$ on $\Y(q)$
with eigenvalue $ \xi^k$ is
\begin{equation}\label{eigen}
u_{ik} = \frac{1}{r} \sum_{j = 0}^{r-1}  \xi^{-jk} t_i^{j}\in  \YY
\end{equation}
that is
$  \{ v \in \YY | t_i v =  \xi^k v \} = u_{ik} \YY.$
\textcolor{black}{For $\blambda=(\lambda^{(1)},\ldots,\lambda^{(r)})
\in \MC$ we define $ U_{\blambda} $ as the product
\begin{equation}U_{\blambda} := \prod_{j=1}^r u_{ i_j , j }
\end{equation}
where $ i_j $ is any number from the $ j$'th component $\T^{\lambda^{(j)}}$ of
$\bT^{\blambda}$ with the convention that $u_{ i_j , j }:=1$
if it is empty.}

We have now gathered all the ingredients of our \textcolor{black}{cellular basis} element
corresponding to $\bT^{\blambda}$.
\begin{defi} Let $\blambda \in \MC $. Then we define $ m_{\blambda} \in \YY$\nomenclature[18]{$m_{\blambda}$}{The Murphy element associated with $\blambda$} via
\begin{equation}\label{dominant}m_{\blambda} := U_{\blambda} \textcolor{black}{\Elambda} x_{\blambda}.
\end{equation}
\end{defi}
The following Lemmas contain some basic properties for $ m_{\blambda}$.
\begin{lem}{\label{l5}}
The following properties for $ m_{{\blambda}} $ are true.
\begin{itemize}
\setlength\itemsep{-.5em}
\item[(1)] \textcolor{black}{ $ U_{\blambda} \textcolor{black}{\Elambda} $ is an idempotent.
  It is independent of the choices of $ i_j$'s and so also
  $ m_{\blambda} $ is independent of the choices of $ i_j$'s.}

\item[(2)] For $ i $ in the $ j$'th component of $\bT^{\blambda}  $ (that is $ p_{\blambda}(i)=j$)
we have $ t_im_{\blambda}= m_{\blambda}  t_i =  \xi^{j} m_{\blambda}$.

\item[(3)] The factors $U_\blambda$, $\textcolor{black}{\Elambda}$ and $x_\blambda$ of $m_\blambda$ commute with each other.

\item[(4)] If $i$ and $j$ occur in the same block of $A_{\blambda}  $ then
$m_{\blambda}e_{ij}= e_{ij} m_{\blambda}=  m_{\blambda} $.

\item[(5)] If $i$ and $j$ occur in two different blocks of $A_{\blambda}  $ then
$m_{\blambda}e_{ij}=0 =e_{ij} m_{\blambda} $.

\item[(6)]  For all $ w\in \Si_{\blambda} $ we have $m_{\blambda}g_w = g_w m_{\blambda} =q^{\ell(w)}m_{\blambda}$.

\end{itemize}

\end{lem}

\begin{demo}
The properties (1) and (2) are consequences of the definitions,
whereas (3) follows from (2) and Lemma \ref{pro1}.
\textcolor{black}{The property (4) follows from (2) and (3) since $U_{\blambda}\textcolor{black}{\Elambda} e_{ij }=U_{\blambda}\textcolor{black}{\Elambda} $
in that case.
Similarly, under the hypothesis of (5) we have that $U_{\blambda}\textcolor{black}{\Elambda} e_{ij }=0 $ and so also (5) follows
from (2) and (3)}.
To show (6) we note that for \textcolor{black}{$ s_i \in \Si_{\norm{ \blambda }}$}
we have that $ \textcolor{black}{\Elambda} g_i^2 = \textcolor{black}{\Elambda} (1 + (q-q^{-1}) g_i) $.
\textcolor{black}{Since $ \Si_{\blambda } $ is a subgroup of $ \Si_{\norm{ \blambda }  }$} the statement of (6) reduces to the similar Iwahori-Hecke algebra statement for $ x_{\blambda} $ which is proved for instance in
\cite[Lemma 3.2]{Mat}.
\end{demo}

\begin{obs}
Note that $ i $ and $ j $ are in the same block of $ \A_{\blambda} $ if and only if they are in the same component of $\bT^{\blambda}  $.
However, the enumerations of the blocks of $ \A_{\blambda} $ and the components of $\bT^{\blambda}  $ are different since $\bT^{\blambda}  $ may have
empty components and so in part (2) of the Lemma we cannot replace one by the other.
\end{obs}

\begin{lem}\label{lem1} Let $\blambda \in \MC $ and suppose
that $w\in \Si_n$. Then $  m_{\blambda} g_w  g_{i}= $
$$
\left\{\begin{array}{l} m_{\blambda} g_{ws_i}\mbox {       if } \ell(ws_i)>\ell(w)\\
m_{\blambda}  g_{ws_i}\mbox {    if } \ell(ws_i)<\ell(w) \mbox{ and
} i,  i+1   \mbox{ are in different blocks of } (A_{\blambda})w
\\
m_{\blambda} (g_{ws_i}+(q-q^{-1}) \textcolor{black}{g_w})\mbox { if }
\ell(ws_i)<\ell(w) \mbox{ and } i,  i+1   \mbox{ are in the same
block of } (A_{\blambda})w.
\end{array}\right.$$
\end{lem}




\color{black}

\begin{demo}
Suppose that $\ell(ws_i) >  \ell(w)$ and let $s_{j_1}\cdots s_{j_k}$
be a reduced expression for $w$. Then $s_{j_1}\cdots s_{j_k}s_i$ is
a reduced expression for $ws_i$ and so $g_{ws_i} = g_wg_{s_i}$ by
definition. On the other hand, if $\ell(ws_i)<\ell(w)$ then $w$ has
a reduced expression ending in $s_i$,
therefore
$$g_wg_{i}=g_{ws_i}g_i^2=g_{ws_i}(1+(q-q^{-1})e_ig_i)=g_{ws_i}+(q-q^{-1})g_we_i.$$
On the other hand, from Lemma
\ref{pro1}
we have that
$ \textcolor{black}{\Elambda} g_{w} e_i = g_w E_{A_{\blambda}w} e_i $ which is equal to $ g_w E_{A_{\blambda} w}  $ or zero depending on
whether $ i $ and $ i+1 $ are in the same block of $ A_{\blambda} $ or not.
This concludes the proof of the
Lemma
\end{demo}




\medskip
With these preparations, we are in position to give the definition of the set of
elements that turn out to contain the cellular basis for $ \YY$.
\begin{defi} Let $\blambda \in \MC $ and suppose that $\mathfrak{s}$ and
$\mathfrak{t}$ are row standard multitableaux of shape $\blambda$. Then we define
\begin{equation}\label{standardbasis}
m_{\Bs\bT}:=g_{d(\Bs)}^*m_\blambda g_{d(\bT)}.
\end{equation}
In particular we have $m_{\blambda} = m_{ \bT^{\blambda} \bT^{\blambda}}  $.
\end{defi}

\textcolor{black}{Recall that Murphy introduced the elements $ x_{\s\T} $ of the Iwahori-Hecke algebra
  $ {\mathcal H}_n(q) $, via
  \begin{equation}{\label{recalMurphy}}
x_{\s\T}:=g_{d(\s)}^{\ast}x_{\lambda} g_{d(\T)}
\end{equation}
for $ \s,\T$ row standard $\lambda$-tableaux. We consider our elements $m_{\Bs\bT}$ as the natural generalization of these $ x_{\s\T} $ to the Yokonuma-Hecke
algebra.}

\medskip
Clearly we have $m_{\Bs\bT}^{\ast}=m_{\bT\Bs}$,
as one sees from the definition of $\ast$.

\medskip
\color{black}
Let $\blambda \in \MC $ and set $\alpha := \norm{\blambda } $.
 We have a canonical decomposition
 of the corresponding Young subgroup
\begin{equation}{\label{canonicaldec}}
  \Si_{\alpha} = \Si_{\alpha_1} \times \Si_{\alpha_2} \times \ldots \times  \Si_{\alpha_r}
  \end{equation}
 where $ \Si_{\alpha_1} $ is the subgroup of $ \Si_n$ permuting $ \{ 1, 2, \ldots, \alpha_1 \} $, whereas
 $ \Si_{\alpha_2} $ is the subgroup permuting $ \{ \alpha_1+1 , \alpha_1+ 2, \ldots, \alpha_1 +\alpha_2 \} $,
 and so on. Note that this notation deviates slightly from the notation introduced above where
 $ \Si_{\alpha_i} $ is the symmetric group on the numbers $ \{1,2,\ldots, \alpha_i \} $;
 this kind of abuse of notation, that we shall use frequently in the following, should not cause confusion.
We now define
\begin{equation}
\YYlambda := \spa_R \{ U_{\blambda} E_{\blambda} g_w  | w \in \Si_{ \alpha }   \}.
\end{equation}
We have the following Lemma.
\begin{lem}{\label{lemmasubalgebra}}
  $ \YYlambda $ is a subalgebra of $ \YY $. Its identity element is given by
  the central idempotent $U_{\blambda} \textcolor{black}{\Elambda} $.
   There is an isomorphism between the Young-Hecke algebra $  \mathcal{H}_{\alpha}(q)$ and $ \YYlambda $
 given by
\begin{equation}
\mathcal{H}_{\alpha}(q)
 \longrightarrow  \YYlambda    , \, \, \, g_w \mapsto  U_{\blambda} \textcolor{black}{\Elambda} g_w \, \, \, \,  \mbox{     where }  w \in \Si_{ \alpha  }.
\end{equation}
Using the canonical isomorphism
$ \mathcal{H}_{\alpha}(q) \cong \mathcal{H}_{\alpha_1}(q) \otimes  \cdots \otimes \mathcal{H}_{\alpha_r}(q) $
it is given by
\begin{equation}
\mathcal{H}_{\alpha_1}(q) \otimes  \cdots \otimes \mathcal{H}_{\alpha_r}(q)
 \longrightarrow  \YYlambda    , \, \, \, a_1 \otimes \cdots \otimes a_r \mapsto  U_{\blambda} \textcolor{black}{\Elambda} \, a_1  \cdots a_r
\, \, \, \,  \mbox{     where }  a_i  \in \mathcal{H}_{\alpha_i}(q).
\end{equation}

\end{lem}

\begin{demo}
  From Lemma \ref{l5} we know that $  U_{\blambda} \textcolor{black}{\Elambda} $ is an idempotent.
  For $ w \in \Si_{ \alpha}   $
  it
  commutes with
  $ g_w $  as can
  be seen by combining the Yokonuma-Hecke algebra relation (\ref{r3}) with  Lemma \ref{l5},
  and hence it is central.
  Moreover, for $ s_i \in \Si_{ \alpha } $ we have that
    $ \textcolor{black}{\Elambda} g_i^2 =   \textcolor{black}{\Elambda}(1+(q-q^{-1})g_i) $, as mentioned in the proof of Lemma
    \ref{l5}, and so we also have $$ U_{\blambda}\textcolor{black}{\Elambda} g_i^2 =   U_{\blambda}\textcolor{black}{\Elambda}(1+(q-q^{-1})g_i)
    \, \, \,     \mbox{  for   } s_i \in \Si_{ \alpha }. $$
  It follows from this that $ \YYlambda $ is a subalgebra of $ \YY$ and that
  there is a homomorphism from $ \mathcal{H}_{\alpha}(q)  $ to $ \YYlambda$
  given by $g_w \mapsto  U_{\blambda} \textcolor{black}{\Elambda} g_w  $.
  On the other hand, it is clearly surjective and using Juyumaya's basis (\ref{juyu}), we get that $\YYlambda $ has the same dimension as
  $ \mathcal{H}_{ \alpha }(q)  $ and so the Lemma follows.
\end{demo}

\begin{obs}{\label{importantremark}}
  Suppose still that $ \blambda \in \MC$ with $ \alpha:= \norm{ \blambda } $.
 If $\Bs$ and $ \bT $  are $\blambda$-multitableaux of
the initial kind, we may view $ m_{\Bs\bT} $ as usual Murphy elements
of the Young-Hecke algebra.
Indeed, in this case $d(\Bs) $ and $d(\bT) $ belong to $ \Si_{\alpha}$
and hence, using the above Lemma, we have that $ g_{d(\Bs)} $ and $ g_{d(\bT)} $ commute with $U_{\blambda} \textcolor{black}{\Elambda} $. In particular, we have that
\begin{equation}{\label{inparticularmst}}
m_{\Bs\bT}= g_{d(\Bs)}^{*}U_{\blambda}\textcolor{black}{\Elambda}x_{\blambda}g_{d(\bT)} =
 U_{\blambda}\textcolor{black}{\Elambda}g_{d(\Bs)}^{*}x_{\blambda}g_{d(\bT)}.
\end{equation}
We have that $\Si_{\blambda}$ is subgroup of $ \Si_{ \alpha }$ compatible with
 ({\ref{canonicaldec}}) in the sense that
\begin{equation}
\Si_{\blambda}  = \Si_{\lambda^{(1)}} \times \Si_{\lambda^{(2)}} \times \cdots \times \Si_{\lambda^{(r)}},
\mbox{ where } \Si_{\lambda^{(i)}} \leq \Si_{\alpha_i}.
\end{equation}
Here $ \Si_{\lambda^{(i)}} $ is subject to the same
abuse of notation as $ \Si_{\alpha_i}$.
We then get a corresponding factorization
\begin{equation}
x_{\blambda} = x_{\lambda^{(1)}} x_{\lambda^{(2)}} \cdots x_{\lambda^{(r)}}
\end{equation}
where $ x_{\lambda^{(i)}} = \sum_{w \in \Si_{\lambda^{(i)}}} g_w $.
Since $ \Bs $ and $ \bT $ are of the initial kind we get decompositions, corresponding to the decomposition in
({\ref{canonicaldec}})
 \begin{equation}{\label{decompositionref}}
d(\Bs)=\left( d(\s^{(1)}), d(\s^{(2)}), \cdots, d(\s^{(r)}) \right) \quad \mbox{ and }\quad
 d(\bT)=\left(d(\T^{(1)}), d(\T^{(2)}), \cdots,  d(\T^{(r)})\right).
\end{equation}
 But then from ({\ref{inparticularmst}}) we get a decomposition of $ m_{\Bs\bT}$ as follows
 \begin{equation}{\label{decompositionasfollows}}
   \begin{array}{l}
m_{\Bs\bT} =
U_{\blambda}\textcolor{black}{\Elambda}g_{d(\s^{(1)})}^{*}x_{\lambda^{(1)}}g_{d(\T^{(1)})}
g_{d(\s^{(2)})}^{*}x_{\lambda^{(2)}}g_{d(\T^{(2)})}
\cdots g_{d(\s^{(r)})}^{*}x_{\lambda^{(r)}}g_{d(\T^{(r)})} \\
\, \, \, \, \, \, \, \, \, \, \, \, \, \, \, \, \, \,\,
= U_{\blambda}\textcolor{black}{\Elambda}
x_{\s^{(1)} \T^{(1)}} x_{\s^{(2)} \T^{(2)}}   \cdots   x_{{\s^{(r)} \T^{(r)}}}
\end{array}
   \end{equation}
where $ x_{\s^{(i)} \T^{(i)}}  := g_{d(\s^{(i)})}^{*}x_{\lambda^{(i)}}g_{d(\T^{(i)})}$.
Under the isomorphism of the Lemma, we then get via ({\ref{decompositionasfollows}}) that
$ m_{\Bs\bT} $ corresponds to
\begin{equation}x_{\s^{(1)} \T^{(1)}} \otimes
x_{\s^{(2)} \T^{(2)}} \otimes  \cdots \otimes  x_{{\s^{(r)} \T^{(r)}}} \in \mathcal{H}_{\alpha}(q)
\end{equation}
where each $ x_{{\s^{(i)} \T^{(i)}}} \in {\mathcal H}_{\alpha_i}(q)$ is a usual Murphy element.
This explains the claim made in the beginning of the Remark.
\end{obs}

\color{black}

\medskip
Our goal is to show that with $ \bs $ and \textcolor{black}{$ \bT $}
running over standard multitableaux for multipartitions, the $
m_{\Bs\bT} $'s form a cellular basis for $ \YY$. A first property of $
m_{\Bs\bT} $ is given by the following Lemma.
\begin{lem}\label{lem0} Suppose that  $\blambda \in \MC $ and that
$\Bs$ and $ \bT $ are $\blambda$-multitableaux. If $ i $ and $ j $ occur
in the same component of $ \bT $ then we have that $ m_{\Bs \bT} e_{i j} =m_{\Bs \bT}$. Otherwise
$ m_{\Bs \bT} e_{i j} = 0 $. A similar statement holds for $e_{i j} m_{\Bs \bT}   $.
\end{lem}

\begin{demo}
\textcolor{black}{From the definitions we have and Lemma \ref{pro1} we have that
$$ m_{\Bs \bT} e_{i j} = g_{ d(\Bs)}^{\ast} x_\blambda \textcolor{black}{\Elambda} U_{\blambda}  g_{d(\bT)} e_{ij} =
g_{ d(\Bs)}^{\ast} m_{\blambda}  e_{ id(\bT)^{-1}  ,j d(\bT)^{-1} } g_{d(\bT)}.
$$
But $ i $ and $ j $ occur in the same component of $ \bT$ iff $ id(\bT)^{-1} \!$ and $ j d(\bT)^{-1} \! $ occur in the same
block of $ A_{\blambda} $ and so the first part of the Lemma follows from (4) and (5) of Lemma \ref{l5}. The second
part is proved similarly or by applying $\ast $ to the first part.}
\end{demo}

\begin{lem}\label{l1}
Let $\blambda \in \MC$
and let $\Bs$
and $\bT$ be row standard $\blambda$-multitableaux.
Then for $h\in \YY$ we have that $ m_{\Bs\bT}h $
is a linear combination of terms of the form
$m_{\Bs\V}$ where $\V$ is a row
standard $\blambda$-multitableau. A similar statement holds for $h m_{\Bs\bT}$.
\end{lem}

\begin{demo}
Using Lemma \ref{lem1} we get that $ m_{\Bs\bT}h $ \textcolor{black}{is}
a linear combination of terms of the form $ m_{\Bs
\bT^{\blambda}}g_{w}$. For each such $ w $ we find a $y\in
\Si_{\blambda}$ and a distinguished right coset representative $d$
of $\Si_{\blambda}$ in $\Si_{n}$ such that $w=yd$ and
$\ell(w)=\ell(y)+\ell(d)$. Hence, via Lemma \ref{lem1} we get that
$$ m_{\Bs\bT}h = q^{\ell(y)}  m_{\Bs \bT^{\blambda}}g_{d}
= q^{\ell(y)}  m_{\Bs \V} $$
where $ \V = \bT^{\blambda}g_{d} $ is row standard. This proves the Lemma in the case
$ m_{\Bs\bT}h $. The case $h m_{\Bs\bT} $ is treated similarly or by applying $\ast$ to the first case.

\end{demo}



\medskip
The proof of the next Lemma is inspired by the proof of Proposition
3.18 of Dipper, James and Mathas' paper \cite{DJM}, although it
should be noted that the basic setup of \cite{DJM} is different
from ours. Just like in that paper, our proof relies on Murphy's Theorem 4.18 in \cite{Mur},
which is a key ingredient for the construction of the standard basis for $ {\mathcal H}_n(q) $.

\begin{lem}\label{l3} Suppose that $\blambda \in  \MC$ and that
$\Bs$ and $\bT$ are row standard $\blambda$-multi-tableaux. Then there are multipartitions $ \bmu \in \MP $ and
standard multitableaux $\Bu$ and $\Bv$ of shape $\bmu$, such that $\Bu \unrhd \Bs$, $\Bv
\unrhd \bT$ and such that $m_{\Bs\bT}$ is a linear combination of the corresponding
elements $m_{\Bu\Bv}$.
\end{lem}
\begin{demo} \color{black}
Let $\alpha$ be the composition $\alpha=(\alpha_1,\alpha_2, \ldots, \alpha_r)  :=\norm{\blambda}$
with corresponding Young subgroup $ \Si_{\alpha} =
\Si_{ \alpha_1 } \times \Si_{ \alpha_2 } \times \cdots \times  \Si_{ \alpha_r }$
(where some of the factors $ \Si_{ \alpha_i}$ may be trivial). Then
there exist $\blambda$-multitableaux $ \Bs_0 $ and $ \bT_0 $ of the
initial kind together with $ w_{\Bs},w_{\bT} \in \Si_n $ such that $ d(\Bs) =
d(\Bs_0) w_{\Bs} $, $ d(\bT) = d(\bT_0) w_{\bT} $ and $ \ell(d(\Bs)) = \ell(d(\Bs_0))+\ell(
w_{\Bs}) $ and $ \ell(d(\bT)) = \ell(d(\bT_0))+\ell( w_{\bT}) $. Thus, $ w_{\Bs} $ and $ w_{\bT} $
are distinguished right coset representatives for $ \Si_{\alpha} $
in $ \Si_n$ and using Lemma \ref{lem1},
together with its left action version obtained via $\ast$, we get that $ m_{\Bs\bT} =
g_{w_{\Bs}}^{\ast} m_{\Bs_0 \bT_0} g_{w_{\bT}} $.
Let $\Bs_0=( \s_0^{(1)}, \s_0^{(2)}, \ldots, \s_0^{(r)}) $ and
$\bT_0=( \T_0^{(1)}, \T_0^{(2)}, \ldots, \T_0^{(r)}) $.
Then
under the isomorphism of Lemma {\ref{lemmasubalgebra}} we have that
$ m_{\Bs_0 \bT_0} $ corresponds to
\begin{equation}{\label{that is we have}}
x_{\s^{(1)} \T^{(1)}} \otimes x_{\s^{(2)} \T^{(2)}} \otimes  \cdots \otimes  x_{{\s^{(r)} \T^{(r)}}} \in \mathcal{H}_{\alpha}(q)
\end{equation}
as explained in Remark {\ref{importantremark}}.
On each of the factors $ x_{\s_0^{(i)}\T_0^{(i)}} $ we now
use Murphy's result Theorem 4.18 of \cite{Mur}
thus concluding that $x_{\s_0^{(i)}\T_0^{(i)}} $  is a linear combination of terms of the
form $ x_{\U_0^{(i)}\V_0^{(i)}} $ where $ \U_0^{(i)} $ and $
\V_0^{(i)} $ are standard $\mu_0^{(i)} $-tableaux on the numbers
permuted by $ \Si_{\alpha_i}$
and satisfying $\U_0^{(i)} \unrhd
\s_0^{(i)}$ and $\V_0^{(i)} \unrhd \T_0^{(i)}$.
Letting $ \bmu:= ( \mu_0^{(1)}, \mu_0^{(2)}, \ldots, \mu_0^{(r)}) $,
$\Bu_0:=( \U_0^{(1)}, \U_0^{(2)}, \ldots, \U_0^{(r)}) $ and
$\Bv_0:=( \V_0^{(1)}, \V_0^{(2)}, \ldots, \V_0^{(r)}) $
and using the isomorphism of Lemma {\ref{lemmasubalgebra}} in the other direction
we then get that $ m_{\Bs_0 \bT_0} $ is a
linear combination of terms $ m_{\Bu_0 \Bv_0} $  where $ \Bu_0$ and $ \Bv_0$
are standard $ \bmu $-multitableaux such that $ \Bu_0 \unrhd \Bs_0 $ and
$ \Bv_0 \unrhd \bT_0 $.
Hence $ m_{\Bs\bT} = g_{w_{\Bs}}^{\ast} m_{\Bs_0 \bT_0} g_{w_{\bT}} $ is a
linear combination of terms $ g_{w_{\Bs}}^{\ast}  m_{\Bu_0 \Bv_0} g_{w_{\bT}} $.
On the other hand, $ \Bu_0 $ and $ \Bv_0 $ are of the initial kind, and so we get
$
g_{w_{\Bs}}^{\ast}  m_{\Bu_0 \Bv_0} g_{w_{\bT}} =
m_{\Bu_0 w_{\Bs}, \Bv_0
w_{\bT}}
$
since $ w_{\Bs} $ and $ w_{\bT} $ are distinguished
right coset representatives for $ \Si_{\alpha} $ in $ \Si_n$. This
also implies that $ \Bu_0 w_{\Bs} \unrhd \Bs_0 w_{\Bs} = \Bs $ and
$ \Bv_0 w_{\Bs} \unrhd \bT_0 w_{\Bs} = \bT $ proving the Lemma.
\end{demo}

\color{black}

\begin{coro}\label{co1}
Suppose that $\blambda \in \MC $ and that $\Bs$
and $\bT$ are row standard $\blambda$-multi-tableaux. If $h\in \YY$, then
$m_{\Bs\bT}h$ is a linear combination of terms of
the form $m_{\Bu\Bv}$ where $\Bu$ and $\Bv$ are standard $\bmu$-multitableaux
for some multipartition $\bmu \in \MP$ and $\Bu\unrhd \Bs$ and $\Bv
\unrhd \bT$. A similar statement holds for $h m_{\Bs\bT}$.
\end{coro}

\begin{demo} This is now immediate from the Lemmas \ref{l1} and \ref{l3}.
\end{demo}

\medskip
So far our construction of the cellular basis has followed the
layout used in \cite{DJM}, with the appropriate adaptions. But to
show that the $m_{\Bs\bT}$'s generate $ \Y(q) $ we shall deviate from
that path. We turn our attention to the $R$-subalgebra
$\mathcal{T}_n$ of $\Y$ generated by $t_1,t_2,\ldots,t_n$.  By the
faithfulness of $ V^{\otimes n} $, it is isomorphic to the
subalgebra $ \TT_n  \subset {\rm End}( V^{\otimes n})   $ considered
above. Our proof that the elements $m_{\Bs\bT}$ generate $ \YY $
relies on the, maybe surprising, fact that $\mathcal{T}_n$ is
compatible with the $\{ m_{\Bs\bT} \}$, in the sense that the elements
\textcolor{black}{$\{ m_{\Bs\Bs}\} $ where $ \Bs$ is a
multitableau corresponding to a  one-column multipartition induce a basis
for $\mathcal{T}_n$.}

\medskip
As already mentioned, we consider our $m_{\Bs\bT}$ as the
natural generalization of Murphy's standard basis to $ \YY$. It is interesting to note that Murphy's standard basis and its generalization have
already before manifested 'good' compatibility properties of the above kind.

\medskip
Let us first define a one-column $r$-multipartition to
\textcolor{black}{be} an element of $\MP$ of the form $ ((1^{c_1}),
\ldots,  (1^{c_r})) $ and let $\MP^1 $ be the set of one-column
$r$-multipartitions. Note that there is an obvious bijection between
$\MP^1 $ and the set of usual compositions in $ r $ parts. We define
$$ \std_{n,r}^1 := \left\{ \Bs \,  | \Bs \in \std(\blambda) \mbox{ for } \blambda \in \MP^1      \right\}.$$
Note that $ \std_{n,r}^1 $ has cardinality $ r^n $ as follows from the multinomial formula.
\begin{lem}
For all $ \Bs
\in \std_{n,r}^1 $, we have that $ m_{ \Bs \Bs} $ belongs to $
{\mathcal T }_n$.
\end{lem}

\begin{demo}
Let $\Bs$ be an element of $  \std_{n,r}^1$. It general, it is useful to think of $d(\Bs) \in \Si_n$ as the row reading of $ \Bs$,
that is the element obtained by
reading the components of $\Bs$ from left to right, and the rows of
each component from top to bottom.

\medskip
We show by induction on
$\ell(d(\Bs))$ that $m_{\Bs\Bs} $ belongs to $\mathcal{T}_n$.
If $\ell(d(\Bs))=0 $ then $ x_{\blambda} = 1 $ and so $m_{\Bs\Bs} = U_{\blambda} \textcolor{black}{\Elambda}$ that certainly
belongs to $\mathcal{T}_n$. Assume that the statement holds for
all multitableaux $ \Bs^{\prime} \in  \std_{n,r}^1 $ such that $ \ell(d(\Bs^{\prime})) < \ell(d(\Bs)) $.
Choose $i $ such that $i$ occurs in $ \Bs   $ to the right of $ i+1$:
such an $i$ exists because $ \ell(d(\Bs)) \neq 0$. Then we can apply the
inductive hypothesis to $ \Bs s_i $, that is $ m_{ \Bs s_i\, \Bs s_i }
\in \mathcal{T}_n$. But then
\begin{equation}{\label{inductive_step}}
m_{ \Bs  \Bs  } = g_{ d(\Bs)}^{\ast} m_{\blambda} g_{ d(\Bs)} = g_i m_{
\Bs s_i\, \Bs s_i } g_i = g_i m_{ \Bs s_i\, \Bs s_i } (g_i^{-1}
+(q-q^{-1})e_i).
\end{equation}
But $ g_i m_{ \Bs s_i\, \Bs s_i } g_i^{-1} $ certainly belongs to
$\mathcal{T}_n$, as one sees from relation (\ref{r3}). Finally, from
Lemma \ref{lem0} we get that $ m_{ \Bs s_i\, \Bs s_i } e_i = 0 $, thus
proving the Lemma.
\end{demo}

\begin{lem}\label{l6}
Suppose that $\blambda \in \MC $ and let $\Bs$ and $ \bT $
be $\blambda$-multitableaux. Then for all $k=1,\ldots,n$ we have that
$$m_{\bT\Bs}t_k=\xi^{p_{\Bs}(k)}m_{\bT\Bs} \, \, \, \,  \mbox{and        }  \, \, \, \, t_k m_{\bT\Bs}=\xi^{p_{\bT}(k)}m_{\bT\Bs}.$$
\end{lem}

\begin{demo}
From (\ref{r3}) we have that $g_w t_k=t_{kw^{-1}}g_w$ for
all $w\in\Si_n$. Then, by Lemma {\ref{l5}}(2) we have
$$\begin{array}{ll}m_{\bT^{\blambda}\Bs}t_k=
m_{\blambda}g_{d(\Bs)}t_k= m_{\blambda}g_{d(\Bs)}t_k= m_{\blambda} t_{ k d(\Bs)^{-1}}  g_{d(\Bs)} =
\xi^{p_\blambda(kd(\Bs)^{-1})}m_{\bT^{\blambda}\Bs}.\end{array}$$ On the
other hand, since $\Bs=\bT^{\blambda}d(\Bs)$ we have that
$p_{\blambda}(kd(\Bs)^{-1})=p_{\Bs}(k)$ and hence
$m_{\bT^{\blambda}\Bs}t_k=\xi^{p_{\Bs}(k)}m_{\bT^{\blambda}\Bs}$. Multiplying
this equality on the left by $g_{d(\bT)}^{\ast}$, the proof of the first formula is completed.
The second formula is shown similarly or by applying $ \ast$ to the first.
\end{demo}

\medskip
Our next Proposition shows that the set $ \{ m_{\Bs \Bs } \} $, where $ \Bs \in \std_{n,r}^1$, forms a basis for $\mathcal{T}_n$, as promised.
We already know that $ m_{\Bs \Bs} \in \mathcal{T}_n $ and that the cardinality of
$ \std_{n,r}^1 $ is $ r^n$ which is the dimension of $\mathcal{T}_n$,
but even so the result is
not completely obvious, since we are working over the ground ring
$ R $ which is not a field.
\begin{propos}\label{pr1} $\{ m_{\Bs \Bs } \, | \, \Bs \in \std_{n,r}^1\}$ is an $R$-basis for $\mathcal{T}_n$.
\end{propos}

\begin{demo}
Recall that we showed in the proof of Theorem {\ref{faithful}} that
$$ V_{i_1, i_2 \ldots, i_n}= \spa_R \{ v_{i_1}^{j_1}\otimes
v_{i_2}^{j_2}\otimes \cdots \otimes v_{i_n}^{j_n}\,  |\,  j_k \in
{\mathbb Z}/r {\mathbb Z} \} $$ is a faithful ${\cal T}_n $-module
for any fixed, but arbitrary, set of lower indices. Let $ \seq $ be
the set of sequences $ \underline{i}=(i_1, i_2, \ldots, i_n) $ of numbers $ 1\leq
i_j\leq n$. Then we have that
\begin{equation}\label{directsumtensor} V^{\otimes n} = \bigoplus_{ \underline{i}
\in \seq } V_{ \underline{i}}
\end{equation}
and of course $ V^{\otimes n} $ is a
faithful ${\cal T}_n $-module, too. For $ \Bs\in \std_{n,r}^1$ and $
\underline{i} \in \seq $ we define
\begin{align} v^{\Bs}_{\underline{i}} := v^{ j_1}_{i_1}  \otimes v^{
j_2}_{i_2} \otimes \ldots \otimes v^{ j_n}_{i_n} \in
V_{\underline{i}}\label{ten}\end{align} where $ (j_1, j_2, \ldots
j_n) := (p_{\Bs}(1), p_{\Bs}(2),  \ldots \, p_{\Bs}(n) )$. Then $ \{
v^{\Bs}_{\underline{i}} \, | \, \Bs \in \std_{n,r}^1, \underline{i} \in
\seq \} $ is an $ R$-basis for $ V^{\otimes n}$. We now claim the
following formula in $ V_{{\underline{i}}}$:
\begin{equation}{\label{wethenhave}}
 v^{\bT}_{\underline{i}}m_{\Bs \Bs} = \left\{ \begin{array}{rl}
v^{\bT}_{\underline{i}} & \mbox{ if } \Bs = \bT \\ 0 & \mbox{
otherwise. } \end{array} \right.
\end{equation}
We show it by induction on $\ell(d(\Bs))$. If $\ell(d(\Bs))=0$, then $\Bs =
\bT^{\blambda} $ where $ \blambda $ is the shape of $ \Bs $. We have $ x_{\blambda} = 1$ and so
$ m_{\Bs \Bs} = m_{\blambda} = U_{\blambda} \textcolor{black}{\Elambda} $.
We then get ({\ref{wethenhave}}) directly from the definitions of
$ U_{\blambda}$ and $ \textcolor{black}{\Elambda} $
together with Lemma {\ref{o1}}.

Let now $\ell(d(\Bs)) \neq 0$ and assume that (\ref{wethenhave}) holds
for multitableaux $ \Bs^{\prime} $ such that $ \ell(d(\Bs^{\prime})) <
\ell(d(\Bs)) $. We choose $j $ such that $j$ occurs in $ \Bs   $ to the
right of $ j+1$. Using (\ref{inductive_step}) we have that $m_{ \Bs
\Bs  } = g_j m_{ \Bs s_j \, \Bs s_j } g_j^{-1}$. On the other hand, $ j
$ and $ j+1 $ occur in different components of $ \Bs$  and so by
Definition \ref{tensoraction} of the $ \Y(q)$-action in $
V^{\otimes n} $ we get that $  v^{\Bs}_{{\underline{i}}}g^{\pm 1}_j
= v^{\Bs s_j }_{ {\underline{i}s_j}} $, corresponding to the first case of ({\ref{operatorG}}).
Hence we get via the
inductive hypothesis that
$$  v^{\Bs }_{{\underline{i}}}m_{ \Bs  \Bs  } =  v^{\Bs }_{{\underline{i}}}g_j m_{ \Bs s_j \, \Bs s_j } g_j^{-1} =
v^{ \Bs s_j}_{  {\underline{i}}s_j} m_{ \Bs s_j \, \Bs s_j }g_j^{-1} =
 v^{\Bs s_j}_{{\underline{i}s_j}}g_j^{-1} = v^{\Bs}_{{\underline{i}}} $$
which shows the first part of (\ref{wethenhave}).

If $ \Bs \neq \bT $ then we essentially argue the same way. We choose
$ j $ as before and may apply the inductive hypothesis to $ \Bs s_j$.
We have that $  v^{\bT }_{{\underline{i}}}m_{ \Bs  \Bs  }  =  v^{\bT
}_{{\underline{i}}}g_j m_{ \Bs s_j \, \Bs s_j } g_j^{-1} $ and so need
to determine $  v^{\bT }_{{\underline{i}}} g_j $. This is slightly
more complicated than in the first case, but using the Definition
\ref{tensoraction} of the $ \Y(q)$-action in $ V^{\otimes n} $ we
get that $v^{\bT }_{{\underline{i}}}g_j$ is always an $R$-linear
combination of the vectors $ v^{\bT s_j}_{{\underline{i}}s_j}$ and $
v^{\bT s_j}_{{\underline{i}}}$: indeed in the cases $ s=t $ of
Definition \ref{tensoraction} we have that $ p_{\bT s_j}(s) =
p_{\bT}(s)$. But $ \Bs \neq \bT $ implies that $ \Bs s_j \neq\bT s_j $
and so we get by the inductive hypothesis that
\begin{equation}\label{aswgby}
v^{\bT }_{{\underline{i}}}m_{ \Bs  \Bs  }  = v^{\bT }_{{\underline{i}}} g_j m_{ \Bs s_j \, \Bs s_j } g_j^{-1}  = 0
\end{equation}
and ({\ref{wethenhave}}) is proved.

From ({\ref{wethenhave}}) we now deduce that $ \sum_{\Bs \in \std_{n,r}^1}  v_{\underline{i}}^{\bT} m_{ \Bs \Bs } = v_{\underline{i}}^{\bT}
$ for any $ \bT$ and $ \underline{i}$, and hence
\begin{equation}\label{sinceisfaithful}
\sum_{\Bs \in \std_{n,r}^1} m_{ \Bs \Bs }  =1
\end{equation}
since $ V^{\otimes n} $ is faithful and the $ \{ v_{\underline{i}}^{\bT} \} $ form a basis for $ V^{\otimes n} $.
We then get that
\begin{equation}\label{t_igenerate}
 t_i = t_i 1 = \sum_{\Bs \in \std_{n,r}^1} t_i m_{ \Bs \Bs } = \sum_{\Bs \in \std_{n,r}^1} \xi^{p_{\Bs}(i)} m_{ \Bs \Bs }
\end{equation}
and hence, indeed, the set $\{  m_{\Bs \Bs} \, | \, \Bs \in \std_{n,r}^1 \} $
generates $ {\cal T }_n$.  On the other hand, the $ R $-independence of $  \{  m_{\Bs \Bs} \} $
follows easily from ({\ref{wethenhave}}), via evaluation on the vectors $
v_{\underline{i}}^{\bT} $. The Theorem is proved.
\end{demo}



\begin{teo}\label{cel} The algebra $\Y(q)$ is a free $R$-module with basis
$$\B=\left\{m_{\Bs\bT}\mid\Bs,\bT\in\std(\blambda) \mbox{ for some multipartition } \blambda \mbox{ of } n \right\}.$$
Moreover, $(\B,\MP)$ is a cellular basis of $\Y(q)$ in the sense of
Definition \ref{cellular}.
\end{teo}

\begin{demo} From Proposition \ref{pr1}, we have
that $1$ is an $R$-linear combination of elements $m_{\Bs\Bs}$ where
$\Bs$ are certain standard multitableaux. Thus, via Corollary
\ref{co1} we get that $\B$ spans $\Y(q)$.
On the other hand, the cardinality of $\B$ is $r^nn!$
since, for example, $ \B $ is the set of tableaux for the Ariki-Koike algebra
whose dimension is $r^nn!$.
But this
implies that $\B$ is an $R$-basis for $\Y(q)$. Indeed, from Juyumaya's basis we know that $\Y(q)$
has rank $N:= r^nn!$ and any surjective homomorphism $ f: R^N \mapsto R^N $
splits since $ R^N $ is a projective \textcolor{black}{$R$-module}.

\medskip
The
multiplicative property that $ \B$ must satisfy in order to be a cellular
basis of $\Y(q)$, can now be shown by repeating the argument of Proposition 3.25 of \cite{DJM}.
For the reader's convenience, we sketch the argument.

\medskip
Let first ${\Y^{\blambda}}(q)$ be the $R$-submodule of $ \Y(q)$ spanned by
$$ \{m_{\Bs\bT} \, |  \,  \Bs, \bT  \in \std(\bmu) \mbox{ for some } \bmu \in \MP \mbox{ and  } \bmu  \rhd \blambda \}.$$
Then one checks using Lemma \ref{l3} that ${\Y^{\blambda}}$ is an ideal of $ \Y(q)$. 
Using Lemma \ref{l3} once again, we get for
$ h \in \Y(q) $
the formula
$$ m_{\bT^{\blambda} \bT} h = \sum_{\Bv} r_{\Bv} m_{\bT^{\blambda} \Bv} \mod {\Y^{\blambda}} $$
where $ r_{\Bv} \in R$. This is so because $ \bT^{\blambda} $ is a maximal element of
$ \std(\blambda)$. Multiplying this equation on the left with $ g_{d(\Bs)}^{\ast} $ we get the formula
$$ m_{\Bs \bT} h = \sum_{\V} r_{\Bv} m_{\Bs \Bv} \mod {\Y^{\blambda}} $$
and this is the multiplicative property that is required for cellularity.
\medskip
\end{demo}

As already explained in \cite{GL}, the existence of a cellular basis in an algebra $A$ has strong consequences
for the modular representation theory of $A$. Here we give an application of our cellular basis $ \B$ that goes in a somewhat different direction,
obtaining from it Lusztig's idempotent presentation of $ \Y(q)$, used in \cite{L1}, \cite{L2}.

\begin{propos}\label{prLusztig}
The Yokonuma-Hecke algebra $ \Y(q) $ is isomorphic to the
associative $R$-algebra
generated by the elements $\{ g_i  | i = 1, \ldots, n-1 \} $
and $ \{ f_{\Bs}   |\,  \Bs \in \std_{n,r}^1 \} $
subject to the following relations:
\begin{alignat}{3}
g_ig_j&=g_jg_i&&\quad\mbox{ for } |i-j|>1\label{rl1}\\
g_ig_{i+1}g_i&=g_{i+1}g_ig_{i+1}&&\quad \mbox{ for all }
i=1,\ldots,n-2\label{rl3} \\
f_{\Bs} g_i &=g_i f_{ \Bs s_i } &&\quad \mbox{ for all } \Bs, i  \label{rl4}\\
g_i^2&=1+(q-q^{-1})\sum_{ \Bs \in \std_{n,r}^1} \delta_{i,i+1}(\Bs) f_{\Bs} g_i && \quad \mbox{ for all } i \label{rl5} \\
\sum_{\Bs \in \std_{n,r}^1} f_{\Bs}& = 1   &&\quad\mbox{ for all } \Bs   \label{rl6}\\
f_{\Bs} f_{\Bs'} &=\delta_{\Bs,\Bs'}f_{\Bs} &&\quad \mbox{ for all }\; \Bs,\Bs'\in\std_{n,r}^1  \label{rl7}
\end{alignat}
where $\delta_{\Bs,\Bs'}$ is the Kronecker delta function on $\std_{n,r}^1$ and
where
we set $ \delta_{i,i+1}(\Bs) := 1 $ if $i $ and $ i+1$ belong to the same component (column) of $ \Bs $, otherwise $ \delta_{i,i+1}(\Bs) := 0 $.
Moreover, we define $ f_{\Bs s_i} := f_{\Bs} $ if $ \delta_{i,i+1}(\Bs) = 0 $.
\end{propos}

\begin{demo}
Let $ \Y^{\prime} $ be the $R$-algebra defined by the presentation of the Lemma.
Then there is an $R$-algebra homomorphism $ \varphi: \Y^{\prime} \rightarrow \Y(q) $, given
by $ \varphi(g_i ):= g_i $ and $ \varphi(f_{\Bs} ):= m_{\Bs \Bs} $. Indeed,
the $ m_{\Bs \Bs } $'s are orthogonal idempotents and have sum $ 1 $ as we see from ({\ref{wethenhave}}) and
(\ref{sinceisfaithful}) respectively. Moreover, using ({\ref{inductive_step}}), (\ref{wethenhave}) and (\ref{sinceisfaithful})
we get that
the relations (\ref{rl4}), (\ref{rl5}), (\ref{rl6}) and (\ref{rl7}) hold with $ m_{\Bs \Bs} $ replacing $ f_{\Bs } $,
and finally the first two relations hold trivially.

On the other hand, using (\ref{t_igenerate}) we get that $ \varphi $ is a surjection and since
$ \Y^{\prime} $ is generated over $ R $ by the set $ \{ g_w f_{\Bs} | w \in \Si_n, \Bs \in \std_{n,r}^1    \} $ of
cardinality $ r^n n ! $, we get that $ \varphi $ is also an injection.
\end{demo}

\begin{obs}
The relations given in the Proposition are the relations, for type $A$, of the algebra $  H_n $ considered in 31.2 of \cite{L1}
see also \cite{MS}. We would like to draw the attention to the sum appearing in the quadratic relation (\ref{rl5}), making it look
rather different than the quadratic relation of Yokonuma's or Juyumaya's presentation.
In 31.2 of \cite{L1}, it is mentioned that $H_n $ is closely related to the convolution algebra associated with a Chevalley group
and its unipotent radical and indeed
in 35.3 of \cite{L2}, elements of this algebra are found that satisfy
the relations of $  H_n $. However, we could not find a Theorem in {\it loc. cit.}, stating explicitly that $  H_n $ is isomorphic to $ \Y(q)$.
\textcolor{black}{(On the other hand, in \cite{JdA1} Jacon and Poulain d'Andecy
have recently given a simple explanation of the isomorphism $ H_n \cong \Y(q) $)}.

\end{obs}
 \section{Jucys-Murphy elements}
 In this section we show that the Jucys-Murphy elements $ J_i $ for $ \Y(q)$, introduced by
Chlouveraki and Poulain d'Andecy in \cite{MCH}, are JM-elements in the abstract sense defined by Mathas, see \cite{Mat2}.
This is with respect to the cellular basis for $ \Y(q) $ obtained in the previous section.

\medskip
We first consider the elements $J_k^{\prime}$ of $\Y(q)$ given by $J_1^{\prime}=0$ and for
 $k\geq 1$
 \begin{align}
 J_{k+1}^{\prime}&=q^{-1}(e_kg_{(k,k+1)}+e_{k-1,k+1}g_{(k-1,k+1)}+\cdots+e_{1,k+1}g_{(1,k+1)})
 \end{align}
where $ g_{(i,k+1)} $ is $ g_w $ for $ w= (i,k+1)$.
These elements are generalizations of the Jucys-Murphy elements for
the Iwahori-Hecke algebra $\mathcal{H}_n(q)$, in the sense that
we have $E_{ \bf n} J_{k}^{\prime}= E_{ \bf n} L_k$,
where $L_k$ are the Jucys-Murphy elements for $\mathcal{H}_n(q)$ defined in \cite{Mat}.

The elements $ J_i $ of $ \Y(q) $ that we shall refer to as Jucys-Murphy elements were introduced by Chlouveraki and Poulain d'Andecy in
 \cite{MCH}
via the recursion
 \begin{align}
 J_1=1 \quad \mbox{and}  \quad J_{i+1}=g_iJ_ig_i \quad \mbox{ for }
 i=1,\ldots,n-1.
 \end{align}
 The relation between $ J_{i} $ and $ J_{i}^{\prime} $ is given by
 \begin{align}
 J_{i}=1+(q^2-1) J_{i}^{\prime}.
 \end{align}
 In fact, in \cite{MCH} the elements $\{J_1,\ldots,J_n \} $, as well as the
 elements $\{ t_1,\ldots,t_n \} $, are called Jucys-Murphy elements for
 the Yokonuma-Hecke algebra.



The following definition appears for the first time in \cite{Mat2}. It formalizes the concept of
Jucys-Murphy elements.
\begin{defi}{\label{JM}}
Suppose that the $\cal R$-algebra $A$ is cellular with antiautomorphism $ \ast$ and cellular basis $\mathcal{C}=\{a_{\s\T}\mid
 \lambda\in \Lambda,\s,\T\in T(\lambda)\}$. Suppose moreover that each set $ T(\lambda) $ is
endowed with a poset structure with order relation $ \rhd_{\lambda} $.
\textcolor{black}{Then we say that a commuting} set
 $\mathcal{L}=\{L_1,\ldots,L_M\}\subseteq A $ is a family of JM-elements for $A$,
 with respect to the basis $\mathcal{C}$, if it satisfies that $ L_i^{\ast} = L_i $ for all $ i $ and if
 there exists a set of
 scalars $\{c_{\T}(i)\mid \T\in T(\lambda),\; 1\leq i\leq M\}$, called the contents of $ \lambda $, such
 that for all $\lambda\in \Lambda$ and $\T\in T(\lambda)$ we have that
\begin{equation}{\label{contents}}
 a_{\s \T}L_i=c_{\T}(i)a_{\s \T}+ \mathop{\sum_{\V\in T(\lambda)}}_{\V \rhd_{\lambda} \T}r_{\s \V}a_{\s \V} \mod A^{\lambda}
\end{equation}
$ \mbox{ for some }r_{\s \V}\in \textcolor{red}{\mathcal{R}}$.
\end{defi}
Our goal is to prove that the set
 \begin{align}
 \mathcal{L}_{\Y}:=\{L_1,\ldots,L_{2n}\mid L_k=J_k,\;L_{n+k}=t_k,\;1\leq
 k\leq n\}
 \end{align}
 is a family of JM-elements for $\Y(q)$
in the above sense.
Let us start out by stating the following Lemma.

 \begin{lem}\label{J1}
 Let $i$ and $k$ be integers such that $1\leq i<n$ and $1\leq k\leq
 n$. Then
\begin{itemize}
\setlength\itemsep{-.2em}
\item[(1)] $g_i \mbox{ and } J_k \mbox{ commute if } i\neq k-1,k. $

\item[(2)]  $\mathcal{L}_{\Y} \mbox{ is a set of commuting elements.} $

\item[(3)] $ g_i \mbox{ commutes with }  J_iJ_{i+1}  \mbox{ and } J_i+J_{i+1}. $

\item[(4)]  $g_iJ_i=J_{i+1}g_i+(q^{-1}-q)e_iJ_{i+1} \mbox{ and }  g_iJ_{i+1}=J_ig_i+(q-q^{-1})e_iJ_{i+1}.$
\end{itemize}
 \end{lem}

 \begin{demo} For the proof of (1) and (2), see \cite[Corollaries 1 and 2]{MCH}. We then prove (3) using
(1) and (2) and induction on $ i$. For $i=1$
 the two statements are trivial. For $i>1$ we have that
 $$\begin{array}{ll}
 g_iJ_iJ_{i+1}&=g_i(g_{i-1}J_{i-1}g_{i-1})(g_ig_{i-1}J_{i-1}g_{i-1}g_i)=
 g_ig_{i-1}J_{i-1}g_ig_{i-1}g_iJ_{i-1}g_{i-1}g_i\\[1mm]
 &=g_ig_{i-1}g_iJ_{i-1}g_{i-1}g_iJ_{i-1}g_{i-1}g_i
 =g_{i-1}(g_{i}g_{i-1}J_{i-1}g_{i-1}g_i)J_{i-1}g_{i-1}g_i\\[1mm]
 &=g_{i-1}J_{i+1}J_{i-1}g_{i-1}g_i
 =(g_{i-1}J_{i-1}g_{i-1})J_{i+1}g_i
 =J_{i}J_{i+1}g_i
 \end{array}$$
 and
 $$\begin{array}{ll}
 g_i(J_i+J_{i+1})=g_iJ_i+g_i^2J_ig_i
 =g_iJ_i+(1+(q-q^{-1})e_ig_i)J_ig_i\\[1mm]
 =J_ig_i+g_iJ_i(1+(q-q^{-1})e_ig_i)
 =J_ig_i+g_iJ_ig_i^{2}
 =(J_i+J_{i+1})g_i.
 \end{array}$$
 Finally, the equalities of (4) are shown by direct computations, that we leave to the reader.

 \end{demo}

Let $ \kk $ be an $R$-algebra as above, such that $ q \in \kk^{\times}$.
Let $\bT$ be a $\blambda$-multitableau and suppose that the node of $\bT$ labelled by $(x,y,k)$ is filled in with $ j $. Then
we define the \textit{quantum content} of $j$ as the element
$c_{\bT}(j):=q^{2(y-x)}\in \kk$.
We furthermore define $ {\rm res}_{\mathfrak{t}}(j):=y-x $ and then have the formula
$c_{\bT}(j)=q^{2 {\rm res}_{\mathfrak{t}}(j)}.$
When $\bT=\bT^{\blambda}$, we write $c_{\blambda}(j)$ for
$c_{\bT}(j)$.

The next Proposition is
the main result of this section.


 \begin{propos}
 $(\Y(q),\B)$ is a cellular algebra with family of JM-elements
 $\mathcal{L}_{\Y}$ and contents given by
$$
 d_{\bT}(k) := \left\{ \begin{array}{ll} c_{\bT}(k) & \mbox{ if } k =1, \ldots, n \\
 \textcolor{black}{\xi^{p_{\bT}(k)}} & \mbox{ if } k =n+1, \ldots, 2n. \end{array} \right.
 $$
 \end{propos}

 \begin{demo}
We have already proved that $ \B $ is a cellular basis for $ \Y(q) $, so we
only need to prove that the elements of $\mathcal{L}_{\Y}$ verify the conditions of Definition
{\ref{JM}}.

\medskip
For the order relation $ \rhd_{\blambda} $ on $\std(\blambda)$ we shall use the dominance order $ \rhd $ on multitableaux that was introduced
above.
By Lemma \ref{l6} the JM-condition ({\ref{contents}}) holds for $ k=n+1, \ldots, 2n $
and so we only need to check the cases $ k=1, \ldots, n $.

\medskip
\color{black}
Let us first consider the case when $\bT$ is a standard
 $\blambda$-multitableau of the initial kind.
\textcolor{black}{
Suppose $ \blambda= (\lambda^{(1)}, \ldots, \lambda^{(r)}) $, $ \bT= (\T^{(1)}, \ldots, \T^{(r)}) $ and}
 $\alpha=\norm{\blambda}$, with corresponding Young subgroup
$ \Si_{\alpha} =  \Si_{ \alpha_1 } \times \cdots \times  \Si_{ \alpha_r} $
and suppose that $k$
belongs to $\T^{(l)}$.
Since $ \bT  $ is of the initial kind we have from ({\ref{decompositionasfollows}})
a corresponding decomposition
\begin{equation}
m_{\textcolor{black}{\blambda} \bT}
= U_{\blambda}\textcolor{black}{\Elambda}
x_{\textcolor{black}{\lambda^{(1)}}\T^{(1)}} x_{\textcolor{black}{\lambda^{(2)}}\T^{(2)}}\cdots
 x_{\textcolor{black}{\lambda^{(r)}}\T^{(r)}}
\end{equation}
\textcolor{black}{where, as before, $ \blambda$ and $ \lambda^{(i)} $ as indices refer to $ \bT^{ \blambda}$ and $ \T^{\lambda^{(i)}} $.}
Hence, by (1) of Lemma
 \ref{J1} we get that
 \begin{equation}
 \begin{array}{l}
 m_{{\blambda}\bT}J_k= U_{\blambda}\textcolor{black}{\Elambda} x_{{\lambda^{(1)}}\T^{(1)}}\cdots
 x_{{\lambda^{(l)}}\T^{(l)}}J_k
 x_{{\lambda^{(l+1)}}\T^{(l+1)}}\cdots
 x_{{\lambda^{(r)}}\T^{(r)}} =\\[3mm]
 x_{{\lambda^{(1)}}\T^{(1)}}\cdots
 x_{{\lambda^{(l)}}\T^{(l)}}U_{\blambda} \textcolor{black}{\Elambda}(1+(q^2-1) J_k^{\prime})
 x_{{\lambda^{(l+1)}}\T^{(l+1)}}\cdots
 x_{{\lambda^{(r)}}\T^{(r)}}
 \end{array}
 \end{equation}
where we used Lemma \ref{pro1} to commute $  U_{\blambda}\textcolor{black}{\Elambda} $ past $ x_{{\lambda^{(1)}}\T^{(1)}}\cdots
 x_{{\lambda^{(l)}}\T^{(l)}}$.
On the other hand, by Lemma
 \ref{lem0} together with the definition of $ J_i^{\prime} $ we have
 that
\begin{equation} U_{\blambda} \textcolor{black}{\Elambda} x_{{\lambda^{(l)}} \T^{(l)}}(1+(q^2-1) J_k^{\prime}) = U_{\blambda} \textcolor{black}{\Elambda}
 x_{\lambda^{(l)}\T^{(l)}}(1+(q^2-1) L_{k}^{l})
\end{equation}
 \textcolor{black}{where $
 L_{k}^{l} = q^{-1}(g_{(k,k+1)}+
 g_{(k-1,k+1)}+\cdots+g_{(m,k+1)})$ is the $k$'th Jucys-Murphy
 element as in \cite{Mat} for the Iwahori-Hecke algebra corresponding to
 $ \Si_{\alpha_l}$, permuting the numbers $ \{ m, m+1, \ldots, m+\alpha_l-1\}$}. Thus under the isomorphism
$ \mathcal{H}_{\alpha_1}(q) \otimes  \cdots \otimes \mathcal{H}_{\alpha_r}(q)
 \cong  \YYlambda $
of Lemma {\ref{lemmasubalgebra}}
we have that the $l$'the factor of $m_{{\blambda}\bT}J_k \in \YYlambda$  is  $ x_{\lambda^{(l)}\T^{(l)}}(1+(q^2-1) L_{k}^{l}) \in \mathcal{H}_{\alpha_l}(q)$
and so we may further manipulate that element inside $ \mathcal{H}_{\alpha_l}(q)$.

Now applying
 \cite[Theorem 3.32]{Mat} we get that $
  x_{{\lambda^{(l)}}\T^{(l)}}(1+(q^2-1) L_{k}^{l}) $ is equal to
 \begin{equation}{\label{JM-initial}}
 \begin{array}{ll}
 x_{{\lambda^{(l)}}\T^{(l)}}+ (q^2-1)[{\rm res}_{\T^{(l)}}(k)]_q
  x_{{\lambda^{(l)}}\T^{(l)}}+\displaystyle\mathop{\sum_{\V\in\std(\lambda^{(l
 )})}}_{\V\triangleright \T^{(l)}}a_{\V} x_{{\lambda^{(l)}}\V}+
 \mathop{\sum_{\aaa_1,\bbb_1\in\std(\mu^{(l)})}}_{\mu^{(l)}\triangleright
 \lambda^{(l)}}r_{\aaa_1\bbb_1}x_{\aaa_1\bbb_1}\\= q^{2({\rm
 res}_{\T^{(l)}}(k))}x_{{\lambda^{(l)}}\T^{(l)}}+\displaystyle\mathop{\sum_{\V\in\std(\lambda^{(l)})}}_{\V\triangleright
 \T^{(l)}}a_{\V}x_{{\lambda^{(l)}}\V}+
 \mathop{\sum_{\aaa_1,\bbb_1\in\std(\mu^{(l)})}}_{\mu^{(l)}\triangleright
 \lambda^{(l)}}r_{\aaa_1\bbb_1}x_{\aaa_1\bbb_1}&
 \end{array}
 \end{equation}
 for some $r_{\aaa_1\bbb_1},a_\V\in R$ where the tableaux $
 \aaa_1,\bbb_1\in\std(\mu^{(l)}) $ involve the numbers permuted by $
 \Si_{\alpha_l}$. For \textcolor{black}{$ \aaa_1,\bbb_1 $} and $ \V$ appearing in the sum set
 $\pmb{\aaa}:=(\T^{\lambda^{(1)}},\ldots ,\aaa_1,\ldots, \T^{\lambda^{(r)}})$,
 $\pmb{\bbb}:=(\T^{(1)},\ldots,\bbb_1, \ldots, \T^{(r)}) $ and
 $\pmb{\ccc}:=(\T^{\lambda^{(1)}},\ldots,\V,\ldots, \T^{\lambda^{(r)}})$.
 Then $ \pmb{\ccc} \in \std(\blambda) $ and $ \pmb{\aaa}, \pmb{\bbb} \in \std(\bmu) $ where
 $\bmu:=(\lambda^{(1)},\ldots ,\mu^{(l)},\ldots\lambda^{(r)})$.
 Moreover, by our definition of the dominance order we have
 $\bmu\triangleright \blambda$, $\pmb{\ccc}\triangleright \bT$ and so $m_{\pmb{\aaa}
\pmb{ \bbb}} \in \Y^{\blambda}$. On the other hand, we have
 $$ U_{\blambda}\textcolor{black}{\Elambda} x_{\lambda^{(1)}\T^{(1)}}\cdots
 x_{\aaa_1\bbb_1} \cdots x_{\lambda^{(r)}\T^{(r)}}
 =U_{\blambda}\textcolor{black}{\Elambda}g_{d(\aaa_1)}^{\ast}x_{\lambda^{(1)}\T^{(1)}}\cdots
 x_{\mu^{(l)}} \cdots
 x_{\lambda^{(r)}\T^{(r)}}g_{d(\bbb_1)} = m_{\pmb{\aaa} \pmb{\bbb}}\nonumber
 $$
and similarly for $ m_{\textcolor{black}{{\blambda}}\bT} $ and $ m_{\textcolor{black}{{\blambda}}\pmb{\ccc}}$.
Using Lemma {\ref{lemmasubalgebra}} in the other direction
together with ${\rm
 res}_{\T^{(l)}}(k)={\rm
 res}_{\bT}(k)$ we then get 
 $$\begin{array}{rcl}m_{\textcolor{black}{{\blambda}}\bT}J_k&=&c_{\bT}(k)m_{\textcolor{black}{\blambda}\bT}+\displaystyle\mathop{\sum_{\pmb{\ccc}\in\std(\blambda)}}_{\pmb{\ccc}\triangleright
 \T}a_{\pmb{\ccc}}m_{\textcolor{black}{{\blambda}}\pmb{\ccc}}
\mod  \Y^{\blambda}
\end{array}$$
which shows the Proposition for $ \bT $ of the initial kind.
\color{black}

\medskip
 For $ \bT $ a general multitableau, there exists a multitableau
 $\bT_0$ of the initial kind together with a distinguished right coset
 representative $w_{\bT}$ of $\Si_\alpha$ in $\Si_n$ such that $\bT=\bT_0 w_{\bT}$.
 Let $w_{\bT}= s_{i_1} s_{i_2} \ldots s_{i_k} $ be a reduced expression
 for $w_{\bT} $. Then we have that $ i_j $ and $ i_j +1 $ are located in
 different blocks of $ \bT_0  s_{i_1} \ldots  s_{i_{j-1}}   $ for all
 $ j \geq 1 $ and that $ \bT_0  s_{i_1} \ldots  s_{i_{j-1}} s_{i_{j}} $ is obtained from $ \bT_1  s_{i_1} \ldots  s_{i_{j-1}}   $ by
 interchanging $ i_j $ and $ i_j +1 $. Using Lemma \ref{lem0} and (4)
 of Lemma \ref{J1} we now get that
 $$m_{\textcolor{black}{{\blambda}}\bT}J_k=
 m_{\textcolor{black}{{\blambda}}\bT_0}g_{w_{\bT}}J_k= m_{\textcolor{black}{{\blambda}}\bT_0}J_{kw_{\bT}^{-1}}g_{w_{\bT}}.$$
 Since $ \bT_0 $ is of the initial kind, we get
 $$\begin{array}{ll  }m_{\textcolor{black}{{\blambda}}\bT}J_k & =
  m_{\textcolor{black}{{\blambda}}\bT_0}J_{kw_{\T}^{-1}}g_{w_{\bT}}
 =
 \left(c_{\bT_0}(kw_{\bT}^{-1})m_{\textcolor{black}{\bT^{\blambda}}\bT_0}+\displaystyle\mathop{\sum_{\Bv_0\in\std(\lambda)}}_{\Bv_0\triangleright
 \T_0}a_{\Bv_0}m_{\textcolor{black}{{\blambda}}\Bv_0}\right)g_{w_{\bT}} \\
  &=  c_{\bT}(k)m_{\textcolor{black}{{\blambda}}\bT}+\displaystyle\mathop{\sum_{\Bv\in\std(\lambda)}}_{\Bv\triangleright
 \T}a_{\Bv}m_{\textcolor{black}{{\blambda}}\Bv}
 \end{array}$$
 where we used that the occurring $ \Bv_0 $ are all of the initial kind
 such that $m_{\Bv}=m_{\Bv_0}g_{w_{\bT}}$ with $ \Bv\triangleright \bT_0$ and
 $a_{\Bv}=a_{\Bv_0}$. This finishes the proof of the Proposition.
 \end{demo}

\medskip
In view of the Proposition, we can now apply the general theory developed in \cite{Mat2}. In particular, we recover the semisimplicity criterion
of Chlouveraki and Poulain d'Andecy, \cite{MCH}, and can even generalize it to the case of ground fields of positive characteristic.
We leave the details to the reader.

\section{Representation theory of the algebra of braids and ties.}
In this final section we once again turn our attention to the
algebra $ \E $ of braids and ties.

In the paper \cite{Rh}, the representation theory of $ \E $ was studied in the generic case, where a parametrizing set for the irreducible
modules was found. On the other hand, the dimensions of the generically irreducible modules were not determined in that paper.
In this section we show that $\E $ is a cellular algebra by giving a concrete combinatorial construction of a cellular basis for it.
As a bonus we obtain a closed formula for the dimensions of the cell modules, which in particular gives a formula for the
irreducible modules in the generic case. Although the construction of the cellular basis for $ \E $ follows the outline of the
construction of the cellular
basis $ \B$ for $ \Y(q)$, the combinatorial details are quite a lot more involved and, as we shall see, involve a couple of new ideas.

\medskip
\textcolor{black}{
It should be pointed out that Jacon and Poulain d'Andecy have recently obtained a very elegant classification of the irreducible
modules for $ \E $ via Clifford theory, see
\cite{JdA1}. Their approach relies on the connection with the Yokonuma-Hecke algebra and
therefore does not work for all fields. Our cellular algebra approach works, at least in principle, for all fields.}

\medskip
Recall that $ \E $ is the $ S:= {\mathbb Z}[q,q^{-1}]$-algebra defined by the generators and relations given in Definition
\ref{braidsties} and
that it was shown in \cite{Rh} that $ \E $ has an $S$-basis of the form $ \{ E_{A} g_w | A \in \mathcal{SP}_n, w \in \Si_n \} $.
(Theorem \ref{embeddding} gave a new proof of this basis).
We shall often need the following relations in $ \E$, that have already appeared implicitly above
\begin{equation}{\label{wehavetherelations}}
E_A g_w = g_w E_{ A w}  \mbox{ and } E_A E_B = E_C \mbox{  for } w\in \Si_n , A,B  \in \mathcal{SP}_n
\end{equation}
where $ C \in \mathcal{SP}_n $ is minimal with respect to $ A \subseteq C, B \subseteq C$.

\subsection{Decomposition of $\E$}
In this subsection we obtain central idempotents of $ \E$ and a corresponding subalgebra decomposition of $ \E$.
This is inspired by I. Marin's recent paper \cite{M}, which in turn is inspired by  \cite{Sol} and \cite{Gre}.

Recall that for a finite poset $ (\Gamma,\preceq)$ there is an associated M\"{o}ebius function $ \mu_{\Gamma}: \Gamma \times \Gamma \rightarrow \mathbb Z$.
In our set partition case $ ({\cal SP}_n, \subseteq)$ the M\"{o}ebius function $  \mu_{ {\cal SP}_n}  $  is given by the formula
\begin{equation}
\mu_{{\cal SP}_n}(A,B)=  \left\{ \begin{array}{ll}
(-1)^{r-s}\prod_{i=1}^{r-1}(i!)^{r_{i+1}}&  \mbox{if }  A \subseteq   B \\ 0 &  \mbox{otherwise  } \end{array} \right.
\end{equation}
where $r$ and $s$ are the number of blocks of $A$ and $B$ respectively, and
where $r_i$ is the number of blocks of $B$ containing exactly $i$ blocks of $A$.

We use the M\"{o}ebius function $ \mu =  \mu_{ {\cal SP}_n}   $ to introduce
a set of orthogonal idempotents elements of $\E$.
This is a special case of the general construction given in {\it loc. cit.}
For $A\in \mathcal{SP}_n$ the idempotent $\mathbb{E}_A \in \E$ \nomenclature[19]{$\mathbb{E}_A$}{The idempotent associated with the set partition $A$} is given by the formula
\begin{equation}{\label{moebius}}
\mathbb{E}_A:= \sum_{ A \subseteq B }\mu(A,B)E_B.
\end{equation}
For example, we have
$\be_{\{\{1\},\{2\},\{3\}\}}=E_{\{\{1\},\{2\},\{3\}\}} -E_{\{\{1, 2\},\{3\}\}}
-E_{\{\{1\},\{2, 3\}\}} -
E_{\{\{1, 3\},\{2\}\}}+2E_{\{\{1,2,3\}\}}$.
We have the following result.
\begin{propos}\label{EGpropiedades} The following properties hold.
\begin{itemize}
\setlength\itemsep{-1.3em}
\item[(1)] $ \{ \mathbb{E}_A | A \in {\cal SP}_n \} $ is a set of orthogonal idempotents of $ \E$. \newline
\item[(2)] For all $w\in\Si_n$ and $A\in \mathcal{SP}_n$ we have $\be_A g_w=g_w\be_{Aw}$. \newline
\item[(3)] For all $A\in \mathcal{SP}_n$ we have $\be_{A}E_B=\left\{\begin{array}{ll}\be_{A}&\mbox{ if } B\subseteq A\\
0&\mbox{ if } B\not\subseteq A.  \end{array}\right.$
\end{itemize}
\end{propos}

\begin{demo}
We have already mentioned (1) so let us prove (2).
We first note that the order relation $ \subseteq $ on $ \mathcal{SP}_n$ is compatible with the action of $\Si_n$ on $\mathcal{SP}_n$ that
is $A\subseteq B$ if and only if $Aw\subseteq Bw$ for all $w\in \Si_n$. This implies that $\mu(Aw,Bw)=\mu(A,B)$ for all $w\in\Si_n$. From this we get, via
(\ref{wehavetherelations}), that
$$\be_A g_w=g_w\sum_{ A \subseteq B }\mu(A,B)E_{Bw}=g_w\sum_{ A \subseteq Cw^{-1}} \mu(A,Cw^{-1})E_{C}=g_w\sum_{Aw \subseteq C }\mu(Aw,C)E_{C}=g_w\be_{Aw}$$
showing (2). Finally, we obtain (3) from the orthogonality of the $\be_A$'s and the formula \textcolor{black}{$E_B = \sum_{B\subseteq A}   \be_{A}$}
which is obtained by inverting ({\ref{moebius}}) \textcolor{black}{(see also \cite{Gre})}.
\end{demo}

\medskip
We say that a set partition
\textcolor{black}{$A=\{I_{1},\ldots,I_{k}\}$ of $\mathbf{n}$ is of}
type $\alpha\in \Par $ if there exists a permutation $\sigma$ such
that $(|I_{i_{1\sigma}}|,\ldots,|I_{i_{k\sigma}}|)=\alpha$. For
example, the set partitions of $\mathbf{3}$ of type $(2,1)$ are
$\{\{1,2\},\{3\}\}$, $\{\{1,3\},\{2\}\}$ and $\{\{2,3\},\{1\}\}$.
For short, we write $|A|=\alpha$\nomenclature[20]{$ \lvert  A \rvert $}{Type of the set partition $A$} if $A\in \mathcal{SP}_n$ is of type $\alpha$.
\textcolor{black}{We also say that a multicomposition
  $\blambda=(\lambda^{(1)},\ldots,\lambda^{(m)}) \in \MC$ is of type $ \alpha$
  if the associated set partition $ A_{\blambda}$ is of type $ \alpha$.}

\medskip
For each $\alpha\in \Par$  we define the following element $\be_{\alpha}$ of  \nomenclature[21]{$\be_{\alpha}$}{The central orthogonal idempotent associated with the partition $\alpha$} {\color{black}$\E$.}
\begin{equation}{\label{Ealpha}}
\be_{\alpha}:=\sum_{\textcolor{black}{A\in \mathcal{SP}_n},\,  |A|=\alpha}\be_{A}.
\end{equation}
Then by Proposition \ref{EGpropiedades}
we have that $\{\be_\alpha |  \alpha \in \Par \}$ is a set of central
orthogonal idempotents of $\E$, which is complete: $ \sum_{ \alpha \in \Par} \be_\alpha = 1 $.
As an immediate consequence we get the following decomposition of $\E$ into a direct sum of two-sided ideals
\begin{equation}\label{descE}
\E=\bigoplus_{\alpha\in\Par} \Ea
\end{equation}
where we define $\Ea:=\be_{\alpha}\E$\nomenclature[22]{$\Ea$}{The
$S$-subalgebra of $\E$ associated with the partition $\alpha$}.
\textcolor{black}{Each $ \Ea $ is an $S$-algebra with identity $
\be_{\alpha} $}.

\medskip
Using the $\{E_Ag_w\}$-basis for $\E$, together with part (3) of Proposition \ref{EGpropiedades}, we get that the set
\begin{equation}\label{basessubalgebras}
 \{\be_{A}g_w\mid w\in\Si_n,\,|A|=\alpha\}
\end{equation}
is an $S$-basis for $\Ea$. In particular, we have that the dimension
of $\Ea$ is $b_n(\alpha)n!$, where $b_n(\alpha)$ is the number of
set partitions of $\mathbf{n}$ having type $\alpha\in\Par$. The
numbers $b_n(\alpha)$ \nomenclature[23]{$b_n(\alpha)$}{The Fa\`{a}
di Bruno coefficients associated with the partition $\alpha$} are
the socalled Fa\`{a} di Bruno coefficients and are given by the
following formula
\begin{equation}\label{faadibruno}
b_n(\alpha)=\dfrac{n!}{(k_1!)^{m_1}m_1!\cdots (k_r!)^{m_r}m_r! }
\end{equation}
where $\alpha=(k_1^{m_1},\ldots,k_r^{m_r})$ and $k_1>\ldots>k_r$.

\subsection{Cellular basis for $\E$}\label{subseccelbtalgebra}
Let us explain the ingredients of our cellular basis for $\E$. The antiautomorphism $ \ast$ is easy to explain, since
one easily checks on the relations for $ \E$ that $\E$ is endowed with an $S$-linear antiautomorphism $ \ast $, satisfying
$ e_i^{\ast} := e_i  $ and $ g_i^{\ast} := g_i $. We have that $ \be_A^{\ast} = \be_A$.

\medskip
Next we explain the poset denoted $ \Lambda $ in Definition
\ref{cellular} of cellular algebras. By general principles, it
should be the parametrizing set for the irreducible modules for $
\E$ in the generic situation, so let us therefore recall this set $
{\mathcal L}_n $ from \cite{Rh}. $ {\mathcal L}_n $
\nomenclature[24]{$ {\mathcal L}_n $}{The parametrizing set for the
irreducible modules for $ \E$} is the set of pairs $ \Lambda =
(\blambda \mid \bmu) $ where $ \blambda = ( \lambda^{(1)}, \ldots ,
\lambda^{(m)})$ is an $m$-multipartition of $ n $. We require that $
\blambda $ be \textit{increasing} by which we mean that $
\lambda^{(i)} < \lambda^{(j)} $ only if $ i < j $ where  $ < $ is
any fixed extension of the usual dominance order on partitions to a
total order, and where we set $ \lambda < \textcolor{black}{\tau} $ if
$ \lambda $ and $ \tau $ are partitions such that $ | \lambda | < |
\tau |$.

\medskip
In order to describe the $ \bmu $-ingredient of $ \Lambda $ we need to introduce some more notation.
The multiplicities of equal $ \lambda^{(i)}$'s give rise to a composition of $ m$.
To be more precise, let
$ m_1 $ be the maximal $ i $ such that $   \lambda^{(1)}=\lambda^{(2)}=\ldots =\lambda^{(i)}   $,
let $ m_2 $ be the maximal $ i $ such that
$   \lambda^{(m_1+1)}  =   \lambda^{(m_1+2)}  = \ldots =   \lambda^{(m_1 +i)}   $, and so on
until $ m_q$. Then we have that $ m = m_1 + \ldots + m_q $.
We then require that $ \bmu $ be
of the form $ \bmu = (\mu^{(1)}, \ldots,  \mu^{(q)}) $ where each $ \mu^{(i)} $ is partition of $m_i$.
\textcolor{black}{This is the description of $ {\mathcal L}_n $ as a set, as given in \cite{Rh}.
If $ \alpha \in \Par $ we use the notation
\begin{equation}
  { \mathcal L}_n (\alpha) := \{ (\blambda \mid  \bmu)  \in { \mathcal L}_n   \big|      {\blambda} \, \,  \mbox{is of type } \, \, \alpha \}.
\end{equation}  }
\nomenclature[25]{${ \mathcal L}_n (\alpha)$}{The parametrizing set
for the
irreducible modules for $\Ea$}\textcolor{black}{We now introduce a poset structure on $ {\mathcal L}_n $.} 
Suppose that $ \Lambda =(\blambda \mid \bmu) $ and $
\overline{\Lambda} =(\overline{\blambda}\mid \overline{\bmu} ) $ are
elements of $  {\mathcal L}_n $ such that $ \norm{ \blambda } = \norm{
\overline{\blambda}} $.
We first write $ \blambda \rhd_1 \overline{\blambda} $ if
$ \blambda = (\lambda^{(1)}, \ldots, \lambda^{(m)}) $ and $ \overline{\blambda} = (\overline{\lambda^{(1)}}, \ldots, \overline{\lambda^{(m)}}) $
and if there exists a permutation $ \sigma $ such that
$ (\lambda^{ (1 \sigma) }, \ldots, \lambda^{ (m\sigma )}) \rhd (\overline{\lambda^{(1)}}, \ldots, \overline{\lambda^{(m)}})$
where $\rhd$ is the dominance order on $m$-multi\-partitions, introduced above.
We then say that $ \Lambda \rhd \overline{\Lambda} $ if $ \blambda \rhd_1 \overline{\blambda} $ or if $ \blambda =\overline{\blambda} $ and
\textcolor{black}{$ \bmu \rhd \overline{\bmu}$.}
As usual we set $ \Lambda \unrhd \overline{\Lambda} $ if $ \Lambda \rhd \overline{\Lambda} $ or if
$ \Lambda = \overline{\Lambda} $.
This is our description of $ {\mathcal L}_n $ as a poset. Note that if $ \norm{ \blambda } \neq \norm{ \overline{\blambda}} $ then
$ \Lambda $ and $ \overline{\Lambda} $ are by definition not comparable.

\begin{obs}
We could have introduced an order '$\succ$' on $ {\mathcal L}_n $ by
replacing '$\rhd_1$' by '$\rhd$' in the above definition, that is $
\Lambda \succ \overline{\Lambda} $ if $ \blambda \rhd
\overline{\blambda} $ or if $ \blambda =\overline{\blambda} $ and
\textcolor{black}{$ \bmu \rhd \overline{\bmu}$.} Then '$ \succ$' is a
finer order than '$ \rhd$', but in general they are different. The
reason why we need to work with '$ \rhd$' rather than '$ \succ$'
comes from the straightening procedure of Lemma
\ref{Lemmatableaustandar} below.

We could also have introduced an order on $ {\mathcal L}_n $ by replacing '$  =   $'  with '$  =_1 $' in the above definition, where '$  =_1 $'  is
defined via a permutation $ \sigma $, similar to what we did for $ \rhd_1 $:
that is $ \Lambda \succ \overline{\Lambda} $ if $ \blambda \rhd_{\textcolor{black}{1}} \overline{\blambda} $ or if $ \blambda =_1\overline{\blambda} $ and
\textcolor{black}{$ \bmu \rhd \overline{\bmu}$.}
On the other hand, since $ \blambda $ and $ \overline{\blambda} $ are assumed to be increasing multipartitions, we get that
'$ =_1 $' is just usual equality '$=$' and hence we would get the same order on $ {\mathcal L}_n $.
\end{obs}

Let us give an example to illustrate our order.

\begin{exa}
We first note that $(3,3,1)\rhd(3,2,2)$ in the dominance order on partitions, but both are incomparable with the partition $(4,1,1,1)$.
Suppose now that $(3,2,2) < (4,1,1,1) < (3,3,1)$ in our extension of the dominance order. We then consider the following increasing multipartitions of 25
$$\blambda=((2),(2), (3,2,2), (4,1,1,1), (3,3,1) )\;\mbox{ and }\;\overline{\blambda}=((2),(2), (3,2,2), (3,2,2) , (4,1,1,1)  ).$$
Then we have that $ \blambda $ and $ \overline{\blambda}$ are increasing multipartitions, but
incomparable in the dominance order on multipartitions. On the other hand
$\blambda \rhd_1 \overline{\blambda}$ via the permutation $\sigma=s_4$ and hence we have the following relation in $ {\mathcal L}_n $
$$\Lambda:=\bigg(\blambda \bigm| \big((2),(1),(1),(1) \big) \bigg) \rhd
\bigg(\overline{\blambda} \bigm| \bigg(\big((1^2),(2), (1)  \big)\bigg)=:\overline{\Lambda}.$$
\end{exa}

For $ \Lambda=(\blambda \mid  \bmu) \in { \mathcal L}_n$ as above,
we next define the concept of $ \Lambda $-tableaux. Suppose that $ \et $ is a pair $\et=(\bT \mid \bu) $.
Then $ \et $ is called a $ \Lambda $-tableau if $ \bT = (\T^{(1)}, \ldots,  \T^{(m)}) $ is a
multitableau of $ n $ in the usual sense, satisfying $ Shape(\bT) = \blambda $, and
$ \bu $ is a \textcolor{black}{$\bmu$-multitableau of the initial kind.}
As usual, if $ \et $ is $ \Lambda $-tableau we define $ Shape(\et) := \Lambda$.

Let $ \Tab(\Lambda) $\nomenclature[26]{$ \Tab(\Lambda) $}{The set of
$ \Lambda$-tableaux} denote the set of $ \Lambda$-tableaux and let $
\Tab_n:= \cup_{\Lambda \in {\mathcal L}_n} \Tab(\Lambda) $. We then
say that $ \et = (\bT \mid \bu) \in \Tab(\Lambda) $ is row standard
if its ingredients are row standard multitableaux in the usual
sense.

We say that $ \et = (\bT \mid \bu) \in \Tab(\Lambda) $ is standard if its
ingredients are standard multitableaux and if moreover $\bT$ is an \textit{increasing}
multitableau.
By increasing we here mean
that whenever $ \lambda^{(i)} = \lambda^{(j)} $ we have that $ i < j $ if and only if
$ \min ( \T^{(i)}) < \min (\T^{(j)}) $ where $ \min (t) $ is the function that reads off the minimal entry of the tableau $ t $.
We define $ {\rm Std}(\Lambda) $ to be the set of all standard $\Lambda$-tableaux.

\begin{exa} For $\Lambda=\bigg(\big( (1,1),(2), (2),  (2,1)\big)\,\big| \,\big((1), (1,1),(1)   \big)\bigg)$ we consider the following $ \Lambda$-tableaux
\begin{equation}\begin{array}{c}\et_1 :={\footnotesize\left(\;\left(\;\young(1,9)\;,\young(35)\;,\young(68)\;,\young(24,7)\;\right)  \bigm|
\left(  \young(1) \, , \young(2,3) \, , \,\young(4)\;\right)\;\right)}\\[4mm]
\et_2 :={\footnotesize\left(\;\left(\;\young(1,9)\;,\young(56)\;,\young(38)\;,\young(24,7)\;\right) \big|
\left(  \young(1) \, , \, \young(2,3)\; ,\young(4) \right)\;\right)}.\end{array}\label{equa1}\end{equation}
Then by our definition, $\et_1$ is a standard $\Lambda$-tableau, but $\et_2$ is not.
\end{exa}

\begin{obs}
The use of the function $ \min ( \cdot) $ is somewhat arbitrary. In fact we could have used any injective function
with values in a totally ordered set.
\end{obs}

\medskip
For $ \bT = (\T^{(1)}, \ldots, \T^{(m)}) $ and $ \overline{\bT} = (\overline{\T^{(1)}}, \ldots, \overline{t^{(m)}}) $
we define
$ \bT  \rhd_1 \overline{\bT} $ if there exists a permutation $ \sigma $ such that
$ (\T^{(1 \sigma )}, \ldots, \T^{(m \sigma) }) \rhd (\overline{\T^{(1)}}, \ldots, \overline{\T^{(m)}})$
in the sense of multitableaux.
We then extend the order on  $ {\mathcal L}_n $ to $ \Tab_n $ as follows.
Suppose that
$ \et =(\bT \mid \bu) \in  \Tab(\Lambda)$ and $  \overline{\et} = (\overline{\bT} \mid \overline{\bu} ) \in \Tab(\overline{\Lambda}) $
and that $ \Lambda \unrhd \overline{\Lambda}$.
Then we say that $ \bT \rhd \overline{\bT}  $ if
$ \bT \rhd_1 \overline{\bT}  $ or if
$ \bT = \overline{\bT}  $ and
\textcolor{black}{$ \bu \rhd \overline{\bu}$.}
As usual we set $ \et \unrhd \overline{\et} $ if $ \et \rhd \overline{\et} $ or
$ \et = \overline{\et} $. This finishes our description of $\Lambda$-tableaux as a poset.

\medskip
From the basis of $ \E $ mentioned above, we have that $\dim {\mathcal E}_n(q) = b_n n! $ where $ b_n $ is the $n$'th Bell number, that is
the number of set partitions on $ \mathbf{n}$.
Our next Lemma is a first strong indication of the relationship between our notion of standard tableaux
and the representation theory of $ {\mathcal E}_n(q) $.

Recall the notation $ d_{\lambda} := | {\rm Std}(\lambda)| $ that we introduced for partitions $ \lambda $.
In the proof of the Lemma, and later on, we shall use repeatedly the formula $ \sum_{ \lambda \in \Par} d_{\lambda}^2 = n! $.

\begin{lem}{\label{cardinal}}
With the above notation we have that
$ \sum_{ \Lambda \in { \mathcal L}_n } | {\rm Std}(\Lambda) |^2 = b_n n! $.
\end{lem}

\begin{demo}
It is enough to prove the formula
\begin{equation}{\label{toprovetheformula}}
 \sum_{ \Lambda \in { \mathcal L}_n (\alpha)} | {\rm Std}(\Lambda) |^2 = b_n(\alpha) n!
\end{equation}
where $ b_n(\alpha) $ is the Fa\`{a} di Bruno coefficient introduced above.
Let us first consider the case $ \alpha = (k^m) $, that is $ n = mk$. Then we have
$$ b_n(m,k):=b_n(\alpha)=\dfrac{1}{m!} \binom{n}{k \cdots k} $$
with  $ k $ appearing $ m $ times in the multinomial coefficient.
Let $ \{ \lambda^{(1)}, \lambda^{(2)}, \ldots, \lambda^{(d)} \} $ be the fixed ordered enumeration of all the partitions of $ k$, \textcolor{black}{introduced above}.
If $  \Lambda= (\blambda \mid \bmu) \in { \mathcal L}_n (\alpha) $
then $ \blambda $ has the form
$$ \blambda= (\overbrace{\lambda^{(1)}, \ldots ,\lambda^{(1)},}^{m_1}\overbrace{\lambda^{(2)}, \ldots ,\lambda^{(2)},}^{m_2} \ldots,
\overbrace{\lambda^{(d)} \ldots ,\lambda^{(d)}}^{m_d})$$
where the $ m_i $'s are non-negative integers with sum $ m$ and
$ \bmu= (\mu^{(1)}, \mu^{(2)}, \ldots, \mu^{(d)} ) $ is a multipartition of type $ \norm{\bmu } = (m_1, m_2, \ldots, m_d)$.
The number of increasing multitableaux of shape $ \blambda $ is
$$
\frac{1}{m_1! \ldots m_d!}
 \binom{n}{  k  \cdots k}
 \prod_{j=1}^d   {d_{\lambda^{(j)}}^{ m_j}}
$$
whereas the number of standard tableaux of shape $ \bmu $ is
$
\prod_{j=1}^d d_{\mu^{(j)}}
$
and so we get
\begin{equation}{\label{whereasthe}}
  |{\rm Std}(\Lambda)|  = \frac{1}{m!}
 \binom{m}{ m_1   \cdots m_d}
 \binom{ n}{k   \cdots k}
 \prod_{j=1}^d   {d_{\lambda^{(j)}}^{ m_j}} d_{\mu^{(j)}}
\end{equation}


\medskip
By first fixing $ \blambda $ and then letting each $ \mu^{(i)} $ vary over all possibilities we get that the square
sum of the above $ |{\rm Std}(\Lambda)|  $'s is the sum of
$$
\binom{n}{k \cdots k}^2
 \prod_{j=1}^d   \frac{d_{\lambda^{(j)}}^{2m_j}}{m_j !} =
\binom{n}{k \cdots k}^2
\frac{1}{m!}
 \binom{m}{ m_1   \cdots m_d}
 \prod_{j=1}^d   {d_{\lambda^{(j)}}^{2m_j}}
$$
with the $ m_i$'s running over the above mentioned set of numbers.
But by the multinomial formula, this sum is equal to
 $$
 \binom{n}{k \cdots k}^2 \frac{1}{m!} \left(  \sum_{j=1}^d   d_{\lambda^{(j)}}^{2} \right)^m
 = \binom{n}{k \cdots k}^2 \frac{1}{m!} \,  k!^m  =
 \dfrac{n!}{m!} \binom{n}{k \cdots k} = b_n(\alpha) n!
 $$
and ({\ref{toprovetheformula}}) is proved in this case.

\medskip
Let us now consider the general case where $ \alpha=(k_1^{M_1}, \ldots, k_r^{M_r}) $,
where $ k_1 > \cdots > k_r $. Set $ n_i = k_i M_i$, $ M := M_1 + \ldots + M_r$. Then $n= n_1 + \ldots + n_r$
and the Fa\`{a} di Bruno coefficient $ b_n(\alpha) $ is given by the formula
\begin{equation}{\label{isgivenbtheformula}}
b_n(\alpha) = \binom{n}{ n_1 \cdots n_r} b_{n_1}(M_1, k_1) \cdots b_{n_r}(M_r, k_r).
\end{equation}
Let us now consider the square sum $ \sum_{ \Lambda \in { \mathcal L}_n(\alpha)} | {\rm Std}(\Lambda) |^2$.
For $ \Lambda= (\blambda \mid  \bmu) \in { \mathcal L}_n(\alpha) $
we split $ \blambda $ into multipartitions $ \blambda_1, \ldots, \blambda_r$, where $ \blambda_1= (\lambda^{ (1)},\ldots,  \lambda^{ (M_1)})$,
$ \blambda_2= (\lambda^{ (M_1+1)},\ldots,  \lambda^{ (M_1+M_2)})$, and so on.
We split $ \bmu $ correspondingly into $ \bmu_i$'s and
set $ \Lambda_i := (\blambda_i \mid \bmu_i) $. Then $ \Lambda_i \in { \mathcal L}_{n_i}( (k_i^{M_i}))$
and we have

\begin{equation}{\label{andmoreover}}
 | {\rm Std}(\Lambda) |=\binom{n}{n_1 \cdots n_r}   | {\rm Std}_{n_1}(\Lambda_1)| \cdots | {\rm Std}_{n_r}(\Lambda_r)|
\end{equation}
where $ {\rm Std}_{n_i}( \Lambda_i) $ means standard tableaux of shape $ \Lambda_i $ on $ {\bf n_i}$.
Combining ({\ref{toprovetheformula}}), ({\ref{isgivenbtheformula}}) and ({\ref{andmoreover}}) we get that
$$ \sum_{ \Lambda \in { \mathcal L}_n(\alpha)} | {\rm Std}(\Lambda) |^2 = n! b_n(\alpha) $$
as claimed.

\end{demo}

\begin{coro}{\label{spelling}}
Suppose that $ \Lambda= (\blambda \mid \bmu) \in {\mathcal L}_n$ is above
with $ \blambda =  (\lambda^{(1)}, \ldots, \lambda^{(m)} )  $ and
$ \bmu =  (\mu^{(1)}, \ldots, \mu^{(q)} )  $ and set $ n_i := | \lambda^{(i)} | $ and $ m_i := | \mu^{(i)} | $.
Then we have that
$$ | \std(\Lambda) | =
 \frac{1}{ m_1!   \cdots m_q!}
 \binom{ n}{ n_1   \cdots n_m}
 \prod_{j=1}^m   {d_{\lambda^{(j)}}}  \prod_{j=1}^q d_{\mu^{(j)}}.
$$
\end{coro}

\begin{demo}
This follows by combining ({\ref{whereasthe}}) and ({\ref{andmoreover}}) from the proof of the Lemma.
\end{demo}

 \medskip

 \textcolor{black}{We fix the following combinatorial notation.
   Let $ \Lambda=(\blambda \mid
\textcolor{black}{\bmu}) =
(( \lambda^{ (1 ) }, \ldots,  \lambda^{ (m ) }) \mid ( \mu^{ (1 ) }, \ldots,  \mu^{ (q ) }))
\in {\cal L}_n(\alpha) $}. With $ \Lambda $
we have associated the set
of multiplicities $ \{ m_{i} \}{\textcolor{black}{_{i=1, \ldots, q}}}$ of equal $ \lambda^{(i)}$'s. We
now also associate with $ \Lambda $ the set of multiplicities $ \{ k_i
\}{\textcolor{black}{_{i=1, \ldots, r}}}$ of equal block sizes $  | \lambda^{(i)} |$.
That is, $ k_1 $ is the maximal $ i $ such that $   |\lambda^{(1)}|
=   | \lambda^{(2)} | = \ldots =  | \lambda^{(i)} |  $, whereas $
k_2 $ is the maximal $ i $ such that $  | \lambda^{(k_1+1)} |  =  |
\lambda^{(k_1+2)} | = \ldots =  | \lambda^{(k_1 +i)} |  $ \textcolor{black}{and so on}.
\textcolor{black}{
We can also describe the $ k_i$'s in terms of the type of $ \blambda $, that is $ \alpha $: indeed we have
$ \alpha = (a_r^{k_r}, \ldots, a_1^{k_1})$
with $ a_r > a_{r-1} > \ldots > a_1$: recall that $ \blambda$ is increasing.
Note that $ m_1+m_2+\cdots + m_q = k_1 +k_2 + \cdots + k_r=m $ and that $ |\mu^{ (j) }| = m_j $
for all $ j $.
}

Let $\Si_{\Lambda} \leq \Si_{n}$
\nomenclature[27]{$\Si_{{\Lambda}}$}{The stabilizer subgroup of the
set partition $A_{\blambda}$, where $\Lambda=(\blambda\mid -)$} be the stabilizer subgroup of the set
partition $ A_{\blambda} \textcolor{black}{=\{I_1, I_2, \ldots,
I_m\}} $ that was introduced in (\ref{fillinginthenumbers}). Then
the two sets of multiplicities give rise to subgroups $
\Si^k_{\Lambda} $ and $ \Si^m_{\Lambda} $ of $ \Si_{\Lambda}$ where
$ \Si^k_{\Lambda} $ consists of the order preserving permutations of
the equally sized blocks of $ A_{\blambda} $, whereas
$\Si^m_{\Lambda}$  consists of the order preserving permutations of
those blocks of $ A_{\blambda} $ that correspond to equal $
\lambda^{(i)}$'s. Clearly we have $ \Si^m_{\Lambda} \leq
\Si^k_{\Lambda}\leq \Si_{{\Lambda}}$.

\color{black} We observe that $\Si^k_{\Lambda}$
\nomenclature[28]{$\Si^k_{\Lambda}$}{The subgroup of $\Si_{\Lambda}$
of the order preserving permutations of the equally sized blocks of
$ A_{\blambda}$} and $\Si^m_{\Lambda}$
\nomenclature[29]{$\Si^m_{\Lambda}$}{The subgroup of $\Si_{\Lambda}$
of the order preserving permutations of those blocks of $
A_{\blambda} $ that correspond to equal $ \lambda^{(i)}$'s} are
products of symmetric groups,
\begin{equation}{\label{areproductsofsymmetricgroup}}
  \Si^k_{\Lambda} \cong \Si_{k_1} \times \ldots \times \Si_{k_r }, \,\, \, \, \,\, \, \, \,\, \, \,
  \Si^m_{\Lambda} \cong \Si_{m_1} \times \ldots \times \Si_{m_q }
\end{equation}
and in fact $\Si^k_{\Lambda}$ is a Coxeter group on generators $ B_i $ that we explain shortly, and
$ \Si^m_{\Lambda} $ is a parabolic subgroup of $\Si^k_{\Lambda}$.
Define subsets $ S_{\Lambda}^m \subseteq S_{\Lambda}^k$ of $ \mathbf{m} $ via
\begin{equation}
S_{\Lambda}^k:= \{ i \in \mathbf{m} | i \neq k_1 + \ldots + k_j \mbox{ for all } j \}, \, \,  \, \,  \, \,
S_{\Lambda}^m:= \{ i \in \mathbf{m} | i \neq m_1 + \ldots + m_j \mbox{ for all } j \}.
\end{equation}
Then for $ i \in S_{\Lambda}^k $ the generator
$ B_i$ of $\Si^k_{\Lambda}$ is
the minimal length element of $\Si_n$
that interchanges the two consecutive blocks $ I_i $ and $ I_{i+1} $ of
$ A_{\blambda} $ (of equal size). Moreover,
$ B_i $ is also a generator for $\Si^m_{\Lambda}$ if and only if
$ i \in \Si^m_{\Lambda}$.
Let us describe $ B_i$ concretely. Letting $ a := | I_i |$ we can
write
\begin{equation} I_i = \{ c +1,c+2, \ldots, c+a \} \mbox{     and    }
  I_{i+1}=\{ c +a+1,c+a+2, \ldots, c+2a \}
  \end{equation}
for some $ c $.
With this notation we have
\begin{equation}{\label{Bi}}
B_i = (c+1,c+a+1) (c+2, c+a+2) \cdots (c+a, c+2a).
\end{equation}
For $ i >  j  $ we set $ s_{ij} := s_{i+c} s_{i-1+c} \ldots  s_{j+c} $ and can then write $ B_i $ in terms of
the $ s_{ij} $'s, and therefore in terms of
simple transpositions $s_i$, as follows
\begin{equation}{\label{expansiontrans}}
B_i = s_{a,1} s_{ a+1,2} \ldots  s_{2a-1,a}.
\end{equation}

Our next step is to show that the group algebras $ S\Si^m_{\Lambda} $ and
$S\Si^k_{\Lambda}$ can be viewed as subalgebras of $ \E$.
\color{black}
For this purpose and inspired by \textcolor{black}{the formula ({\ref{expansiontrans})} for} $B_i$,
we define $ \BB_i \in \Ea$ as follows
\begin{equation}{\label{wedefine}}
\BB_i := \textcolor{black}{\beLambda} g_{a,1} g_{ a+1,2}
\ldots g_{2a-1,a},\, \, \, g_{i j } := g_{i+c} g_{i-1+c} \ldots  g_{j+c}
\end{equation}
\textcolor{black}{where we from now on use the notation
\begin{equation}{\label{beLambda}}
\beLambda :=\be_{\textcolor{black}{A_{\blambda}}}.
\end{equation}
}
\nomenclature[30]{$\beLambda$}{The idempotent $\be_{A_{\blambda}}$, where $\Lambda=(\blambda\mid -)$}We can now
state our next result.

\color{black}

\begin{lem}{\label{embedding}}
Suppose that $ \Lambda=(\blambda \mid \bmu) \in  {\mathcal
L}_n(\alpha) $. Then we have $S$-algebra embeddings
\begin{itemize}
\setlength\itemsep{-1.5em}
\item[(1)] $\iota:S\Si^m_{\Lambda} \hookrightarrow \Ea, \,\, \, \,\,\, \, \,  \mbox{ via }  B_i \mapsto \BB_i \mbox{ for } i \in S_{\Lambda}^m $.\newline
\item[(2)] $\iota: S\Si^k_{\Lambda} \hookrightarrow \Ea, \,\, \, \,\,\, \, \,  \mbox{ via }  B_i \mapsto \BB_i \mbox{ for } i \in S_{\Lambda}^k $.
\end{itemize}\end{lem}

\begin{demo}
  It is enough to prove part (2) of the Lemma since $ \Si^m_{\Lambda}$ is simply the parabolic subgroup of
  $ \Si^k_{\Lambda}$ corresponding to $ S_{\Lambda}^m $.
  Now $ \blambda $ is of type $ \alpha $
   and so the presence of the factor $ \beLambda $ in $ \BB_i $ gives via
   (\ref{basessubalgebras}) that $ \BB_i \in \Ea $.
  Hence, in order to show the Lemma we need to check the following three identities
  \begin{equation}{\label{threeprop}}\begin{array}{ll} a) \,\,\, \,
  \BB_i \BB_{i+1} \BB_i = \BB_{i+1} \BB_{i} \BB_{i+1} \mbox{ for } i,i+1 \in S_{\Lambda}^k
      & b) \,\, \, \,    \BB_i^2 =  \beLambda  \mbox{ for } i \in S_{\Lambda}^k
   \\
   c )  \, \, \, \,\, \,   \BB_i \BB_j  = \BB_j \BB_i \, \,   \mbox{ for }  \, \, \,  |i-j | > 1 \mbox{ and } i,j \in S_{\Lambda}^k.
&
  \end{array}
  \end{equation}

Let us first take a closer look at the expansion
$ B_i = s_{i_1} s_{i_2} \cdots s_{i_p} $
according to the definitions of $ B_i $ and $ s_{ij} $. We claim that this expansion
is a reduced expression in the $s_i$'s. Indeed, let
$ \underline{i}:= (1,2,\ldots, n) \in \seq$.  Then the right action of $ B_i $
on $ \underline{i} $, according to $ B_i = s_{i_1} s_{i_2} \cdots s_{i_p} $, changes
at the $ j $'th step $ \ldots i_p  \ldots i_p+1 \ldots $ to
$ \ldots i_p+1  \ldots i_p \ldots $, as one easily checks, and from this we
conclude, via the inversion description of the length function on $ \Si_n$,  that $ s_{i_1} s_{i_2} \cdots s_{i_p} $
indeed is a reduced expression for $ B_i$.
On the other hand, by the description of $ B_i $ in ({\ref{Bi}}) we also have that
$ B_i = s_{i_p} \cdots s_{i_2} s_{i_1} $. By length considerations this must be a reduced expression for $ B_i $ as well,
and hence via Matsumuto's Theorem we get that
\begin{equation}{\label{afterall}}
  \BB_i = \beLambda g_{i_p} \cdots g_{i_2} g_{i_1}
\end{equation}
since, after all, the $ g_i$'s verify the braid relations.

\medskip
In order to show $ a) $ and $ c) $,
we now first observe, acting once again on the sequence $ \underline{i}$ above, that the expansions of
each side of these identities in terms of $ s_i$'s
are also reduced expressions. On the other hand, by Proposition
\ref{EGpropiedades}(2) we can commute $ {\beLambda} $ to the right of $ \BB_i $ that is
$ \BB_i :=  g_{a,1} g_{ a+1,2}
\cdots g_{2a-1,a} \, {\beLambda}
$
and so we get $a) $ and $ c) $ using Matsumuto's Theorem directly on the corresponding reduced
expressions.

\medskip
In order to show $ b) $ we have to argue a bit differently. It is enough to show that
\begin{equation}{\label{bitdiffrently1}}
  \BB_i^2 = {\beLambda} g_{i_1}    \cdots g_{i_{p-1}} g_{i_p} \beLambda g_{i_p}g_{i_{p-1}} \cdots  g_{i_1} = {\beLambda}
\end{equation}
since we can use ({\ref{afterall}}) for the second expression for $ \BB_i$.
Commuting $ g_{i_p} $ past $ \beLambda $ this becomes
\begin{equation}{\label{pastthisbecomes}}
  {\beLambda} g_{i_1}    \cdots g_{i_{p-1}} g_{i_p}^2 \be_{{(A_{\blambda}})s_{i_p}} g_{i_{p-1}} \cdots  g_{i_1}
  =
  {\beLambda} g_{i_1}    \cdots g_{i_{p-1}} \left(1+ (q-q^{-1})g_{i_p}e_{i_p}\right) \be_{{(A_{\blambda}})s_{i_p}} g_{i_{p-1}} \cdots  g_{i_1}
  \end{equation}
by Proposition
\ref{EGpropiedades}(2). But $ i_i $ and $i_p+1$ are
in different blocks of $ (A_{\blambda})s_{i_p} $ and so
we have $ e_{i_p}\be_{{(A_{\blambda}})s_{i_p}}  = 0 $
by Proposition
\ref{EGpropiedades}(3). Hence (\ref{pastthisbecomes}) is equal to
\begin{equation}
  {\beLambda} g_{i_1}    \cdots g_{i_{p-1}}  \be_{{(A_{\blambda}})s_{i_p}} g_{i_{p-1}} \cdots  g_{i_1}.
\end{equation}
With the same reasoning we move $ g_{i_{p-1}} $ past $ \be_{{(A_{\blambda}})s_{i_p}} $ to arrive at
\begin{equation}
  {\beLambda} g_{i_1}    \cdots g_{i_{p-2}}  \be_{{(A_{\blambda}})s_{i_p} s_{i_{p-1}}} g_{i_{p-2}} \cdots  g_{i_1}
\end{equation}
and so on until
\begin{equation}{\label{bitdiffrently2}}
  \beLambda   \be_{{(A_{\blambda}})s_{i_p} s_{i_{p-1}} \ldots s_{i_{1} }} =
\beLambda   \be_{{A_{\blambda}}}
    = \beLambda
\end{equation}
via Proposition
\ref{EGpropiedades}(2). This proves $ b) $.

For
$ B_y $ an element of $ \Si^k_{\Lambda}$ written in reduced form as
$ B_y = B_{i_1} B_{i_2} \cdots B_{i_k} $ with $ i_j \in S_{\Lambda}^k$, we define
\begin{equation}\label{notationBy}
 \BB_y := \BB_{i_1}  \BB_{i_2} \cdots \BB_{i_k}\in    \Ea.
\end{equation}
Then by the above, $ \BB_y\in \Ea $ is
independent of the chosen reduced expression.
Since $ \beLambda$ commutes with $ \iota(S\Si^k_{\Lambda}) $ and since $ \iota(B_i) = \BB_i$
we have that
$  \iota(B_y) = \BB_y $.

It only remains to show that the induced homomorphism
$\iota: S\Si^k_{\Lambda} \rightarrow \Ea $ is an enbedding.
But this follows directly from the basis $ {\cal B}:=  \{\be_{A}g_w \} $ for $ \Ea $ given in (\ref{basessubalgebras}).
Indeed, we have that $ \iota(B_y) = \BB_y  = \beLambda g_{y} \in {\cal B} $ where
$ y \in \Si_n $ is the element obtained by expanding
$ B_y = B_{i_1} B_{i_2} \cdots B_{i_k} $ completely in terms of $ s_i$'s.
From this it also follows that $ \iota(B_y) =  \iota(B_w) $ iff $ B_y =  B_w $, proving
the injectivity of $ \iota$.
\end{demo}

\begin{obs}
The identity element of $ \iota(S\Si^k_{\Lambda}) $ is $  \beLambda $ whereas
the identity element of $ \Ea $ is $\be_\alpha $, as was seen implicitly in the proof. In particular, $ \iota$ does not preserve identity elements.
\end{obs}

\color{black}

\color{black}

\color{black}
\textcolor{black}{Recall that for any $ S$-algebra $ \cal A $, the wreath product algebra $ {\cal  A } \wr \Si_f$ is defined as
the semidirect product
$  {\cal A }^{ \otimes f} \rtimes \Si_f $ where $ \Si_f $ acts on $ {\cal A }^{ \otimes f} $ via place permutation.
If $ \cal A $ is free over $ S$ with basis $B $ then
$ {\cal  A } \wr \Si_f$ is also free over $ S $ with basis $ (b_{i_1} \otimes \cdots \otimes b_{i_f}) \otimes w $
where $ b_{i_j} \in B  $ and where $ w \in \Si_f$. 
There are canonical algebra embeddings
$ i_{{\cal A},f} : {\cal A}^{ \otimes f} \hookrightarrow {\cal  A } \wr \Si_f $ and $ j_{{\cal A}, f}: S \Si_f
\hookrightarrow {\cal  A } \wr \Si_f$ whose images generate ${\cal  A } \wr \Si_f $, subject to the following relations
\begin{equation}{\label{universal}}
j_{{\cal A}, f}(w)  i_{{\cal A}, f}(b_{i_1} \otimes \cdots \otimes b_{i_f}) =
  i_{{\cal A}, f}(b_{i_{1w^{-1}}} \otimes \cdots \otimes b_{i_{fw^{-1}}}) j_{{\cal A}, f}(w).
\end{equation}
}
Recall $ \Si_{{\Lambda}}  \leq \Si_n$, the stabilzer subgroup of the set partition $ A_{\blambda}$.
With the above notation we have an isomorphism
\begin{equation}{\label{productwreath}}
S \Si_{{\Lambda}} \cong S \Si_{ a_1 } \wr \Si_{ k_1} \otimes \cdots\otimes S \Si_{ a_r } \wr \Si_{ k_r}.
\end{equation}
\color{black}
We are interested in the following deformation of $ S \Si_{{\Lambda}}$
\begin{equation}\mathcal{H}^{wr}_{\alpha}(q):=
\HH_{a_1}(q)  \wr \Si_{k_1} \otimes \cdots\otimes \HH_{a_r}(q) \wr \Si_{ k_r}.
\end{equation}
\nomenclature[31]{$\mathcal{H}^{wr}_{\alpha}$}{The tensor product of wreath algebras associated with the partition $\alpha$}Recall that we have
$ S \Si^k_{\Lambda} = \textcolor{black}{S}\Si_{k_1} \otimes \cdots\otimes \textcolor{black}{S}\Si_{k_r }$ by ({\ref{areproductsofsymmetricgroup}}).
Let
\begin{equation} j:  S \Si^k_{\Lambda} =\textcolor{black}{S} \Si_{k_1} \otimes \cdots\otimes \textcolor{black}{S}\Si_{k_r }   \hookrightarrow \mathcal{H}^{wr}_{\alpha}(q)
\end{equation}
be the embedding
induced by the $ j_{\textcolor{black}{\mathcal{H}}_{a_i}, k_i}$'s
and let
\begin{equation}
 i: \HH_{a_1}(q)^{\otimes k_1} \otimes \cdots\otimes  \HH_{a_r}(q)^{\otimes k_r} \hookrightarrow \mathcal{H}^{wr}_{\alpha}(q)
\end{equation}
be the embedding induced by the $ i_{\Si_{a_i}, k_i}$'s.

Now
$ \HH_{a_1}(q)^{\otimes k_1} \otimes \cdots\otimes  \HH_{a_r}(q)^{\otimes k_r} $ is canonically isomorphic to the Young-Hecke algebra
$ \mathcal{H}_{\alpha^{op}}(q) $. Moreover, the multiplication map $ g_{w_1}  \cdots  g_{w_m} \mapsto \beLambda g_{w_1}  \cdots  g_{w_m}$
induces an embedding of $ \mathcal{H}_{\alpha^{op}}(q) $ in $ \Ea$.
{\color{black}Indeed, for $i $ and $i+1$ belonging to the same block of $ A_{\blambda} $ we have by
  Proposition \ref{EGpropiedades}(3) that
$$ \beLambda g_i^2 = \beLambda(1+(q-q^{-1})e_i g_i) = \beLambda(1+(q-q^{-1}) g_i)  $$
  and so the multiplication map certainly is an algebra homomorphism. The elements
  $ \beLambda g_{w_1}  \cdots  g_{w_m} $ belong to the basis
of $ \Ea$ given in  (\ref{basessubalgebras}) and so it is also injective, as claimed.}

Combining,
we get an embedding
\begin{equation}{\label{epsilonembedding}}
 \epsilon: \HH_{a_1}(q)^{\otimes k_1} \otimes \cdots\otimes  \HH_{a_r}(q)^{\otimes k_r} \hookrightarrow \Ea.
\end{equation}

With these preparations we can now extend Lemma \ref{embedding} to $\mathcal{H}^{wr}_{\alpha}(q)$ as follows.
\begin{lem}{\label{embeddingwreath}}
There is a unique embedding
\begin{equation}
\upsilon: \mathcal{H}^{wr}_{\alpha}(q) \hookrightarrow \Ea
\end{equation}
such that $ \epsilon: \HH_{a_1}(q)^{\otimes k_1} \otimes \cdots\otimes  \HH_{a_r}(q)^{\otimes k_r} \hookrightarrow \Ea$
factorizes as $ \epsilon = \upsilon \circ i $
and such that $ \iota: S \Si^k_{\Lambda} \hookrightarrow \Ea$ from the previous Lemma factorizes as $ \iota= \upsilon \circ j $.
\end{lem}

\begin{demo}
Let $ \epsilon_i: \HH_{a_i}(q)^{\otimes k_i} \rightarrow \Ea$ be the composition of
the canonical embedding $$ \HH_{a_i}(q)^{\otimes k_i} \hookrightarrow
\HH_{a_1}(q)^{\otimes k_1} \otimes \cdots\otimes \HH_{a_i}(q)^{\otimes k_i} \otimes \cdots\otimes  \HH_{a_r}(q)^{\otimes k_r}$$
with $ \epsilon$ and let $ \iota_i:  S\Si_{k_i} \rightarrow \Ea$ be the composition of
the canonical embedding
$$ S \Si_{k_i} \longrightarrow \textcolor{black}{S}\Si_{k_1} \otimes \cdots\otimes  S \Si_{k_i}  \otimes \cdots\otimes \textcolor{black}{S}\Si_{k_r } =
\textcolor{black}{S}\Si^k_{\Lambda}
\longrightarrow \mathcal{H}^{wr}_{\alpha}(q)
 $$ with $\iota$.
The existence and uniqueness of $ \upsilon$ follows from the universal property of the wreath product.
In other words, by ({\ref{universal}}) we must check that
\begin{equation}
\iota_i(w )  \epsilon_i( g_{y_1} \otimes \cdots\otimes g_{{y_{k_i}}}) =
\epsilon_i( g_{y_{1w^{-1}}} \otimes \cdots\otimes g_{y_{k_iw^{-1}}}) \iota_i(w )
\end{equation}
where $ w \in \Si_{k_i} $ and
the $ g_{y_j}$'s belong to $\HH_{a_i}(q) $.
By the definitions, this becomes the following equality in $ \Ea$
\begin{equation}{\label{universalcheck}}
\BB_w   g_{y_1} \cdots   g_{{y_{k_i}}} =
g_{y_{1w^{-1}}}  \cdots  g_{y_{k_iw^{-1}}} \BB_w
\end{equation}
where the $g_{y_j}$'s belong to the distinct Hecke algebras given
by the $ k_i $ distinct Hecke algebra factors of
$ \epsilon_i(\HH_{a_i}(q)^{\otimes k_i}) $ and similarly for the $ g_{y_{jw^{-1}}}$'s.

\medskip
To verify this we may assume $ i= 1 $ and $ r=1$. Let $ k_1:= k $ and $ a_1 := a$.
Assume
that $ g_{y_1} := g_s $ with $ s \in \{1, \ldots, ka \}$ and $ a \nmid s$ and let
$ \BB_w = \BB_j $ where $ 1 \le j < k$.
Then ({\ref{universalcheck}}) reduces to proving
\begin{equation}{\label{blackucestoproving1}}
\begin{array}{ll}
a) \, \,    \BB_j g_s = g_{s+a} \BB_j  &  \mbox{ if } s \in \{(j-1)a +1, (j-1)a +2, \ldots, ja -1 \}   \\
b) \, \,    \BB_j g_s = g_{s-a} \BB_j   &  \mbox{ if } s \in \{ ja+1, j +2, \ldots, ja +a -1 \}   \\
c) \, \,  \BB_j g_s = g_{s} \BB_j  &   \mbox{ otherwise.  }
\end{array}
\end{equation}
Let us assume that $ j= 1 $, the other cases are treated similarly.
Then in the notation of ({\ref{Bi}}) we have that $ c=0$ and by ({\ref{wedefine}})
$ g_{i j } := g_{i} g_{i-1} \cdots g_{j}$ and
$$ \BB_i := \textcolor{black}{\beLambda} g_{a,1} g_{ a+1,2}
\cdots g_{2a-1,a}.
$$
We then have to prove
  \begin{equation}{\label{wefirstprove}}
  \begin{array}{ll}
   a1)  \, \,  g_{a,1} g_{ a+1,2} \cdots g_{2a-1,a} g_s= g_{s + a  } g_{a,1} g_{ a+1,2} \cdots g_{2a-1,a} &
    \mbox{ for } s \in \{1, 2,  \ldots, a-1 \} \\
   b1)  \, \,  g_{a,1} g_{ a+1,2} \cdots g_{2a-1,a} g_s= g_{s-a  } g_{a,1} g_{ a+1,2} \cdots g_{2a-1,a} &
    \mbox{ for } s \in \{a+ 1, a+ 2,  \ldots, 2a-1 \} \\
   c1)  \, \, g_{a,1} g_{ a+1,2} \cdots g_{2a-1,a} g_s= g_{s } g_{a,1} g_{ a+1,2} \cdots g_{2a-1,a} &
    \mbox{ otherwise. }
  \end{array}
     \end{equation}
But using only braid relations one checks that
$ g_{a,b} g_s = g_{s-1} g_{a,b} $ if $  s \in \{ b+1, \ldots, a \} $, which gives $ b1)$.
On the other hand, as mentioned above we have that
$$  (g_{a,1} g_{ a+1,2} \cdots g_{2a-1,a}  )^{\ast} =  g_{a,1} g_{ a+1,2} \cdots g_{2a-1,a} $$
(actually this can also be shown directly using only the commuting braid relations)
and hence the $ a1) $ case follows by applying $ \ast$ to the $ b1) $ case. The
remaining case $ c1) $ is
easy.

\medskip
Now the general $  g_{y}$-case of ({\ref{universalcheck}}) follows from $ a), b) $ and $ c) $ by expanding
$ g_y = g_{s_1} \cdots g_{s_l} $ in terms of simple $ g_s$'s and pulling $ \BB_i $ through all factors.
Finally the general $ \BB_w $-case is obtained the same way by expanding $ \BB_w =
\BB_{i_1} \cdots \BB_{i_l} $ and pulling all factors through.
\medskip

To show that $ \upsilon$ is an embedding we argue as in the previous Lemma. Indeed,
by construction, the images under $ \upsilon$ of of the canonical basis vectors of
$ \mathcal{H}^{wr}_{\alpha}(q) $ belong to the basis $  \cal B$ for $ \Ea$ and are pairwise distinct, proving
that $ \upsilon$ is an embedding.
\end{demo}
\color{black}

\medskip
We are now finally ready
to give the construction of the cellular basis for $ \Ea$.
As in the Yokonuma-Hecke algebra case, we first construct, for each $ \Lambda \in {\cal L}_n(\alpha) $, an
element $ m_{\Lambda} $ that acts as the starting point of the basis. Suppose that $ \Lambda = (\blambda \mid \bmu ) $ is as above
with $ \blambda = (\lambda^{(1)}, \ldots, \lambda^{(m)}) $ and $ \bmu = (\mu^{(1)}, \ldots,\mu^{(q)})$. We then define $ m_{\Lambda} $ as follows
\begin{equation}{\label{startingbasiselement}}
m_{\Lambda}:= \textcolor{black}{\beLambda}   x_{\blambda}     b_{\bmu}.
\end{equation}
Let us explain the factors of the product.
Firstly, \textcolor{black}{$\beLambda$} is the idempotent defined in
(\ref{beLambda}).
Secondly, $ x_{\blambda} \in \Ea $ \textcolor{black}{is an analogue
for $ \Ea $ of the element $ x_{\blambda}$ for the Hecke
algebra}, or the element $ m_{\blambda} $ in  the Yokonuma-Hecke algebra case. It is given as
\begin{equation} x_{\blambda} :=\textcolor{black}{\beLambda  \sum_{w\in\Si_{\blambda}}q^{\ell(w)}g_{w}.}
\end{equation}
Mimicking the argument in (6) of Lemma \ref{l5} we get that
\begin{equation}{\label{justasinYH}}
x_{\blambda} g_w = g_w  x_{\blambda} = q^{l(w)} x_{\blambda} \, \, \, \, \mbox{ for } w \in \Si_{\blambda}.
\end{equation}
Finally, in order to explain the factor $ b_{\bmu} $ we \textcolor{black}{recall from
  ({\ref{areproductsofsymmetricgroup}}) the decomposition
  \begin{equation}{\label{recalldecom}} \Si^m_{\Lambda}  \cong \Si_{m_1} \times \cdots\times \Si_{m_q}
  \end{equation}
 where $ m_i = |\mu^{(i)} |  $.
Let $ x_{\bmu}(1)  $
be the $ q= 1 $ specialization of the Murphy element corresponding to the multipartition $\bmu$,
it may be viewed as an element of $S \Si^m_{\Lambda}$.
Then $   b_{\bmu} $ is defined as
\begin{equation}
  b_{\bmu}:=\iota(x_{\bmu}(1)) \in \Ea
\end{equation}
where $ \iota:S\Si^m_{\Lambda} \hookrightarrow \Ea $ is the embedding from Lemma \ref{embedding}.
}
Let $ \et^{\Lambda} $ be the $ \Lambda$-tableau given in the obvious
way as $  \et^{\Lambda} :=(\bT^{\blambda} \mid \bT^{\bmu}) $. Then $ \et^{\Lambda}$ is a maximal $
\Lambda$-tableau, that is the only standard $ \Lambda$-tableau $ \et
$ satisfying $ \et \unrhd \et^{\Lambda} $ is $ \et^{\Lambda} $
itself. For $ \es =(\Bs \mid  \bu) $ a $
\Lambda$-tableau we define $ d(\es) := (d(\textcolor{black}{\Bs}) \mid \iota(d(\bu))) $ where $ {d(\Bs)}  \in \Si_n $
\textcolor{black}{as usual is given by
  $ \bT^{ \blambda } d(\Bs) =  \Bs $ and $ d(\bu) \in \Si^m_{\Lambda}$
  by $ \bT^{ \bmu} d(\bu) = \bu $. For simplicity, we often write $(d(\textcolor{black}{\Bs}) \! \!\mid d(\bu)) $
for $ (d(\textcolor{black}{\Bs}) \mid \iota(d(\bu))) $.}
Note that since
$\bu=(\U_1, \ldots,\U_q)$ is always of the initial kind, we have a decomposition
$ d(\bu) = (d(\U_1),  \ldots, d(\U_q)) $, according to ({\ref{recalldecom}}),
and also
$$ \BB_{d(\bu)}=  \BB_{d(\U_1)}\cdots \BB_{d(\U_q)}.$$
Finally, we define the main object of this section. For $\es=(\Bs \mid \bu),\, \et=(\bT \mid \bv) $ row standard $\Lambda$-tableaux we define
\begin{equation}{\label{maindef}}
m_{\es\et}:=  g_{d(\Bs)}^{\ast}  \textcolor{black}{ \beLambda}
\BB_{d(\bu)}^{\ast}x_{\blambda} b_{\bmu}\BB_{d(\bv)}  \,
g_{d(\bT)}.
\end{equation}

Our aim is to prove that the $  m_{\es\et}  $'s, with $\es$ and $\et$ running over standard $\Lambda$-tableaux,
form a cellular basis for $ \Ea$.
To achieve this goal we first need to work out commutation rules between the various ingredients of $ m_{\es \et} $.
The rules shall be formulated in terms of a certain $ \circ $-action on tableaux that we explain now.

\medskip
    {\it \textcolor{black}{Let $ B_y \in \Si^k_{\Lambda} $.
        From now on, when confusion should not be possible,
        we shall write $ \Bs y $ for $ \Bs B_y $ where $ \Bs  $ is the first part of a $ \Lambda$-tableau
        and where the action of $ B_y$ is given by the complete expansion
        of $ B_y $ in terms of $ s_i$'s. }}
\medskip

Let $ \es = (\Bs \mid \bu) $ be a $ \Lambda$-tableau.
We then define a new multitableau $ y \circ \Bs$ as follows.
Set first $ \Bs_1:= \Bs {y^{-1}} = (\s_1^{(1)} , \ldots, \s_1^{(m)}) $. Then $ y \circ \Bs$ is given by the formula
\begin{equation}{\label{circleaction}} y \circ \Bs := (\s_1^{ (1) y}, \ldots, \,  \s_1^{ (m) y}).
\end{equation}
With this notation we have the following Lemma which is easy to verify.
\begin{lem}{\label{easytoverify}}
The map $ (y, \Bs ) \mapsto y \circ \Bs $ defines a left action of $ \Si_{\Lambda}^k $ on the set of multitableaux $ \Bs$
such that $  Shape(\Bs)  =_1  \blambda  $ where $ \blambda $ is the first part \textcolor{black}{of a $ \Lambda$-tableau};
that is
$ Shape(\Bs)  $ and $ \blambda $ are equal multipartitions up to a permutation.
Moreover, if $ \Bs $ is of the initial kind then also $y \circ \Bs$ is of the initial kind, and if
$ y \in \Si^m_{\Lambda} $ then $y \circ \Bs =   \Bs $.
\end{lem}

\begin{exa}
We give an example to illustrate the action. As can be seen, it permutes the partitions of the multitableau, but keeps the
numbers.
Consider
$$\Bs:={\footnotesize\left(\;\young(12)\,,\,\young(3,4)\,,\,\young(5,6)\,,\,\young(79,8)
\,,\,\ytableausetup{centertableaux,boxsize=1.3em}\begin{ytableau}10&11&12\end{ytableau}\,,\,\begin{ytableau}13\\14\\15\end{ytableau}\;\right)}\quad \mbox{ and }\quad \textcolor{black}{B_y:=B_1B_2B_1B_4B_5}.$$
We first note that $\Bs y^{-1}={\footnotesize\left(\;\young(56)\,,\,\young(3,4)\,,\,\young(1,2)\,,\,\ytableausetup{centertableaux,boxsize=1.3em}\begin{ytableau}10&12\\11\end{ytableau}
\,,\,\ytableausetup{centertableaux,boxsize=1.3em}\begin{ytableau}13&14&15\end{ytableau}\,,\,\begin{ytableau}7\\8\\9\end{ytableau}\;\right)}$. Then we have
$$y\circ \Bs={\footnotesize\left(\;\young(1,2)\,,\,\young(3,4)\,,\,\young(56)\,,\,\young(7,8,9)
\,,\,\ytableausetup{centertableaux,boxsize=1.3em}\begin{ytableau}10&12\\11\end{ytableau}\,,\,\begin{ytableau}13&14&15\end{ytableau}\;\right)}.$$
\end{exa}
Let $\Bs $ and $ \bT $ be $ \blambda$-multitableaux. Then we define $  x_{\Bs\bT} \in \Ea$, just as for the Yokonuma-Hecke algebra, that is
\begin{equation}
 x_{\Bs\bT} := g_{d(\Bs)}^{\ast} x_{\blambda} g_{d(\bT)} \in \Ea.
\end{equation}

\textcolor{black}{
The following remark
is an analogue of Remark
{\ref{importantremark}} for the Yokonuma-Hecke algebra.
\begin{obs}{\label{importantremark2}}
Let $ \Bs $ and $ \bT$ be multitableaux of the initial kind and let
$$
d(\s)=\left( d(\s^{(1)}), d(\s^{(2)}), \cdots, d(\s^{(m)}) \right) \quad \mbox{ and }\quad
 d(\T)=\left(d(\T^{(1)}), d(\T^{(2)}), \cdots,  d(\T^{(m)})\right)
$$
be the decompositions given in ({\ref{decompositionref}}).
Then, under the embedding from ({\ref{epsilonembedding}}) of the Young-Hecke algebra
$$  \epsilon:  \mathcal{H}_{\alpha}(q)
 \hookrightarrow \Ea $$
 we have that
$ \epsilon(x_{\s^{(1)} \T^{(1)}} \otimes
x_{\s^{(2)} \T^{(2)}} \otimes  \cdots \otimes  x_{{\s^{(m)} \T^{(m)}}})  = x_{\Bs\bT}.$
\end{obs}
}

The next Lemma gives the promised commutation formulas.
\begin{lem}{\label{importantcommutation}}
Suppose $ \es=(\Bs \mid \bu) $ and $ \et=(\bT \mid \bv) $ are $\Lambda$-tableaux such that $ \Bs $ and $ \bT $ are of the initial kind
and suppose that $ B_y \in \Si^k_{\Lambda} $.
Then we have the following formulas in $ \Ea$.
\begin{itemize}
\setlength\itemsep{.1em}
\item[(1)]  $\textcolor{black}{\beLambda}\BB_{y}  g_{d(\Bs)} = \textcolor{black}{\beLambda} g_{d(y \circ \Bs)} \BB_{y}.$
\item[(2)]  $\textcolor{black}{\beLambda} \BB_y x_{\Bs\bT}    = \textcolor{black}{\beLambda}  x_{y \circ \Bs, y \circ \bT }  \BB_y$.
\end{itemize}
\end{lem}

\begin{demo}
\textcolor{black}{In order to prove (1) we may assume that $ \BB_y = \BB_i $,
since
$ y \mapsto y \circ \Bs $ is a left action.
Now $ \Bs= (\s^{(1)}, \ldots, \s^{(m)})$ is of the initial kind and so we have a decomposition
$ g_{ d(\Bs) } = g_{ d(\s^{(1)})}   \cdots g_{ d(\s^{(m)})}$ with the $ g_{ d(\s^{(i)})}$'s belonging to
Hecke algebras running over the distinct indices given by the symmetric group factors $ \Si_{k_i}$ of $\Si^k_{\Lambda}$.
Then by ({\ref{blackucestoproving1}}) we have that
$$\textcolor{black}{\beLambda}\mathbb{B}_i g_{d(\Bs)}=\textcolor{black}{\beLambda}\mathbb{B}_y g_{d(\s^{(1)})}    \cdots g_{d(\s^{(m)})} =
\textcolor{black}{\beLambda}
g_{d(\textcolor{black}{\R}^{(1)})   }    \cdots g_{d(\textcolor{black}{\R}^{(m)})   }   \mathbb{B}_i $$
where $ d(\s^{(k)}) = d(\textcolor{black}{\R}^{(k)}) $ for $ k \neq i, i+1$
and where $ d(\s^{(i)}), d(\textcolor{black}{\R}^{(i+1)}), d(\s^{(i+1)}), d(\textcolor{black}{\R}^{(i)}) $ are related as in
({\ref{blackucestoproving1}}): each factor $ g_s $ of $ d(\s^{(i)}) $ is replaced by $ g_{s+k_i}$ to arrive at $ d(\textcolor{black}{\R}^{(i+1)})$ and
similarly for $ d(\s^{(i+1)})$ and $ d(\textcolor{black}{\R}^{(i)}) $.
But this means exactly that
$$g_{ d( B_i \circ \Bs)} = g_{d(\textcolor{black}{\R}^{(1)}) } \cdots g_{d(\textcolor{black}{\R}^{(m)})}$$
and so (1) follows.}

\medskip
On the other hand, applying $ \ast $ to (1) and using that $ \BB_{y}^{\ast} = \BB_{y^{-1}} $,  we find
 $ \textcolor{black}{\beLambda} \BB_{y^{-1}} g_{d(y \circ \Bs)}^{\ast} =  \textcolor{black}{\beLambda} g_{d(\Bs)}^{\ast} \BB_{y^{-1}}    $,
that is
$$ \textcolor{black}{\beLambda} \BB_{y} g_{d(\Bs)}^{\ast} =  \textcolor{black}{\beLambda} g_{d( y \circ \Bs)}^{\ast} \BB_{y} .   $$
(Alternatively, one can also repeat the argument for (1)).
Now (\ref{universalcheck}) can be formulated as follows
 \begin{equation}
{ \beLambda}\mathbb{B}_yg_k={\beLambda}g_{kB_{y}^{-1}}\mathbb{B}_y
 \end{equation}
and hence we get  $ \textcolor{black}{\beLambda} \BB_{y} x_{\blambda} = \textcolor{black}{\beLambda}
x_{\bmu} \BB_y $, where $\bmu=Shape(y\circ \bT^{\blambda})$.
In view of the definitions
this shows (2).
\end{demo}

\begin{coro}\label{commute}
The factor $x_{\blambda}$ of $m_{\es\et}$ commutes with each of
the factors $ \BB_{d(\bu)}^{\ast} $, $ b_{\bmu}$ and $\BB_{d(\bv)}$ of $ m_{\es\et} $.
Furthermore,
\begin{equation}\label{invarianteE}
m_{\es\et}^{*}=m_{\et\es}.
\end{equation}
\end{coro}

\begin{demo}
  Setting $ \Bs = \bT = \bT^{\blambda}$ in part (2) of the Lemma we get for $ \BB_y \in \Si^m_{\Lambda} $ that
  \begin{equation}
    \textcolor{black}{\beLambda} \BB_y x_{\blambda}    = \textcolor{black}{\beLambda}  x_{y \circ \bT^{\blambda}, y \circ \bT^{\blambda} }  \BB_y = \textcolor{black}{\beLambda}  x_{\blambda}  \BB_y
\end{equation}
  since $  \bT^{\blambda} $ is of the initial kind and therefore $ y   \circ \bT^{\blambda} =  \bT^{\blambda}$
  by Lemma \ref{easytoverify}. This shows the first claim. To show the second claim, we use the first claim
  together with $\textcolor{black}{\beLambda}\BB_y=\BB_y\textcolor{black}{\beLambda}$ for all $\textcolor{black}{y}\in\Si_{\Lambda}^{k}$,
  as follows from Proposition \ref{EGpropiedades}\textcolor{black}{(2)} and the definition of  $\BB_i$, to get
$$ m_{\es\et}^{*}
=g_{d(\bT)}^{\ast}\BB_{d(\bv)}^{\ast}b_{\bmu}^{\ast}x_{\blambda}^{\ast}\BB_{d(\bu)}\textcolor{black}{\beLambda^{\ast}}g_{d(\Bs)}=
g_{d(\bT)}^{\ast}   \textcolor{black}{\beLambda}
\BB_{d(\bv)}^{\ast}x_{\blambda} b_{\bmu}\BB_{d(\bu)}g_{d(\Bs)}
= m_{\et\es}
$$
as claimed.
  \end{demo}

\medskip
\color{black}
We need the following technical Lemma.
\begin{lem}{\label{technical}}
  Suppose that $ \Lambda = (\blambda \mid \bmu)  $ such that $ \Bs $ is a $\blambda $-multitableau.
  Let $ w_{{\Bs}}$ be
the distinguished representative for $ d(\Bs) $ with respect to $ \Si_{\norm{ \blambda } }$, that is we have the decomposition $ d(\Bs) = d(\Bs_0) w_{{\Bs}}$, as in \textcolor{black}{(\ref{decompositionintialkind})}.
 Let $ B_y \in \Si_{\Lambda}^k$. Then,
  in $ \Ea$
  we have the identity
  $ \BB_y g_{w_{{\Bs}}} = \textcolor{black}{\beLambda} g_{B_y w_{{\Bs}}} $ (even though in general
  $ l( B_y w_{{\Bs}} ) \neq l( B_y) + l( w_{{\Bs}} ) $). Moreover, for any multitableau $ \bT_0 $ of the initial kind with respect to $ \blambda$, we have
that $ \textcolor{black}{\beLambda}g_{ d( \bT_0) B_y    w_{{\Bs}} } = \textcolor{black}{\beLambda} g_{ d( \bT_0) } \BB_y   g_{w_{{\Bs}}}$.
  \end{lem}

\begin{demo}
The ingredients of the proof are already present in the proof of part b) of Lemma {\ref{embedding}}.
As before we set $ A_{\blambda} = \{ I_1, I_2, \ldots, I_q \} $, with blocks $ I_i$. Let $ B_y = s_{i_1 } \ldots s_{i_r } $
be the expansion of $ B_y $ according to the definitions and let  $ w_{{\Bs}} = s_{j_1 } \ldots s_{j_s } $ be a reduced
expression. The action of $ B_y $ involves at each step distinct blocks, that
is $ i_k  $ and $ i_k +1 $ occur in distinct blocks of $ (A_{\blambda}) s_{i_1 } \ldots s_{i_{k-1 }}$ for all $ k $.
A similar property holds for $ w_{{\Bs}} $ since it is the distinguished coset representative for $ d(\Bs) $
with respect to $ \Si_{\blambda} $.
  But the blocks of $ (A_{\blambda}) B_y$ are a permutation of the blocks of $ A_{\blambda}$, that is
  $ (A_{\blambda}) B_y = A_{\blambda} $ as set partitions, and so also the action of the concatenation
  $ s_{i_1 } \ldots s_{i_r } s_{j_1 } \ldots s_{j_s } $ on $ A_{\blambda} $ involves at each step distinct blocks.

  We now transform $ s_{i_1 } \ldots s_{i_r } s_{j_1 } \ldots s_{j_s } $ into a reduced expression for $ B_y  w_{{\Bs}} $
  using the Coxeter relations of type $A$. We claim that these Coxeter relations map a sequence
  $ s_{\iota_1 } \ldots s_{\iota_t }  $ having the property of acting at each step in distinct blocks
  to another sequence having the same
  property. This is clear for the commuting Coxeter relations $ s_i s_j = s_j s_i$ and also for the quadratic
  relations $ s_i^2 =1 $. In the case of the braid relations $ s_i s_{i+1}s_i = s_{i+1} s_{i}s_{i+1} $
  we observe that both $ s_i s_{i+1}s_i$ and $ s_{i+1} s_{i}s_{i+1} $ have the above property with respect to
  $ A=\{ J_1, \ldots, J_u \} $ exactly when all three numbers $ i, i+1 $ and $ i+2$ occur in three
  distinct blocks $ J_i $ of $A$
  and so the claim follows also in that case.

  Now by definition $ \BB_y g_{w_{{\Bs}}} =  \textcolor{black}{\beLambda}g_{i_1 } \ldots g_{i_r } g_{j_1 } \ldots g_{j_s } $
  and so the above sequence of Coxteter relations will transform $ \BB_y g_{w_{{\Bs}}} $ to
  $  \textcolor{black}{\beLambda} g_{B_y w_{{\Bs}}}$. Indeed, for each occurence of the relation $ s_i^2 =1 $
  we have by part (3) of Proposition \ref{EGpropiedades} a corresponding relation
  \begin{equation}
   \mathbb{E}_{A}  g_i^2 = \mathbb{E}_{A} (1+ (q-q^{-1})e_i  g_i ) =  \mathbb{E}_{A}
\end{equation}
whenever $ i $ and $ i+1 $ are in distinct blocks of $ A$. This proves the first statement of the Lemma.
The second statement follows from the first since $ B_y w_{\Bs} $ is the distinguished
representative for its class with respect to $ \Si_{\norm{ \blambda }}$ as follows from the characterization of
distinguished representatives as row standard tableaux, of shape $ \beta= \norm{ \blambda}^{op}$ in this case. Indeed,
$ \T^{\beta}B_y w_{\Bs} $ is obtained from $ \T^{\beta} w_{\Bs} $ by permuting some rows
   and so one is row standard iff the other is row standard.
\end{demo}

\color{black}

\medskip
The following Lemma is the $\Ea$-version of Lemma \ref{l1} in the Yokonuma-Hecke
algebra case.
\begin{lem}{\label{secondlast}}
Suppose that $\Lambda \in  {\mathcal  L}_n(\alpha)$ and that
$\es = \textcolor{black}{(\Bs \mid \bu)}$ and $\et=\textcolor{black}{(\bT \mid \bv)} $
are row standard $\Lambda$-tableaux.
Then for every $h\in \Ea$  we have that $ m_{\es \et}h $
is a linear combination of terms of the form
$m_{\es\ev}$ where $\ev$ is a row
standard $\Lambda$-tableau. A similar statement holds for $h m_{\es\et}$.
\end{lem}

\begin{demo}
\color{black}
The idea is to repeat the arguments of Lemma \ref{l1}. It is enough to consider the $ m_{\es \et}h $ case.
Using \textcolor{red}{Corollary \ref{commute}} we have that
\begin{equation}{\label{mtimesh}}m_{\es\et}=  g_{d(\Bs)}^{\ast}   \textcolor{black}{\beLambda}
\BB_{d(\bu)}^{\ast}x_{\blambda} b_{\bmu}\BB_{d(\bv)}  \,
g_{d(\bT)} =
g_{d(\Bs)}^{\ast}
\BB_{d(\bu)}^{\ast} b_{\bmu}\BB_{d(\bv)} x_{\blambda} \textcolor{black}{\beLambda}  \,
g_{d(\bT)}.
\end{equation}
Since $ h $ is general we reduce to the case $ g_{d(\bT)} = 1 $, that is $ \bT = \bT^{\blambda}$.
We may assume that $ h = \be_A g_v $ since such elements form a basis for $ \Ea$.
But
the $ \be_A$'s are orthogonal idempotents, as was shown in Proposition
\ref{EGpropiedades}\textcolor{black}{(1)},
and so we may further reduce to the case $ h = g_v$.
We have a decomposition $ v = v_0 d(\Bv) $ with $ v_0 \in \Si_{\blambda} $ and $ \Bv $
a row standard $ \blambda $-multitableau such that $l(w) =  l(v_0) + l({d(\Bv)})  $. Hence
via ({\ref{justasinYH}}) we get
that $ m_{\es\et}g_v $ is a multiple of
\begin{equation}{\label{mtimesh}}
g_{d(\Bs)}^{\ast}
\BB_{d(\bu)}^{\ast} b_{\bmu}\BB_{d(\bv)}  \textcolor{black}{\beLambda}  x_{\blambda}g_{ d(\Bv)} = g_{d(\Bs)}^{\ast}
\BB_{d(\bu)}^{\ast} \textcolor{black}{\beLambda} b_{\bmu} x_{\blambda} \BB_{d(\bv)}    g_{ d(\Bv)}
= m_{\es\ev}
\end{equation}
where $ \ev = (\Bv \mid \bv) $.
\end{demo}

\medskip
\color{black}
Our next Lemma is the analogue for $ \Ea $ of Lemma \ref{l3}. It is the key Lemma for our results on $ \Ea$.
\begin{lem}\label{Lemmatableaustandar}
Suppose that $\Lambda \in  {\mathcal L}_n(\alpha)$ and that
$\es$ and $\et $ are row standard $\Lambda$-tableaux. Then there are
standard tableaux $\eu$ and $\ev$ such that $\eu \unrhd \es, \ev \unrhd \et$  and  such that $m_{\es\et}$ is a linear combination of the
elements $m_{\eu\ev}$.
\end{lem}

\begin{demo}
Let $ \Lambda= (\blambda \mid \bmu) $,
$ \es= (\Bs \mid \bu) $ and $ \et = (\bT \mid \bv) $. Then we have
\begin{equation}{\label{takeintoaccount}}
m_{\es\et}=  g_{d(\Bs)}^\ast    \textcolor{black}{\beLambda} \BB_{ d(\bu)}^{\ast} b_{ \bmu}
x_{\blambda}
\BB_{ d(\bv)} g_{d(\bT)}.
\end{equation}
Suppose first that standardness fails for \textcolor{black}{$ \Bs $ or $\bT.$}
The basic idea is then to proceed as in the proof of Lemma \ref{l3}. There
exist multitableaux $ \Bs_0 $ and $ \bT_0 $ of the
initial kind together with $ w_{\Bs},w_{\bT} \in \Si_n $ such that $ d(\Bs) =
d(\Bs_0) w_{\Bs} $, $ d(\bT) = d(\bT_0) w_{\bT} $ and $ \ell(d(\Bs)) = \ell(d(\Bs_0))+\ell(
w_{\Bs}) $ and $ \ell(\bT)) = \ell(d(\bT_0))+\ell( w_{\bT}) $. That is, $ w_{\Bs} $ and $ w_{\bT} $
are distinguished right coset representatives for $ d(\Bs) $ and $ d(\bT) $ with respect to
$ \Si_{\norm{\blambda}} $
and ({\ref{takeintoaccount}}) becomes
\begin{equation}{\label{takeintoaccountbecomes}}
m_{\es\et}=  g_{w_{\Bs}}^{\ast} g_{d(\Bs_0)}^\ast
 \textcolor{black}{\beLambda} \BB_{ d(\bu)}^{\ast}
x_{\blambda} b_{ \bmu}
\BB_{ d(\bv)} g_{d(\bT_0)}g_{w_{\bT}}
\end{equation}
since the two middle terms commute. Note that the factor $ \textcolor{black}{\beLambda} $ commutes with all other except the two extremal
factors of ({\ref{takeintoaccountbecomes}}). Expanding $ b_{ \bmu}
\BB_{ d (\bv)}  $ completely as a linear combination of
\textcolor{black}{$\BB_{y_{\bT}} \!$'s with $ B_{y_{\bT}} \in  \Si_{\Lambda}^{m} $} and
\textcolor{black}{setting $ \BB_{ y_{\Bs}} := \BB_{ d (\bu)} $} we get via \textcolor{black}{Lemma \ref{easytoverify} and Lemma
{\ref{importantcommutation}}} that ({\ref{takeintoaccountbecomes}})
is a linear combination of terms
\begin{equation}{\label{becomesalinear}}
g_{w_{\Bs}}^{\ast}   \BB_{y_{\Bs}}^{\ast}
x_{  \Bs_0  \bT_0  }
\BB_{y_{\bT}}
g_{w_{\bT}}  =
g_{w_{\Bs}}^{\ast}
x_{  \Bs_0  \bT_0  }  \BB_{y}
g_{w_{\bT}}
\end{equation}
where $ \BB_{y} =  \BB_{y_{\Bs}}^{\ast} \BB_{y_{\bT}}    $.
For each appearing \textcolor{black}{$B_y \in \Si_{\Lambda}^{m}  $ we have
by Lemma {\ref{technical}} that}
$$ \textcolor{black}{\beLambda}  g_{B_{y}} g_{w_{\bT}} = \textcolor{black}{\beLambda}  g_{B_{y} w_{\bT}}.
$$
Thus (\ref{becomesalinear}) becomes a linear combination of terms
\begin{equation}
g_{w_{\Bs}}^{\ast}   \textcolor{black}{\beLambda}
x_{    \Bs_0     \bT_0  }
g_{y_{\bT,1}}
\end{equation}
where $y_{\bT,
1} := \textcolor{black}{B_{y}}w_{\bT} $. {We now proceed as in the Yokonuma-Hecke algebra case. Via Remark \ref{importantremark2}}
we apply Murphy's result
\cite[Theorem 4.18]{Mur} on $ x_{\Bs_0    \bT_0  }  $,
thus rewriting it as a linear combination of $ x_{ \Bs_1
\bT_1 } $ where $ \Bs_1  $ and $\bT_1 $ are standard
$\bnu$-multitableaux of the initial kind, \textcolor{black}{satisfying} $ \Bs_1 \unrhd   \Bs_0
$ and $\bT_1 \unrhd  \bT_0 $. We then get that
(\ref{becomesalinear}) is a linear combination of such terms
\begin{equation}{\label{HenceBecomesnow}}
g_{w_{\Bs}}^{\ast}  \textcolor{black}{\beLambda}   x_{  \Bs_1  \bT_1   }  g_{y_{\bT,1}}.
\end{equation}
Let $ \bnu = (\nu^{(1)},
\ldots, \nu^{(m)})$. It need not be an increasing multipartition and our task is to fix this problem.

We determine a $ B_{\sigma} \in \Si^{k}_{\Lambda} $ such that the multipartition $ \bnu^{ord} := (\nu^{ (1)\sigma }, \ldots, \nu^{(m)\sigma}) $
is increasing.
Then, using (2) of Lemma {\ref{importantcommutation}} we get
that
({\ref{HenceBecomesnow}}) is equal to
\begin{equation}{\label{104}}
g_{  w_{\Bs}  }^{\ast}  \BB_{\sigma}^{\ast} \BB_{\sigma} x_{  \Bs_1 \bT_1   }   g_{   y_{\bT,1} } =
g_{  w_{\Bs}  }^{\ast}  \BB_{\sigma}^{\ast}  x_{\sigma \circ   \Bs_1, \sigma \circ  \bT_1   } \BB_{\sigma}  g_{   y_{\bT,1} } =
g_{  y_{\Bs,2}  }^{\ast}    x_{\sigma \circ   \Bs_1, \sigma \circ  \bT_1   }   g_{   y_{\bT,2} }
\end{equation}
where $ y_{\Bs,2} := B_{\sigma} w_{\Bs} $ and $ y_{\bT,2} :=
B_{\sigma} y_{\bT,1} $, and where we used Lemma {\ref{technical}} once again. Here $ \bT^{\bnu^{ord}} y_{\Bs,2} $ and $
\bT^{\bnu^{ord}} y_{\bT,2} $ are standard $ \bnu^{ord}
$-multitableaux but not necessarily increasing, and so we must now fix this problem.
Let therefore  $
\Si^{m^{\prime}} $ be the subgroup of $ \Si_{\Lambda}^{k} $ that
permutes equal $ \nu^{(i)}$'s.  We can then find $ \sigma_1, \sigma_2 \in
\Si_{\Lambda}^{m^{\prime}} $ such that $ \bT^{\bnu^{ord}}B_{\sigma_1}
y_{\Bs,2}  $ and $ \bT^{\bnu^{ord}}B_{\sigma_2} y_{\bT,2} $ are
increasing $ \bnu^{ord}$-tableaux. With these choices, ({\ref{104}})
becomes via Lemma {\ref{technical}}
\begin{equation}{\label{becomesequal}}
g_{  \textcolor{black}{y_{\Bs,3}}  }^{\ast}    x_{\sigma \circ   \Bs_1,
\sigma \circ \bT_1 } \BB_{\sigma_1}  \BB_{\sigma_2}^{\ast}  g_{
\textcolor{black}{y_{\bT,3}} }
\end{equation}
where $ y_{\Bs,3} := B_{\sigma_1} y_{\Bs,2} $ and $
y_{\textcolor{black}{\bT},3} := B_{\sigma_2} y_{\bT,2} $, and where we
used (2) of Lemma {\ref{importantcommutation}} to show that $ \sigma
\circ   \Bs_1 $ and $ \sigma \circ  \bT_1  $ are unchanged by the
commutation with $ \BB_{\sigma_1}  $. We now
set $ \Bs_3:= \bT^{\bnu^{ord}} d( \sigma \circ   \Bs_1 ) y_{\Bs,3}$
where $d( \sigma \circ   \Bs_1 ) $ is calculated with
respect to $ Shape( \sigma \circ   \Bs_1) = \bT^{\bnu^{ord}} \! \! \!\! \! \! $ of course, and similarly $
\bT_3:= \bT^{\bnu^{ord}} d( \sigma \circ   \bT_1 ) y_{\bT,3}$.
Then $ \Bs_3  $ and $ \bT_3 $ are increasing standard
multitableaux of shape $ \bnu^{ord} $ and we get via Lemma {\ref{technical}} that
({\ref{becomesequal}}) is equal to
\begin{equation}{\label{becomesequal2}}
g_{ d( \Bs_3 ) }^{\ast}    x_{\bnu^{ord}   } \BB_{\sigma_1}  \BB_{\sigma_2}^{\ast}  g_{   d(\bT_3)}
\end{equation}
since
$g_{d(\sigma \circ   \bT_1 )} $ and $ \BB_{\sigma_1}^{\ast}  \BB_{\sigma_2}$ commute by Lemma {\ref{importantcommutation}}.

\color{black}

\medskip
In order to show that ({\ref{takeintoaccount}}) has the form $ m_{\eu \ev} $ stipulated by the Lemma, we must now treat the
factor $ \BB_{\sigma_1} \BB_{\sigma_2}^{\ast}  $.
But since $\Si^{m^{\prime}}$ is a product of symmetric groups, $ B_{\sigma_1}  B_{\sigma_2}^{\ast} $
can be written as a linear combination of $ x_{ \bu^{\prime}  \bv^{\prime} }(1)$, with $  \bu^{\prime} $ and $  \bv^{\prime} $ running
over multitableaux of the initial kind according to the factors of $\Si^{m^{\prime}}$
and where once again $ x_{ \bu^{\prime}  \bv^{\prime} }(1)$ is the usual
Murphy standard basis element, evaluated at $1$, for that product. Thus
\textcolor{black}{(\ref{becomesequal2}) becomes
\begin{equation}
g_{ d( \Bs_3 ) }^{\ast}   \BB_{\bu^{\prime}}^{\ast} \textcolor{black}{b_{\bmu}}x_{\bnu^{ord}   }   \BB_{\bv^{\prime}}g_{   d(\bT_3)} = m_{\eu\ev}
\end{equation}
where $ \eu = (\Bs_3 \mid \bu^{\prime}) $ and $ \ev = (\bT_3 \mid \bv^{\prime}) $.
Note that by the constructions we have that the shape of $ \Bs_3 $ and $ \bT_3$ is of type $ \alpha $ and that
$\eu \unrhd \es, \ev \unrhd \et $ and so the Lemma is proved in this case.}

Finally,  we must now treat the case where standardness holds for $ \Bs $ and $\bT$, but fails for $\bu $ or $\bv$. But this case is much easier, since we can here
apply Murphy's theory directly, thus expanding the nonstandard terms in terms of standard terms.

\end{demo}

We are now ready to state and prove the main Theorem of this section.
\begin{teo}{\label{mainTheorem}}
Let $ { \cal BT}_{ \! \! \! n} := \{  m_{ \es \et} \,  | \, \es, \et \in \std(\Lambda), \Lambda \in {\cal L}_n \} $ and
${\cal BT}^{\alpha}_{ \! \! \! n} := \{  m_{ \es \et} \,  | \, \es, \et \in \std(\Lambda), \Lambda \in {\cal L}_n(\alpha)\} $
for $ \alpha \in \Par $.
Then $ {\cal BT}$ is a cellular basis for $ \E$ and
$ {\cal BT}^{\alpha}_{ \! \! \! n}$ is a cellular basis for $ \Ea$.
\end{teo}

\begin{demo}
  By the decomposition in (\ref{descE}) it is enough to show that
  $ {\cal BT}^{\alpha}_{ \! \! \! n}$ is a cellular basis for $ \Ea$.
  Let $ \textcolor{black}{\beLambda} $ be the idempotent corresponding to
  any element of $ \Lambda \in {\cal L}_n(\alpha)$: in fact $ \textcolor{black}{\beLambda} $ is
  independent of the choice of $ \Lambda \in {\cal L}_n(\alpha)$. Then
  the set \textcolor{black}{$ \{ g_w \textcolor{black}{\beLambda}   g_{w^1} \, | w, w^{1} \in \Si_n\}   $}
generates $ \Ea $ over $ S $.
Thus letting $ \Lambda = (\blambda \mid \bmu )
\in {\cal L}_n(\alpha) $
vary over pairs of one-column multipartitions with $ \blambda $ of type $ \alpha$ and letting
$ \es, \et $ vary over row standard $ \Lambda$-tableaux,
we get that the corresponding $  m_{\es \et}  $ generate $ \Ea $ over $ S$.
\textcolor{black}{Indeed, for such $ \Lambda$ we have that
$ \blambda  $ is a one-column multipartition and therefore
$ \bT^{\blambda} w$ is row standard
  for any $ w$. Moreover, for such $ \Lambda $
  the row stabilizer of $ \bmu$ is trivial and therefore $ b_{\bmu} $ is just the identity element of $\Si^m_{\Lambda}$.
  In other words, any $ g_w \textcolor{black}{\beLambda}  g_{w^1}$ can be realized in the form
  $m_{\es \et}$ for $ \Lambda$-tableaux $ \es$ and $ \et $.
}

But then, using the last two Lemmas, we deduce that the elements from $ {\cal BT}_{ \! \! \! n}^{\alpha}$
generate $ \Ea$ over $ S$. On the other hand, by the proof of Lemma {\ref{cardinal}}
these elements have cardinality equal to $\dim \Ea$, and so they indeed form a basis for $ \Ea$,
as can be seen by repeating the argument of Theorem \ref{cel}.

\medskip
The $\ast$-condition for cellularity has already been checked above in (\ref{invarianteE}).
Finally, to show the multiplication condition for $ {\cal BT}_{ \! \! \! n}^{\alpha} $ to be cellular, we can repeat the argument from the Yokonuma-Hecke algebra case.
Indeed, to $ \Lambda = (\blambda \mid \bmu) \in {\cal L}_n(\alpha) $ we have associated the $ \Lambda$-tableau $ \et^{\Lambda} $ and have noticed
that the only standard $\Lambda$-tableau $ \et $ satisfying $\et \unrhd \et^{\Lambda} $ is $ \et^{\Lambda} $ itself. The Theorem
follows from this just like in the Yokonuma-Hecke algebra case.
\end{demo}

\begin{coro}{\label{dim}}
The dimension of the cell module $ C(\Lambda) $ associated with $ \Lambda \in {\cal L}_n $ is given by the formula of Corollary \ref{spelling}.
\end{coro}

\textcolor{black}{Unlike $\Si_{\Lambda}^m $, the group $\Si_{\Lambda}^k $ has so far not played any important role in the article,
but now it enters the game. We need the following definition.}

\begin{defi}\label{wreathtypedefinition} Let $ \Lambda \in {\cal L}_n(\alpha) $
\textcolor{black}{for $ \alpha \in \Par$}
and let $ \es=(\Bs \mid \bu) $ be a $ \Lambda$-tableau. Then we say
that $ \es $ (and $ \Bs$) is of wreath type for $ \Lambda $ if $\Bs = \Bs_0 \textcolor{black}{y} $
for some $ \textcolor{black}{B_y} \in \Si_{\Lambda}^k$ where $\Bs_0$ a multitableau of the
initial kind. Moreover we define $$ {\cal BT}^{\alpha, wr}_{ \! \!
\! n} := \{  m_{ \es \et} \,  | \, \es, \et \in \std(\Lambda) \mbox{
of wreath type for }  \Lambda \in {\cal L}_n(\alpha)\}. $$
\end{defi}

The next Corollary should be compared with the results of Geetha and Goodman, \cite{GG}, who show that
$ {\cal A } \wr \Si_m $ is a cellular algebra whenever $ \cal A $ is a {\it cyclic} cellular algebra; \textcolor{black}{by definition this means}
that all cell modules all cyclic.

\begin{coro}\label{wreathcellular}
We have that $ {\cal BT}^{\alpha, wr}_{ \! \! \! n} \! \!  $ is a cellular basis for the subalgebra
$\mathcal{H}^{wr}_{\alpha}(q) $ of $\Ea$,
\textcolor{black}{given}
by Lemma {\ref{embeddingwreath}}.
\end{coro}
\begin{demo}
\textcolor{black}{
Let us first check that $   {\cal BT}^{\alpha, wr}_{ \! \! \! n}  \subseteq \mathcal{H}^{wr}_{\alpha}(q) $.
This is an argument similar to the one used in the beginning of Lemma \ref{Lemmatableaustandar}.
Supposing $ \es=(\Bs \mid \bu) $ and $ \et=(\bT \mid \bv) $ are of wreath type we may use Lemma {\ref{technical}} to write
\begin{equation}
\begin{array}{rl} m_{ \es \et} =&
g_{d(\Bs)}^\ast   \textcolor{black}{\beLambda} \BB_{ d(\bu)}^{\ast} b_{ \bmu}
x_{\blambda}
\BB_{ d(\bv)} g_{d(\bT)} \\
= & g_{y_{\Bs}}
g_{d(\Bs_0)}^\ast   \textcolor{black}{\beLambda} \BB_{ d(\bu)}^{\ast} b_{ \bmu}
x_{\blambda}
\BB_{ d(\bv)} g_{d(\bT_0)} g_{y_{\bT}}
 \\
= & \BB_{y_{\Bs}}
g_{d(\Bs_0)}^\ast   \textcolor{black}{\beLambda} \BB_{ d(\bu)}^{\ast} b_{ \bmu}
x_{\blambda}
\BB_{ d(\bv)} g_{d(\bT_0)} \BB_{y_{\bT}}
\end{array}
\end{equation}
where $ B_{y_{\Bs}}, B_{y_{\bT}}$ belong to $\Si_{\Lambda}^k$ and $ \Bs_0, \bT_0 $ are multitableaux of the initial kind.
Expanding
$ b_{ \bmu} $ out as a linear combination of $ \BB_y $'s with $ B_y \in \Si_{\Lambda}^m $ this becomes
via Lemma {\ref{easytoverify}}
and
Lemma {\ref{importantcommutation}}
a linear combination of
\begin{equation}{\label{thisbecomesvia}}
\BB_{y_{\Bs}}
g_{d(\Bs_0)}^\ast   \textcolor{black}{\beLambda} \BB_{ d(\bu)}^{\ast}  \BB_y
x_{\blambda}
\BB_{ d(\bv)} g_{d(\bT_0)} \BB_{y_{\bT}}
=
\BB_{y_{\Bs}} \BB_{y_1}
  \textcolor{black}{\beLambda}  x_{  \Bs_0    \bT_0    } \BB_{y_2}  \BB_{y_{\bT}}
\end{equation}
where $ B_{y_1}, B_{y_2} \in \Si_{\Lambda}^m$.
Since
$  \Bs_0  $ and $   \bT_0  $ are of the initial kind
we now get
from
Lemma {\ref{embeddingwreath}} that ({\ref{thisbecomesvia}}), and hence also $ m_{ \es \et}$, belongs to $\mathcal{H}^{wr}_{\alpha}(q) $, as claimed.
}

\medskip
Next
it follows from Geetha and Goodman's results in \cite{GG}, \textcolor{black}{or via a direct counting argument}, that the cardinality of $ {\cal BT}^{\alpha, wr}_{ \! \! \! n}
\! \!  $
is equal to the dimension of $ \mathcal{H}^{wr}_{\alpha}(q)$. On the other hand, one easily checks that
Lemma {\ref{secondlast}} holds for $  {\cal BT}^{\alpha, wr}_{ \! \! \!  n } \! \!    $ with respect to $ h \in \mathcal{H}^{wr}_{\alpha}(q) $.
Moreover, applying the
straightening procedure of Lemma \ref{Lemmatableaustandar}
on $ m_{ \es \et}  $ for $ \es, \et $ tableaux of wreath type, the result is a linear combination of $ m_{ \eu \ev}  $
where $ \eu, \ev $ are standard tableaux and still of wreath type. Thus the proof of Theorem {\ref{mainTheorem}}
also gives a proof of the
Corollary.
\end{demo}

\begin{obs}
Recall that \textcolor{black}{we have $  \alpha = (a_r^{k_r}, \ldots,a_1^{k_1}) \in \Par$}
with $\Si^k_{\Lambda} = \Si_{ k_1} \times \cdots\times \Si_{ k_r} $.
From Geetha and Goodman's cellular basis
for $ \mathcal{H}^{wr}_{\alpha}(q)$
one may have expected $ {\cal BT}^{\alpha, wr}_{ \! \! \! n} $ to be slightly different,
namely given by pairs $ (\Bs \mid \bu) $ such that $ \Bs $ is a multitableau of the initial kind whereas $ \bu $
is an $ r$-tuple of multitableaux on the numbers $ \{a_i k_i \} $.
For example
for $\Lambda=\bigg(\big( (1,1),(2), (2),  (2,1)\big)\,\big| \,\big((1), (1,1),(1)   \big)\bigg)$
we would have expected tableaux of the following form
\begin{equation}\begin{array}{c}\et :={\footnotesize\left(\;\left(\;\young(1,2)\;,\young(34)\;,\young(56)\;,\young(79,8)\;\right)  \bigg|
\left( \big(\, \young(2) \, , \young(1,3)\, \right), \left( \young(4) \right)\;\right)}
\end{array}
\end{equation}
where the shapes of the multitableaux occurring in $ \bu$ are given by the equally shaped tableaux of $ \Bs$.
On the other hand, there is an obvious bijection between our standard tableaux of wreath type and
the standard tableaux appearing in Geetha and Goodman's basis and so the cardinality of our basis is correct, which is enough for
the above argument to work.

\end{obs}

\subsection{$\E$ is a direct sum of matrix algebras}
In this final subsection we use the cellular basis for $ \Ea $ to show that $\E$ is isomorphic to a direct sum of matrix algebras in the
spirit of Lusztig and Jacon-Poulain d'Andecy's result for the Yokonuma-Hecke algebra.


Suppose that $\Lambda = (\blambda \mid \bmu)  \in{\cal L}_n(\alpha)$ and that $\es=(\Bs
\mid\bu)$ is a standard $\Lambda$-tableau.
Recall the decomposition $ d(\Bs) = d(\Bs_0)w_{\Bs} $ such
that
$ \Bs_0 $ is a multitableau of the initial kind
and such that
$ \ell(d(\Bs)) = \ell(d(\Bs_0))+\ell(w_{\Bs}) $.
We remark that if
$ \sigma \in \Si_{\alpha} $ permutes the numbers inside the components of $ \Bs $ then
$ w_{\Bs } = w_{\Bs \sigma } $. Indeed, for such $ \sigma $ we have that
$ d({\Bs \sigma })   = \sigma_0d({\Bs  })  $ where
$ \sigma_0 \in \Si_n $ is an element permuting the numbers inside the components of $ \bT^{\lambda} $, that
is $ \bT^{\lambda} \sigma_0  $ is of the initial kind.
But then
$ d( \Bs \sigma )  =(\sigma_0 d(\Bs_0)) w_{\Bs} $ is the decomposition of
$ d( \Bs \sigma )$ and so $ w_{\Bs } = w_{\Bs \sigma } $, as claimed.

We now explain a small variation of this decomposition.
Since $ \Bs  $ is increasing we
have that $ i < j $ if and only if
$ \min ( \s^{(i)}) < \min (\s^{(j)}) $
whenever $ \lambda^{(i)} = \lambda^{(j)} $. We now choose $ B_y \in \Si_{\Lambda}^k $ such that
$ \overline{\Bs }:=\bT^{\blambda} B_y d(\Bs)  $ is increasing in the {\it stronger} sense that $ i < j $ iff
$ \min ( \overline{\s}^{(i)}) < \min ( \overline{\s}^{(j)}) $ whenever
$ |  \lambda^{(i)} | = | \lambda^{(j)} |$. Clearly such a $ B_y $ exists and is unique.
We then consider the decomposition $ d( \overline{\Bs}) = d( \overline{\Bs_0}) w_{ \overline{\Bs}}$.  Since
$ d( \overline{\Bs}) = B_y d( {\Bs}) $
we have
\begin{equation}{\label{decomwreath}}d( {\Bs}) = d( {\Bs}_1) z_{  {\Bs}}
\end{equation}
where $ z_{  {\Bs}} := w_{  \overline{\Bs}} $ and
where $ \Bs_1 := \bT^{\blambda} B_y^{-1} d( \overline{\Bs}_0) = \overline{\Bs}_0 B_y^{-1}$
is a tableau of wreath type.
This gives us the promised decomposition of
$ d( {\Bs}) $. The numbers within the components of $ \bT^{\blambda} z_{  {\Bs}} $ are
just the numbers within the components of $ \overline{\Bs}$ and so
$ \bT^{\blambda} z_{  {\Bs}} $ is
an increasing multitableau in the strong sense defined above.

\begin{lem}\label{lemadecompwreath}
With the above notation we have the following properties.
  \begin{itemize}
\setlength\itemsep{-0.2em}
\item[(1)] The decomposition in ({\ref{decomwreath}}) is unique subject to $ \Bs_1 $ being of wreath type and
  $ \bT^{\blambda} z_{  {\Bs}} $ being increasing in the strong sense.
\item[(2)]   The $ z_{\Bs} $'s appearing in (\ref{decomwreath}) are representatives for distinct coset classes of
  $ \Si_{\Lambda} \backslash \Si_n$ where $ \Si_{\Lambda}$ is the stabilizer group of the set partition
  $ A_{\blambda} $, as introduced above.
\end{itemize}
\end{lem}

\begin{demo}
  Recall that $ \Si_{\Lambda} $ is a product of groups
$  \left(\Si_{a_i} \times \cdots \times \Si_{a_i} \right) \rtimes \Si_{k_i} $.
     Let us prove (2).
    Supposing that $ z_{\Bs} $ and $ z_{\bT}$ belong to the same $ \Si_{\Lambda}$-coset, we have that
    $ z_{\Bs}= \sigma_0 B z_{\bT} $ where the components of $ \sigma_0 $ belong to the
$  \Si_{a_i} \times \cdots \times \Si_{a_i}  $'s and the components of
$ B $ belong to the $ \Si_{k_i} $'s, according to the above description of $ \Si_{\Lambda} $.
    Now both $ \bT^{\blambda} z_{\Bs} $ and $ \bT^{\blambda} z_{\bT} $ are increasing multitableaux in the strong sense
    which implies that $ B = 1 $. Hence
    $ \bT^{\blambda} z_{\Bs} $ and $ \bT^{\blambda} z_{\bT} $ are equal up to a permutation of the numbers inside their
    components and so, by the remark prior to the Lemma, we must have $ z_{\Bs} = z_{\bT}$. This proves
    part (2) of the Lemma. The uniqueness statement of (1) is shown by a similar argument.
\end{demo}

\medskip
For any $\Lambda$-tableau $\es=(\Bs \mid\bu)$ we define the $\Lambda$-tableau $\es_1$ via
\begin{equation}
\es_1=(\Bs_1 \mid\bu).
\end{equation}

\begin{lem}\label{lemamulti}
  Suppose that $\Lambda=\textcolor{black}{(\blambda \mid \bmu)},\,\overline{\Lambda}=
  \textcolor{black}{(\overline{\blambda} \mid\overline{\bmu})}\in{\cal L}_n(\alpha)$, \textcolor{black}{that} $\es=(\Bs \mid\bu)$ is a standard $\Lambda$-tableau
and that $\et=(\bT \mid \bv)$ is a standard $\overline{\Lambda}$-tableau. \textcolor{black}{Then we have that}
\begin{equation}{\label{statement}}
m_{\et^{\Lambda}\es}m_{\et\et^{\overline{\Lambda}}}=\left\{\begin{array}{ll}
m_{\et^{\Lambda}\es_1}m_{\et_1\et^{\overline{\Lambda}}}&\mbox{ if }
{z}_{\Bs} = {z}_{\bT} \\
0&\mbox{ otherwise. } \end{array}\right.
\end{equation}
\end{lem}

\begin{demo}
Both $ \blambda $ and $ \overline{\blambda} $ are of type $ \alpha$ and so
$ \be_{{\Lambda}} = \be_{\overline{\Lambda}} $.
In the decomposition $d( {\Bs}) = d( {\Bs}_1) z_{  {\Bs}}   $ from ({\ref{decomwreath}})
we have in general that $ l(d( {\Bs})) \neq l(d( {\Bs}_1))+l( z_{  {\Bs}} ) $, but even so
$ m_{\et^{\Lambda}\es} = m_{\et^{\Lambda}\es_1} g_{z_{\bs}}$ by Lemma {\ref{technical}}.
Similarly we have that $ m_{\et^{\overline{\Lambda}}\et} = m_{\et^{\overline{ \Lambda}}\et_1} g_{z_{\bT}}$.
Hence we get via Proposition \ref{EGpropiedades}\textcolor{black}{(2)} that
\begin{equation}\label{lemazero1}
\begin{array}{lll}
m_{\et^{\Lambda}\es} m_{\et\et^{\overline{\Lambda}}} &= m_{\et^{\Lambda}\es_1} g_{z_{\Bs}}  g_{z_{\bT}}^{\ast} m_{\et_1 \et^{\Lambda}}                     =
m_{\et^{\Lambda}\es_1} \be_{{\Lambda}} g_{z_{\Bs}}  g_{z_{\bT}}^{\ast} \be_{{\Lambda}}m_{\et_1 \et^{\Lambda}}                  \\  & =
m_{\et^{\Lambda}\es_1}g_{z_{\Bs}}  \be_{   (A_{\blambda})z_{\Bs} }   \be_{(A_{\blambda})z_{\bT} } g_{z_{\bT}}^{\ast} m_{\et_1 \et^{\Lambda}}.
\end{array}
\end{equation}
We now apply the previous Lemma to deduce that $ \be_{   (A_{\blambda})z_{\Bs} }   \be_{(A_{\blambda})z_{\bT} } =0 $  if
$ z_{\Bs} \neq z_{\bT} $ and hence also $ m_{\et^{\Lambda}\es} m_{\et\et^{\overline{\Lambda}}} = 0$
if $ z_{\Bs} \neq z_{\bT} $, thus showing the second part of the Lemma.
Finally, if $ z_{\Bs} = z_{\bT} $ we have that
\begin{equation}\label{lemazero2}
g_{z_{\Bs}}  \be_{   (A_{\blambda})z_{\Bs} }   \be_{(A_{\blambda})z_{\bT} } g_{z_{\bT}}^{\ast} =
\be_{{\Lambda}}    g_{z_{\Bs}}   g_{z_{\bT}}^{\ast}  \be_{{\Lambda}} = \be_{{\Lambda}}
\end{equation}
as can be seen, once again, by expanding $ z_{\Bs} $ out in terms of simple transpositions and noting that
the action at each step involves different blocks. The first part of the Lemma now follows by combining
(\ref{lemazero1}) and (\ref{lemazero2}).
\end{demo}

\medskip
Recall that for
any algebra $ \cal A $ we denote by $ {\rm Mat}_{N}( {\cal A})$ the algebra of
$N\times N$-matrices with entries in $ \cal A$.

\medskip
The cardinality of $ \{  z_{\textcolor{black}{\Bs}} \} $ is
$ b_n(\alpha) $, the
Fa\`{a} di
Bruno number.
We introduce an
arbitrary total order on $ \{  z_{\Bs} \} $
and denote by $M_{\Bs \bT}  $ the elementary matrix of $ {\rm
Mat}_{b_n(\alpha)}\big(\mathcal{H}^{wr}_{\alpha}(q)\big) $ which is
equal to 1 at the intersection of the row and column indexed by $ z_{\Bs}$ and $ z_{\bT}$, and $0$ otherwise.

We can now prove the following
isomorphism Theorem. {\color{black}Via the decomposition in (\ref{descE}) we deduce from it that $ \E$ is isomorphic to a direct sum of matrix algebras, as promised.}

\begin{teo}{\label{promisediso}}
Let $\alpha$ be a partition of $n$.
The $S$-linear map $ \Psi_\alpha $ given by
$$ \Ea \longrightarrow{\rm Mat}_{b_n(\alpha)}\big(\mathcal{H}^{wr}_{\alpha}(q) \big), \, \, \, m_{ \es \et}
\mapsto m_{\es_1\et_1}M_{\Bs\bT}
$$
is an isomorphism of $S$-algebras.
A similar statement holds for the specialized algebra over $ \kk $.
\end{teo}

\begin{demo}
  Note first that by Corollary \ref{wreathcellular} we have that
  $ m_{\es_1\et_1} \in \mathcal{H}^{wr}_{\alpha}(q) $. Furthermore,
  by the uniqueness statement of the previous Lemma we have that
$\Psi_{\alpha}$ maps an $S$-basis to an $S$-basis and so
we only need to show that it is a homomorphism, preserving the multiplications on both sides.

For this suppose that $\Lambda, \overline{\Lambda} \in{\cal L}_n(\alpha) $.
Given a pair of standard $\Lambda$-tableaux $\es=(\Bs \mid\bu_1)$,  $\et=(\bT \mid \bu_2)$
and a pair of standard $\overline{\Lambda}$-tableaux  $\eu=(\Bu \mid\bv_1)$, $\ev=(\Bv \mid \bv_2)$ we have by the
previous Lemma that
\begin{equation}
\begin{array}{rl}
  m_{\es\et}m_{\eu\ev} & = \left\{\begin{array}{ll}m_{\es\et_1}m_{\eu_1\ev} &
\, \, \, \,\, \, \, \,\, \, \, \,\, \, \, \, \, \, \, \,\, \, \, \,\, \, \, \, \, \, \, \,\, \, \, \,\, \, \, \mbox{ if } z_{\bT}={z}_{\U}\\
0& \, \, \, \,\, \, \, \,\, \, \, \,\, \, \, \, \, \, \, \,\, \, \, \,\, \, \, \, \, \, \, \,\, \, \, \,\, \,\,  \mbox{  otherwise }  \end{array}\right. \\
& = \left\{\begin{array}{ll} g_{z_{\Bs}}^{\ast} ( m_{\es_1\et_1}m_{\eu_1\ev_1})  g_{z_{\Bv}}  &
\mbox{ if } {z}_{\bT}={z}_{\Bu}\\
0&\mbox{ otherwise. }  \end{array}\right.
  \end{array}
\end{equation}
Expanding $ m_{\es_1\et_1}m_{\eu_1\ev_1} \in \mathcal{H}^{wr}_{\alpha}(q) $ out as a linear
combination of cellular basis elements $ m_{\ea_1\eb_1} $ of $ \mathcal{H}^{wr}_{\alpha}(q) $ we have that
$ m_{\es\et}m_{\eu\ev}  $ is the corresponding linear combination of
$ g_{z_{\Bs}}^{\ast} m_{\ea_1\eb_1} g_{z_{\Bv}}\! $'s, and so
\begin{equation}\label{compare1}
  \Psi_{\alpha}(m_{\es\et}m_{\eu\ev} )  = \left\{
  \begin{array}{ll}
      m_{\es_1\et_1}m_{\eu_1\ev_1} M_{\Bs \Bv } &
    \mbox{ if }
{z}_{\bT}={z}_{\Bu}\\
0&\mbox{ otherwise. }  \end{array}\right.
\end{equation}
On the other hand, by the matrix product formula $M_{\Bs \bT}M_{\Bu \Bv }=\delta_{z_{\bT} z_{ \Bu} }M_{\Bs \Bv }$
we have that
\begin{equation}\label{compare2}
\Psi_{\alpha}(m_{\es\et})\Psi_{\alpha}(m_{\eu\ev})=\left\{\begin{array}{ll}m_{\es_1\et_1}m_{\eu_1\ev_1}M_{\Bs\Bv}&\mbox{ if }
{z}_{\bT}={z}_{\Bu}\\
0&\mbox{ otherwise. }  \end{array}\right.
\end{equation}
Comparing (\ref{compare1}) and (\ref{compare2}) we conclude that $ \Psi_{\alpha} $ is an algebra homomorphism as
claimed. The Theorem is proved.
\end{demo}

\color{black}
\newpage

{\small

Throughout the paper we adopt the following conventions:
\begin{itemize}
\item We use the normal frak font, like $\s$, to denote tableaux whose shape is a composition.
\item We use the boldfrak font, like $\Bs$, to denote multitableaux whose shape is a multicomposition.
\item We use the mathematical doble-struck font, like $\es$, to denote tableaux whose shape is an element of
  $ \mathcal{L}_n$.
\item For $\lambda$ a composition (resp. $\blambda $ a multicomposition,
  resp. $ \Lambda $ an element of $ \mathcal{L}_n $) we denote
  by $ \T^\lambda $ (resp. $\bT^{\blambda} $, resp. $\et^{\Lambda}$)
  the maximal tableau of shape $ \lambda$ (resp. shape $\blambda$, resp. shape $\Lambda$)
  as introduced in the text. Note that $\bT^{\blambda} $ and $\et^{\Lambda}$ are not the unique maximal tableaux
  of their shape.
\item
For $\lambda$ a composition (resp. $\blambda $ a multicomposition,
  resp. $ \Lambda $ an element of $ \mathcal{L}_n $) we denote
  by $ \std(\lambda) $ (resp. $ \std(\blambda) $, resp. $\std(\Lambda)$)
  the set of standard $ \lambda$-tableaux (resp. $\blambda$-tableaux, resp. $ \Lambda$-tableaux).
\end{itemize}

\printnomenclature[1.5cm]
}

\sc Instituto de Matem\'atica y F\'isica, Universidad de Talca, Chile, mail joespinoza@utalca.cl,
steen@inst-mat.utalca.cl,

\end{document}